\documentclass{article}
\usepackage[utf8]{inputenc}
\usepackage{amssymb,mathtools}
\usepackage{amsthm}
\usepackage[hidelinks]{hyperref}
\usepackage{xr}
\usepackage{cleveref}
\usepackage{comment}
\usepackage{graphicx}
\usepackage{makecell}
\usepackage{float}
\usepackage[cal=cm]{mathalfa}
\usepackage{algorithm}
\usepackage[normalem]{ulem}
\usepackage{algpseudocode}
\usepackage{pgfplots,pgfplotstable}
\usepgfplotslibrary{external}
\usepackage{booktabs,multirow}
\usepackage[nomessages]{fp}
\usepgfplotslibrary{groupplots}
\pgfplotsset{
    every axis/.append style={
        line width=1pt,
        label style={font=\footnotesize},
        tick label style={font=\footnotesize},
        legend style={font=\scriptsize},
    },
    compat=1.18
}

\allowdisplaybreaks
\usepackage[
    top=2.5cm,
    left=3cm,
    right=3cm,
    bottom=4.0cm
]{geometry}
\usepackage{silence}

\usepackage[table]{xcolor}

\usepackage{listings}
\usepackage{color}
\usepackage{titlesec}

\titleformat{\paragraph}[runin]
  {\normalfont\normalsize\bfseries}
  {}
  {0pt}
  {}
  [.]

\numberwithin{equation}{section}

\definecolor{dkgreen}{rgb}{0,0.6,0}
\definecolor{gray}{rgb}{0.5,0.5,0.5}
\definecolor{mauve}{rgb}{0.58,0,0.82}

\newenvironment{nalign}{
    \begin{equation}
    \begin{aligned}
}{
    \end{aligned}
    \end{equation}
    \ignorespacesafterend
}

\newcommand{\bb}[1]{\mathbf{#1}}
\newcommand{\bd}[1]{\dot{\mathbf{#1}}}

\newcommand{\R}[0]{\mathbb{R}}

\newcommand{\pb}[0]{\mathbb{P}}
\newcommand{\p}[0]{\mathcal{P}}

\newcommand{\JM}{\widehat{J}}
\newcommand{\YM}{\widehat{Y}}
\newcommand{\J}{J}
\newcommand{\Jr}{J_r}
\newcommand{\Jop}{J_{\mathrm{op}}}
\newcommand{\Jophr}{J_{\mathrm{op,hr}}}
\newcommand{\barJ}{\overline{J}}
\newcommand{\barJhr}{\overline{J}_{\mathrm{hr}}}
\newcommand{\barY}{\overline{Y}}
\newcommand{\barYhr}{\overline{Y}_{\mathrm{hr}}}
\newcommand{\Yophr}{Y_{\mathrm{op,hr}}}
\newcommand{\GM}{\widehat{G}}
\newcommand{\G}{G}
\newcommand{\Gr}{G_r}
\newcommand{\Gop}{G_\mathrm{op}}
\newcommand{\Gophr}{G_\mathrm{op,hr}}
\newcommand{\barG}{\overline{G}}
\newcommand{\barGhr}{\overline{G}_{\mathrm{hr}}}

\newcommand{\Xophr}{V_\mathrm{op,hr}}
\newcommand{\QX}{\mathcal{Q}_X}
\newcommand{\QS}{\mathcal{Q}_S}
\newcommand{\QShr}{\mathcal{Q}_{\Seim}}
\newcommand{\QH}{\mathcal{Q}_H}
\newcommand{\QHhr}{\mathcal{Q}_{\Heim}}
\newcommand{\ccS}{\bb{c}_{S}}
\newcommand{\ccH}{\bb{c}_{H}}
\newcommand{\ccX}{\bb{c}_{X}}
\newcommand{\ccHhr}{\bb c_{H_{\eim}}}
\newcommand{\cS}{\ccS(\mu)}
\newcommand{\cH}{\ccH(\mu)}
\newcommand{\cX}{\ccX(\mu)}
\newcommand{\cShr}{\bb{c}_{\Seim}(\mu)}
\newcommand{\cHhr}{\bb{c}_{\Heim}(\mu)}

\newcommand{\bS}{\kappa_{S}}
\newcommand{\bH}{\kappa_{H}}
\newcommand{\bX}{\kappa_{X}}
\newcommand{\bShr}{\kappa_{\Seim}}
\newcommand{\bHhr}{\kappa_{\Heim}}

\newcommand{\PS}{\Pi_{S,\ccS}}
\newcommand{\PX}{\Pi_{X,\ccX}}
\newcommand{\PHhr}{\Pi_{\Heim,\bb{c}_{\Heim}}}
\newcommand{\PShr}{\Pi_{\Seim,\bb{c}_{\Seim}}}
\newcommand{\uf}{\bb{u}}
\newcommand{\ur}{\bb{u}_r}
\newcommand{\dur}{\dot{\bb{u}}_r}
\newcommand{\uhr}{\bb{u}_{\hr}}
\newcommand{\duhr}{\dot{\bb{u}}_{\hr}}
\newcommand{\jf}{\bb{j}}
\newcommand{\yr}{\bb{j}_{r}}

\newcommand{\jeim}{\bb{j}_{\eim}}
\newcommand{\yf}{\bb{y}}
\newcommand{\yeim}{\bb{y}_{\eim}}
\newcommand{\gf}{\bb{g}}
\newcommand{\zr}{\bb{g}_{r}}

\newcommand{\geim}{\bb{g}_{\eim}}
\newcommand{\bw}{\bb{w}}
\newcommand{\bv}{\bb{v}}

\newcommand{\norm}[1]{\lVert#1\rVert}
\newcommand{\inp}[1]{\left\langle #1\right\rangle}
\newcommand{\SM}{\mathbb{S}} 
\newcommand{\MtoM}{M} 

\newcommand{\wrec}{\bw^{(r)}}

\newcommand{\malpha}{\bb{m}^{\alpha}}

\newcommand{\ahat}{\widehat{\bb{a}}^\alpha}

\newcommand{\Ahat}{\widehat{A}^\alpha}
\newcommand{\aj}{\alpha_J}
\newcommand{\ag}{\alpha_G}
\newcommand{\ay}{\alpha_Y}
\newcommand{\kS}{\beta_{S}}
\newcommand{\kH}{\beta_{H}}
\newcommand{\Cu}{C_{u}}

\newcommand{\LX}{\mathcal{L}_X}
\newcommand{\NX}{\mathcal{N}_X}
\newcommand{\LH}{\mathcal{L}_H}

\newcommand{\NH}{\mathcal{N}_H}

\newcommand{\NHhr}{\mathcal{N}_{H,\hr}}

\newcommand{\NXhr}{\mathcal{N}_{X,\hr}}
\newcommand{\Dx}{D_{X}}
\newcommand{\Dh}{D_{H}}

\newcommand{\cNh}{\bb{a}_H}
\newcommand{\cNx}{\bb{a}_X}
\newcommand{\tcNx}{\widetilde{\bb{a}}_X}

\newcommand{\Heim}{H_{\eim}}
\newcommand{\Seim}{S_{\eim}}
\newcommand{\Xeim}{X_{\eim}}
\newcommand{\nJmeim}{m_J}

\newcommand{\nYmeim}{m_Y}
\newcommand{\nHeim}{m_H}
\newcommand{\nSeim}{m_S}
\newcommand{\nXeim}{m_X}
\newcommand{\idx}{\tau}
\newcommand{\errmu}{\mathcal{E}_\Lambda^\mu}
\newcommand{\errmuF}{\mathcal{E}_{\Lambda,X}^\mu}

\newcommand{\errmutest}{\mathcal{E}_\Lambda^{\mu_\text{test}}}

\newcommand{\errmutestH}{\mathcal{E}_{\Lambda,H}^{\mu_\text{test}}}
\newcommand{\errmutestS}{\mathcal{E}_{\Lambda,S}^{\mu_\text{test}}}

\newcommand{\Dmu}{\Omega_{\mu}}
\newcommand{\DmuH}{\Omega_{H,\ccH}}
\newcommand{\DmuS}{\Omega_{S,\ccS}}
\newcommand{\DmuX}{\Omega_{X,\ccX}}
\newcommand{\Dmur}{\Omega_{\mu,r}}
\newcommand{\DmuSr}{\Omega_{S,r}}
\newcommand{\Dmueim}{\Omega_{\mu}^{\eim}}
\newcommand{\DmuHeim}{\Omega_{\Heim,\bb{c}_{\Heim}}}
\newcommand{\DmuSeim}{\Omega_{\Seim,\bb{c}_{\Seim}}}
\newcommand{\Dmuhr}{\Omega_{\mu,\hr}^{\eim}}

\newcommand{\Dmuh}{\mathcal{D}_{\Lambda,H}^\mu}
\newcommand{\DmutestH}{\mathcal{D}_{\Lambda,H}^{\mu_\text{test}}}

\DeclareMathOperator{\rank}{rank}

\DeclareMathOperator{\range}{Range}
\DeclareMathOperator{\hr}{hr}
\DeclareMathOperator{\vect}{vec}
\DeclareMathOperator{\eim}{EIM}
\DeclareMathOperator{\diag}{diag}
\DeclareMathOperator{\train}{train}

\newtheorem{theorem}{Theorem}[section]
\newtheorem{lemma}[theorem]{Lemma}
\newtheorem{proposition}[theorem]{Proposition}

\theoremstyle{definition}

\theoremstyle{remark}
\newtheorem{remark}[theorem]{Remark}

\AddToHook{env/theorem/begin}{\crefalias{section}{theorem}}
\AddToHook{env/lemma/begin}{\crefalias{section}{lemma}}
\AddToHook{env/proposition/begin}{\crefalias{section}{proposition}}
\AddToHook{env/corollary/begin}{\crefalias{section}{corollary}}
\AddToHook{env/definition/begin}{\crefalias{section}{definition}}
\AddToHook{env/remark/begin}{\crefalias{section}{remark}}
\AddToHook{env/ex/begin}{\crefalias{section}{ex}}

\title{Structure-Preserving Model Reduction of Metriplectic Systems via Oblique Projections and Block Structure}

\author{
  A. Benvenuti\thanks{Department of Mathematics, University of Pisa, Pisa, Italy (\url{andrea.benvenuti@phd.unipi.it})}
  \and
  C. Pagliantini\thanks{Department of Mathematics, University of Pisa, Pisa, Italy (\url{cecilia.pagliantini@unipi.it})}
}
\date{}

\begin{document}

\maketitle
\begin{abstract}
    Metriplectic systems combine conservative and dissipative dynamics and are characterized by energy conservation and entropy production. The degeneracy conditions underlying these properties are generally not preserved by standard projection-based model order reduction.
    We develop two structure-preserving reduction strategies. The first enforces degeneracy conditions through oblique projections and combines gradient-preserving and matrix empirical interpolation for hyper-reduction.
    The second constructs reduced bases adapted to available operator block structure, preserving degeneracy without online oblique corrections. We establish sufficient conditions for preserving the Jacobi identity and derive a priori error estimates. Numerical tests on benchmark problems demonstrate accuracy, thermodynamic consistency, and computational efficiency.
\end{abstract}


\section{Introduction}

Many dynamical systems combine reversible motion with irreversible processes, such as friction and heat transfer
\cite{GO97,grmela1997dynamics,suzuki20,DPZ12,EGB20}. In a thermodynamically consistent description of a closed system, these processes conserve total energy while producing entropy. The metriplectic formalism~\cite{morrison1986paradigm} and the general equation for the nonequilibrium reversible-irreversible coupling (GENERIC)~\cite{grmela1997dynamics} encode this behavior through the coupling of a Hamiltonian vector field and a generalized gradient flow. In finite dimensions and in a parametric setting, the dynamics is modeled by the evolution
\begin{equation}\label{eq:pbm}
    \dot{\uf}(t;\mu)=\J (\uf(t;\mu))\nabla H(\uf(t;\mu);\mu)+\G(\uf(t;\mu))\nabla S(\uf(t;\mu);\mu),\quad \forall\,t\in\mathcal{T},\, \mu\in\p,
\end{equation}
with initial condition $\uf(0;\mu)=\uf_0(\mu), $ where $\mathcal{T}:=[0,T]$ is a temporal interval and $\p\subset\R^p$ is a compact set of parameters. Here $\uf:\mathcal{T}\times\p\to\R^N$ denotes the state variable, $\J(\uf):\R^N\to\R^{N\times N}$ is a Poisson operator, and $\G(\uf):\R^N\to\R^{N\times N}$ is symmetric positive semi-definite. The functionals $H, S:\R^N\times\p\to\R$ 
represent the Hamiltonian and entropy, respectively, although the physical interpretation of $S$ may vary between applications.
%
A metriplectic system is also characterized by the so-called \emph{mutual degeneracy conditions}
\begin{equation}\label{eq:degCond}
    \G(\bv)\nabla H(\bv;\mu)=0,\qquad \J(\bv)\nabla S(\bv;\mu)=0,\qquad\forall\,\bv\in\R^N,\,\mu\in\p.
\end{equation}
These, together with the symmetry properties of the operators,
imply that the Hamiltonian is conserved 
and the entropy is non-decreasing.
Retaining this compatibility is essential when
constructing thermodynamically consistent approximations.


In a general setting, the operators $\J$ and $\G$ in \eqref{eq:pbm} may have the form
$\J(\bv)= M^{-1}\JM(\bv)M^{-1}$ and
$\G(\bv)= M^{-1}\GM(\bv)M^{-1}$, for any $\bv\in\R^N$,
where $M$ is a symmetric positive definite matrix (typically a mass matrix), and for some suitable $\JM$ and $\GM$. 
This is the case, for example, when \eqref{eq:pbm} stems from the spatial discretization of a partial differential equation characterized by a Hamiltonian term and a gradient flow dissipation, as in \cite{lombardi2024conformal}.
Also in this case, we assume that the operator $J$ is Poisson, that is, it is skew-symmetric and satisfies the Leibniz rule and the Jacobi identity.
%
Throughout, we use the term almost-Poisson when the Jacobi identity is not satisfied, and almost-metriplectic for a system
with $J$ almost-Poisson and satisfying the degeneracy conditions \eqref{eq:degCond}.


\paragraph{Model order reduction of Hamiltonian systems with gradient flow dissipation}

Repeated simulation of high-dimensional metriplectic systems can become computationally prohibitive in parameter studies, optimization, and uncertainty quantification. Model order reduction addresses this issue by approximating the solution trajectory in a space of dimension $n\ll N$, typically constructed from instances of the solution at certain values of time and parameter. However, state accuracy alone does not ensure thermodynamic consistency: projection-based reduction need not retain energy conservation or
entropy production.

For problems of the form \eqref{eq:pbm} characterized by a Hamiltonian part and a dissipative part, a
possible strategy to preserve the structure of the model is to approximate the Hamiltonian and entropy gradients in the reduced basis space. This corresponds to performing model order reduction on a mixed formulation of \eqref{eq:pbm} where auxiliary variables are introduced so that it reads:
given $\mu\in\p$, and $\uf(0;\mu)=\uf_0(\mu)$, find $\big(\uf(\cdot;\mu),\jf(\cdot;\mu),\gf(\cdot;\mu)\big)\in C^1(\mathcal{T};\R^N)\times C^0(\mathcal{T};\R^N)\times C^0(\mathcal{T};\R^N)$ such that
\begin{nalign}\label{eq:MF}
    M\bd{u}(t;\mu) &= \JM(\uf(t;\mu))\jf(t;\mu) + \GM(\uf(t;\mu))\gf(t;\mu), && \quad\forall\,t \in \mathcal{T}, \\
    M\jf(t;\mu) &= \nabla H(\uf(t;\mu);\mu), && \quad\forall\,t \in \mathcal{T}, \\
    M\gf(t;\mu) &= \nabla S(\uf(t;\mu);\mu), && \quad\forall\,t \in \mathcal{T}.
\end{nalign}

Let $U\in\R^{N\times n}$, with $n\leq N$, be a reduced basis matrix such that $U^\top MU=I_n$, let $P:\R^N\to\R^n$ be the surjective map $P:=U^{\top} M$, and $\Pi_U:=UP$ the $M$-orthogonal projection onto $\range(U)$. 
Galerkin projection of the mixed formulation \eqref{eq:MF}
onto the $\range(U)\subseteq\R^N$ results in the following problem:
given $\mu\in\p$, and $\ur(0;\mu)=P\uf_0(\mu)$, find $\big(\ur(\cdot;\mu),\yr(\cdot;\mu),\zr(\cdot;\mu)\big)\in C^1(\mathcal{T};\R^n)\times C^0(\mathcal{T};\R^n)\times C^0(\mathcal{T};\R^n)$ such that
\begin{nalign}\label{eq:rMF}
    \bd{u}_r(t;\mu) &= \Jr(\ur(t;\mu))\yr(t;\mu) + G_r(\ur(t;\mu))\zr(t;\mu), &&\quad \forall\,t \in \mathcal{T},\\
    \yr(t;\mu) &= \nabla H_r(\ur(t;\mu);\mu), && \quad\forall\,t \in \mathcal{T},\\
    \zr(t;\mu) &= \nabla S_r(\ur(t;\mu);\mu), && \quad\forall\,t \in \mathcal{T},
\end{nalign}
where, for any $\bw \in \mathbb{R}^n$, $H_r(\bw;\mu) \coloneq H(U\bw;\mu)$ and $S_r(\bw;\mu) \coloneq S(U\bw;\mu)$ are the reduced Hamiltonian and entropy, respectively, and
\begin{equation}\label{eq:JGr}
\Jr(\bw) \coloneq U^\top  \JM(U\bw) U\in\R^{n\times n},\qquad G_r(\bw) \coloneq U^\top  \GM(U\bw) U\in\R^{n\times n}.
\end{equation}

It can be easily verified that the reduced order model \eqref{eq:rMF} preserves the geometric structure of full order model \eqref{eq:MF} in the following sense.
\begin{proposition}\label[proposition]{prop:rMFproperties}
For any $\mu\in\p$ the reduced order model \eqref{eq:rMF} is such that the reduced operators $\Jr$ and $\Gr$ retain skew-symmetry and positive semi-definiteness, respectively. In particular, if $J$ is constant, $\Jr$ is Poisson.
Moreover, $H_r$ is conserved when $G\equiv0$, and $S_r$ is non-decreasing when $J\equiv0$. In the Hamiltonian case, the difference between the Hamiltonian evaluated at the full and reduced solutions is constant and determined by the initial condition.
\end{proposition}

Note that a direct Galerkin projection of \eqref{eq:pbm} does not generally retain the gradient representation \eqref{eq:rMF}, and hence need not inherit the properties above.

In this framework, the reduced basis $U$ should be constructed from snapshots of $(\uf,\jf,\gf)$, as motivated by the error analysis in \Cref{sec:convergence}. This means that the training set also includes instances of the Hamiltonian and entropy gradient.


    We also point out that, owing to \Cref{prop:rMFproperties},
 if the initial condition $\uf(0;\mu)=\uf_0(\mu)$ belongs to $\range(U)$, then the Hamiltonian along reduced trajectories coincides with the full order Hamiltonian.
    If the initial condition of the full order model depends on the parameter $\mu$, enforcing this condition may require a large reduced basis space.
For parameter-dependent initial conditions, the ansatz $\uf_0(\mu)+U\ur(t;\mu)$, with $\ur(0;\mu)=0$, reproduces the initial state of the full order model exactly \cite{hesthaven2021structure}.

\paragraph{Literature review}
Structure-preserving projection-based methods have been extensively studied for Hamiltonian systems \cite{peng2016,afkham2017structure,gong2017structure,buchfink2019symplectic,hesthaven2021structure} and port-Hamiltonian problems, e.g., \cite{gugercin12,polyuga12}.
 
For coupled reversible and irreversible dynamics,
preserving skew-symmetry and positive semi-definiteness under projection is insufficient to preserve the degeneracy conditions: the projected gradients need not belong to the kernels of the projected operators. Therefore, even a reduced model that retains a formal separation into conservative and dissipative terms, as \eqref{eq:rMF}, can fail to conserve energy or to produce entropy. 
For this reason, the reduced order model should satisfy a reduced version of the degeneracy conditions, that is
\begin{equation}\label{eq:rdegCond}
    \Jr(\bw)\nabla S_r(\bw;\mu)=0,\qquad G_r(\bw)\nabla H_r(\bw;\mu)=0,\qquad\forall\, \bw\in\R^n,\,\mu\in\p.
\end{equation}

Öttinger \cite{Ottinger2015} studied the preservation of thermodynamic structure under model reduction, constructing reduced GENERIC brackets using an invariant manifold.
Within projection-based model reduction, a POD-based approach that enforces the degeneracy conditions through exterior-algebraic representations of the reversible and dissipative operators was proposed in \cite{gruber2023energetically}. This method preserves the algebraic properties underlying energy conservation and entropy production. However, evaluating the reduced vector field may still incur a computational cost that depends on the full order dimension.

Related constructions have also been explored in data-driven models, by penalizing degeneracy violation in the training objective function \cite{Hernandez2021} or by enforcing the degeneracy conditions algebraically through operator parameterizations \cite{Zhang2022,gruber2025efficiently}.
These are based on projections onto the orthogonal complements of the energy and entropy gradients and thus encode the required kernel relations. This perspective connects structure-preserving model reduction to the parametrization of the reversible and dissipative operators, and motivates the oblique-projection framework developed in this work.


\paragraph{Main contributions}
Other than structure-preservation, a further difficulty arises in the presence of nonlinear functionals and state-dependent operators, whose evaluation 
may hinder the computational advantage of model order reduction.
Hyper-reduction techniques, such as the Empirical Interpolation Method (EIM)~\cite{barrault2004empirical}, reduces this cost but need not preserve gradient structures or degeneracy.

We develop two structure-preserving model order reduction strategies for parametric metriplectic systems. In the first one, starting from the mixed formulation in which the energy and entropy gradients are represented by auxiliary variables, we introduce oblique projections to enforce the reduced degeneracy conditions while retaining the symmetry properties of the operators.
%
This construction is combined with gradient-preserving EIM \cite{pagliantini2023gradient} for the Hamiltonian and entropy and a matrix form of EIM \cite{negri2015efficient} for the state-dependent operators. Positive semi-definiteness is enforced through a correction performed entirely at the reduced level. The hyper-reduced operators are then employed in the construction of the oblique projections. The resulting hyper-reduced model retains the degeneracy conditions associated with the hyper-reduced functionals and, under suitable sparsity of the nonlinear terms, has arithmetic complexity independent of the full order dimension. We also establish a priori error estimates where the contributions of state and gradient projection from those of hyper-reduction are separated.

The second strategy applies to systems with a suitable block structure. It relies on a block-diagonal reduced basis whose projector preserves the operators' kernel, ensuring the degeneracy conditions without online oblique corrections. We provide a constructive procedure for state-independent kernels and establish approximation bounds in terms of snapshot singular values.

Both strategies preserve the thermodynamic structure and identities, while preservation of the Jacobi identity is guaranteed if the reduced conservative operator $\Jr$ is constant. Moreover, after hyper-reduction, exact Hamiltonian conservation and entropy production are guaranteed for the hyper-reduced functionals. 


The remainder of the paper is organized as follows. \Cref{sec:red-OP} introduces the oblique-projection reduced model and discusses its properties. \Cref{sec:OP-HROM} pertains to the hyper-reduced model based on oblique projections constructed on hyper-reduced operators. The properties of the model, its computational complexity, and convergence estimate are established. \Cref{sec:block} discusses the block-diagonal approach.
In \Cref{sec:num-exp} two numerical tests are shown. Concluding remarks are drawn in \Cref{sec:conclusions}.

\section{Model order reduction via oblique projections}\label{sec:red-OP}

Let $\mu\in\p$ be fixed and let $X\in\{H,S\}$. Given a $\mu$-dependent vector $\cX\in\R^N$, define
\begin{equation}\label{eq:beta}
\bX(\bv;\mu):=\cX^\top \nabla X(\bv;\mu)\qquad
\forall\,\bv\in \R^N,
\end{equation}
and the open set
\begin{equation}\label{eq:DmuX}
\DmuX\coloneq\{\bv\in\R^N \mid\ \bX(\bv;\mu)\neq 0\}\subseteq\R^N.
\end{equation}
Moreover, define, for any $\bv\in\DmuX$ and $\mu\in\p$, the operators
\begin{equation}\label{eq:PinQ}
    \QX(\bv;\mu)\coloneq I_n -\PX(\bv;\mu),\quad
    \PX(\bv;\mu) := U^{\top} \frac{M\cX \nabla X(\bv;\mu)^\top}{\bX(\bv;\mu)}U.
\end{equation}

\begin{lemma}\label[lemma]{lem:oblproj}
    Let $\mu\in\p$ and $\bv\in \DmuX$ be fixed.
    If $\cX\in\range(U)$, then
    \begin{itemize}
    \item[(i)] $\PX(\bv;\mu)\circ \PX(\bv;\mu)=\PX(\bv;\mu)$;
    \item[(ii)] $\PX(\bv;\mu)P\cX=P\cX$;
    \item[(iii)] $\PX(\bv;\mu)^{\top}U^\top\nabla X(\bv;\mu)=U^\top \nabla X(\bv;\mu)$.
    %
    \end{itemize}
In particular, $\PX(\bv;\mu)$ is an orthogonal projection with respect to the Euclidean inner product if $P\cX$ and $U^\top\nabla X(\bv;\mu)$ are parallel, and it is oblique otherwise.
\end{lemma}
\begin{proof}
Let $\bb{\omega}:=P\cX\in\R^n$ and $\bb{\sigma}:=U^\top\nabla X(\bv;\mu)\in\R^n$.
Since $\cX\in\range(U)$, it holds $U\bb{\omega}=\cX$ and $\bb{\sigma}^\top\bb{\omega}=\bX(\bv;\mu)$.
Then, the properties follow by writing $\PX$ as
$$\PX(\bv;\mu)=\dfrac{\bb{\omega}\, \bb{\sigma}^\top}{\bb{\sigma}^\top\bb{\omega}}.$$
\end{proof}

Let $\Dmu\coloneq\DmuH\cap\DmuS$, and
$\Dmur\coloneq \{\bw\in\R^n\mid\, U\bw\in \Dmu\}\subseteq\R^n.$
The oblique-projection reduced order model (OP-ROM) reads: given $\mu\in\p$ and $\ur(0;\mu)=P\bb{u}_0(\mu)\in \Dmur$, find $\ur(\cdot;\mu)\in C^1(\mathcal{T};\R^n)$ such that, for all $t\in\mathcal{T}$,
\begin{equation}\label{eq:OP-ROM}
    \dur(t;\mu)=\Jop(U\ur(t;\mu);\mu)\nabla H_r(\ur(t;\mu);\mu)+\Gop(U\ur(t;\mu);\mu)\nabla S_r(\ur(t;\mu);\mu),
\end{equation}
where, for all $\bv\in \Dmu$,
\begin{equation}\label{eq:JGop}
    \Jop(\bv;\mu)\coloneq\QS(\bv;\mu)\barJ(\bv)\QS(\bv;\mu)^\top,\quad \Gop(\bv;\mu)\coloneq \QH(\bv;\mu)\barG(\bv)\QH(\bv;\mu)^\top,
\end{equation}
the operators $\QS$ and $\QH$ are defined as in \eqref{eq:PinQ},
and
\begin{equation}\label{eq:JGbar}
        \barJ(\bv)\coloneq U^\top\JM(\bv)U,\qquad\barG(\bv)\coloneq U^\top\GM(\bv)U,\qquad\forall\,\bv\in\R^N.
\end{equation}
Observe that $\barJ(U\bw)=\Jr(\bw)$ for any $\bw\in\R^n$.

\begin{remark}[Well-posedness of the reduced order model]\label[remark]{rmk:OPROMwellposed}
Let $\mu\in\p$ be fixed. For sufficiently regular $\JM$, $\GM$, $\nabla H(\cdot;\mu)$, and $\nabla S(\cdot;\mu)$, the reduced velocity field of \eqref{eq:OP-ROM} is locally Lipschitz continuous on $\Dmur$ since $\bH(\cdot;\mu)$ and $\bS(\cdot;\mu)$ are continuous and nonzero in $\Dmu$. Hence, there exists a time $T_{\mu}>0$ such that every initial condition $\ur(0;\mu)=P\uf_0(\mu)$ in $\Dmur$ determines a unique solution $t\in[0,T_{\mu})\mapsto \ur(t;\mu)\in\Dmur$.
In the remaining of the work, we assume that $T<T_{\mu}$ for every $\mu\in\p$, so that the solution exists and belongs to $\Dmur$ for all $t\in\mathcal{T}=[0,T]$.
%
\end{remark}

The introduction of the operators $\QH$ and $\QS$ is instrumental to ensure that the reduced order model \eqref{eq:OP-ROM} is (almost-)metriplectic, as shown in the next result.

\begin{lemma}\label[lemma]{lem:OPmetriplectic}
Let $\mu\in\p$ be fixed and assume that the reduced order model \eqref{eq:OP-ROM} is well-posed in $\Dmur$ in the sense of \Cref{rmk:OPROMwellposed}.
If $\cX\in\range(U)$, for $X\in\{H,S\}$, then
problem \eqref{eq:OP-ROM} is almost-metriplectic.
The Hamiltonian $H$ evaluated at the reduced trajectories is conserved and the entropy is non-decreasing.

Moreover, if $\barJ\in\R^{n\times n}$ defined in \eqref{eq:JGbar} is constant in $\R^N$, then $\Jop(\cdot;\mu)$ is a Poisson operator in $\DmuS$, and thus \eqref{eq:OP-ROM} is metriplectic. 
\begin{proof}
The fact that $\Jop(\cdot;\mu)$ is skew-symmetric and $\Gop(\cdot;\mu)$ is symmetric positive semi-definite follows from their definition \eqref{eq:JGop}  and from the properties of $\widehat{J}$ and $\widehat{G}$.
Therefore, problem \eqref{eq:OP-ROM} is of the form \eqref{eq:pbm}, that is, it is a Hamiltonian system with gradient flow dissipation.
Moreover, by \Cref{lem:oblproj} item (iii), 
$\QX(U\bw;\mu)^\top \nabla X_r(\bw;\mu)=0$, with $X\in\{H,S\}$
and for all $\bw\in\Dmur$; this implies that
the degeneracy conditions corresponding to \eqref{eq:degCond} are satisfied, that is
\begin{equation*}
    \Gop(U\bw;\mu)\nabla H_r(\bw;\mu)=0,\qquad \Jop(U\bw;\mu)\nabla S_r(\bw;\mu)=0,\quad\forall\,\bw\in\Dmur.
\end{equation*}

Let $\DmuSr:=\{\bw\in\R^n\mid U\bw\in\DmuS\}$.
For smooth functions $f,g:\DmuSr\to\R$, define
$$\{f,g\}_{\mu}(\bw):=\nabla f(\bw)^\top\Jop(U\bw;\mu)\nabla g(\bw),\qquad\forall\,\bw\in\DmuSr.$$
It can be easily verified that $\{\cdot,\cdot\}_{\mu}$ is bilinear, skew-symmetric and it satisfies the Leibniz rule
$\{f,gh\}_{\mu}=g\{f,h\}_{\mu}+h\{f,g\}_{\mu}$, for any $f,g,h\in C^{\infty}(\DmuSr;\R)$.
To show that
$\{\cdot,\cdot\}_{\mu}$ satisfies the Jacobi identity
we consider a basis
$\boldsymbol{\chi}_1,\ldots,\boldsymbol{\chi}_{n-1}\in\R^n$ of the orthogonal complement to $\operatorname{span}(P\cS)$, and set $\boldsymbol{\chi}_n:=\nabla S_r(\cdot;\mu)$. The vectors $\boldsymbol{\chi}_1,\ldots,\boldsymbol{\chi}_n$ form a basis of $\R^n$, since
$(P\cS)^\top\boldsymbol{\chi}_n(\bw)=\bS(U\bw;\mu)\neq 0$ for any $\bw\in\DmuSr$. Moreover, they can be written as gradients of the functions $\xi_i:\R^n\to\R$,
$$\xi_n(\bw)=S_r(\bw;\mu),\qquad \xi_i(\bw):=\boldsymbol{\chi}_i^\top\bw,\qquad\forall\, 1\leq i\leq n-1,\,\forall\, \bw\in\R^n.$$
%
For any $1\leq i\leq n-1$, since $(P\cS)^\top\boldsymbol{\chi}_i=0$, it holds
$\QS(U\bw;\mu)^\top \boldsymbol{\chi}_i=\boldsymbol{\chi}_i$. Moreover, from item (iii) in \Cref{lem:oblproj}, it follows that
\begin{equation*}    
    \QS(U\bw;\mu)^\top \boldsymbol{\chi}_n=\QS(U\bw;\mu)^\top \nabla S_r(\bw;\mu)=0,\qquad\forall\,\bw\in\DmuSr.
\end{equation*}
Therefore, for any $1\leq i,j\leq n-1$ and $\bw\in\DmuSr$,
\begin{equation*}
    \{\xi_i,\xi_j\}_{\mu}(\bw)=(\QS(U\bw;\mu)^\top \boldsymbol{\chi}_i)^\top \barJ(\QS(U\bw;\mu)^\top \boldsymbol{\chi}_j) = \boldsymbol{\chi}_i^\top \barJ \boldsymbol{\chi}_j.
\end{equation*}
Furthermore $\{\xi_i,\xi_j\}_{\mu}(\bw)=0$,
if either $i=n$ or $j=n$.
This implies that, for any $1\leq i,j,k\leq n$, and $\bw\in\DmuSr$, it holds
$\{\xi_i,\{\xi_j,\xi_k\}_{\mu}\}_{\mu}(\bw)=0$.
\end{proof}
\end{lemma}

\paragraph{Choice of $\ccH$ and $\ccS$}
The choice of $\ccX$, with $X\in\{H,S\}$, plays an important role in the derivation of the reduced order model.
One can ensure that $\ccX$ belongs to $\range(U)$ by reduced basis enrichment. However,for a large set of test parameters, allowing $\ccX$ to be parameter dependent may enlarge the domain of definition $\Dmu$ but require a ``large'' reduced space.
In practical applications, as shown in the numerical tests of \Cref{sec:num-exp}, it is often possible to set $\ccX$ independent of the parameter.

We also point out that a particular choice of $\ccH$ and $\ccS$ makes the formulation \eqref{eq:OP-ROM}, with $M=I_N$, equivalent to the energetically consistent reduced order model proposed in \cite{gruber2023energetically}.
In fact, taking $\cH=\bb e_{k_0}, \cS=\bb e_{k_1}$, for some suitable $k_0,k_1\leq N$, and assuming that the reduced basis $U$ contains the vectors $\bb{e}_{k_0},\bb e_{k_1}$, standard calculations show that the two methods lead to the same approximate solution trajectory provided both models are well-defined and the same initial condition is considered.
However, even in this particular case, the reduced order model \eqref{eq:OP-ROM} proposed in this work is naturally amenable to the hyper-reduction of the nonlinear terms, as it will be shown in \Cref{sec:OP-HROM}. Applying hyper-reduction to the reduced order model of \cite{gruber2023energetically} in its plain form requires instead a factorization of the semi-metric tensor $\GM(\bv)$, at any $\bv\in\R^N$, resulting in an online arithmetic complexity of $O(N^3)$.

\section{Hyper-reduction of the OP-ROM}\label{sec:OP-HROM}

Although the introduction of oblique projections ensures that the resulting reduced order model is almost-metriplectic, the presence of nonlinear functions might hinder its computational efficiency.
Indeed, evaluating \eqref{eq:OP-ROM} generally requires the reconstruction of the reduced state (so-called \textit{lifting} to the full space) and thus a computational cost that may still depend on the full order dimension. In order to mitigate this issue, we combine gradient-preserving EIM \cite{pagliantini2023gradient} for the hyper-reduction of $\nabla H_r,\nabla S_r,$ with Matrix DEIM (MDEIM) \cite{negri2015efficient} for the hyper-reduction of $\Jop,\Gop$, and construct oblique corrections from the hyper-reduced gradients. Under suitable sparsity assumptions, this allows for efficient evaluation while retaining the degeneracy conditions associated with the hyper-reduced functionals.



\subsection{Gradient-preserving hyper-reduction}\label{sec:HSeim}
We summarize the method proposed in \cite[Section 3]{pagliantini2023gradient} for the hyper-reduction of the gradient of $X_r$, with $X\in\{H,S\}$ denoting either the Hamiltonian or the entropy function, and we refer to the aforementioned reference for further details.

First, for any $\mu\in\p$, we write $X:\R^N\to\R$
as the sum of a quadratic function $\LX(\cdot;\mu):\R^N\to\R$ and a nonlinear term $\NX(\cdot;\mu):\R^N\to\R$, namely
\begin{equation*}
    X(\bv;\mu)= \LX(\bv;\mu)+ \NX(\bv;\mu),
    \qquad\forall\,\bv\in\R^N.
\end{equation*}

As discussed in \cite[Remark 3.1]{pagliantini2023gradient}, in most cases of interest, it is possible to decompose the nonlinear term as
\begin{equation}\label{eq:Nh}
    \NX(\bv;\mu)=\cNx^\top\,\bb{x}(\bv;\mu),
    \qquad\forall\,\bv\in\R^N,
\end{equation}
where $\bb{x}:\R^N\times\p\to\R^{d}$, $\bb{x}\in\{\bb{h},\bb{s}\}$,
is a smooth function, $\cNx\in\R^{d}$
is a constant vector, and for some $d\in\mathbb{N}$ that typically scales with $N$.

From the expression in \eqref{eq:Nh}, the gradient of the nonlinear component of $X$ reads
$\nabla \NX(\bv;\mu)=\Dx(\bv;\mu)^\top\cNx$,
where $\Dx(\cdot;\mu):\R^N\to\R^{d\times N}$ denotes the Jacobian matrix of $\bb{x}$.
The idea of the method proposed in \cite{pagliantini2023gradient} is to exploit this decomposition to hyper-reduce the Jacobian matrix of $\bb{h}$ mapped into the reduced space via EIM.
More in details, consider the EIM oblique projection
$ \pb_X\coloneq \Xi_X(P_X^\top\Xi_X)^{-1}P_X^\top\in\R^{d\times d}$,
where $\Xi_X\in\R^{d\times \nXeim}$, $\nXeim\leq d$, is the EIM basis obtained from the snapshots of $\Dx(\Pi_U\uf(t;\mu))U$ and $P_X\in\R^{d\times \nXeim}$ is the interpolation matrix which can be obtained as in, e.g., \cite[Algorithm 1]{chaturantabut2010nonlinear}.
The gradient of the nonlinear function $\NX$ mapped into the reduced space is then approximated as
\begin{equation*}
    U^\top \nabla\NX(\bv;\mu)
    =U^\top \Dx(\bv;\mu)^\top\cNx
    \approx U^\top \Dx^\top(\bv;\mu)\pb_X^\top\cNx,\qquad\forall\,\bv\in\R^N.
\end{equation*}
This can be equivalently seen
as an approximation of the nonlinear term $\NX$ in \eqref{eq:Nh} which generates 
the hyper-reduced function 
\begin{equation}\label{eq:Xhr}
    X_{\hr}(\bw;\mu)\coloneq \LX(U\bw;\mu)+ \NXhr(\bw;\mu),\quad
    \NXhr(\bw;\mu)\coloneq\cNx^\top\pb_X \bb{x}(U\bw;\mu).
\end{equation}
For the subsequent analysis, it is useful to introduce the quantity
$\Xeim:\R^N\times \p\to\R$, with $X\in\{H,S\}$, defined as
$\Xeim(\bv;\mu):=\LX(\bv;\mu)+ \cNx^\top\pb_X \bb{x}(\bv;\mu)$. This allows us to write the hyper-reduced function as $X_{\hr}= \Xeim \circ U$, and 
\begin{align*}
    \nabla X_{\hr}(\bw;\mu)=U^\top \nabla \Xeim(U\bw;\mu),\qquad\forall\,\bw\in\R^n.
\end{align*}

\begin{remark}\label[remark]{oss:GPeim}
To assess the arithmetic complexity of evaluating $\nabla\NXhr$, with $X\in\{H,S\}$,
assume that each entry of the vector-valued function $\bb{x}(\cdot;\mu):\R^N\to\R^d$ in \eqref{eq:Nh} depends on $s_{\bb{x}}\leq N$ entries of its variable. Denote by $\tcNx\in\R^{\nXeim}$ the constant vector $\tcNx\coloneq(P_X^\top \Xi_X)^{-\top}\Xi_X^\top\cNx$. This can be computed in, at most, $O(\nXeim^3)+O(\nXeim d)$ operations. The gradient of the hyper-reduced nonlinear part of the function $X$ reads
\begin{equation*}
   \nabla\NXhr(\bw;\mu)=U^\top \Dx^\top(U\bw;\mu)P_X\,\tcNx,\qquad\forall\,\bw\in\R^n.
\end{equation*}
Since $P_X=\begin{bmatrix}
\bb e_{i_1}&\ldots &\bb e_{i_m}
    \end{bmatrix},$ for some indices $\{i_1,\ldots,i_m\}\subset\{1,\ldots,d\}$, it follows that
    \begin{equation*}
        \Dx^\top(U\bw;\mu)P_X=\begin{bmatrix}
            \nabla \bb x_{i_1}(U\bw;\mu)&\ldots& \nabla \bb x_{i_m}(U\bw;\mu)
        \end{bmatrix}\in\R^{N\times \nXeim}.
    \end{equation*}
Each $\nabla \bb{x}_{i_j}$ requires at most $s_{\bb{x}}$ components of $U\bw$ to be reconstructed, and thus can be computed in $O(n s_{\bb{x}})$ operations. Since this is done for $\nXeim$ terms, the total online cost of this evaluation amounts to $O(n\nXeim s_{\bb{x}})$ operations.
Moreover, since each column of the above matrix has at most $s_{\bb{x}}$ non zero entries, the product $U^\top \left(  \Dx^\top(U\bw;\mu)P_X\right)$ can be performed in $O(n\nXeim s_{\bb{x}})$ operations as well. The final product with $\tcNx$ incurs in $O(n\nXeim)$ extra operations.
We also point out that, if the linear term of the gradient of $X$ does not admit an affine or another computationally efficient decomposition, it can be included in the nonlinear term and undergo hyper-reduction.
%
\end{remark}

\subsection{Hyper-reduction of the Poisson and semi-metric operators}


To preserve the metriplectic structure and degeneracy conditions while enabling evaluation at a cost independent of the full order dimension $N$, we replace the nonlinear operators in \eqref{eq:OP-ROM} with suitable hyper-reduced approximations, as follows.
Let $\DmuX$, with $X\in\{\Heim,\Seim\}$ be defined as in \eqref{eq:DmuX}. We introduce the sets
$$\Dmueim\coloneq\DmuHeim\cap\DmuSeim,\qquad
\Dmuhr\coloneq\{ \bw\in\R^n \mid \ U\bw\in \Dmueim \}.$$
We construct a hyper-reduced system, dubbed OP-HROM, of the form:
given $\mu\in\p$ and $\uhr(0;\mu)=P\bb{u}_0(\mu)\in\Dmuhr$, find $\uhr(\cdot;\mu)\in C^1(\mathcal{T};\R^n)$ such that
\begin{equation}\label{eq:hrOP-ROM}
    \duhr(t;\mu)=\Xophr(U\uhr(t;\mu);\mu),\qquad\forall\,t\in\mathcal{T}.
\end{equation}
where
\begin{equation*}
    \Xophr(U\uhr;\mu)\coloneq\Jophr(U\uhr;\mu)\nabla H_{\hr}(\uhr;\mu)+\Gophr(U\uhr;\mu)\nabla S_{\hr}(\uhr;\mu).
\end{equation*}
The operators $H_{\hr}$ and $S_{\hr}$ are given in \eqref{eq:Xhr}, while the hyper-reduced Poisson and semi-metric tensor are defined, for any $\bv\in \Dmueim$, as
\begin{equation}\label{eq:JGophr}
\begin{aligned}
    \Jophr(\bv;\mu)&\coloneq \QShr(\bv;\mu)\barJhr(\bv)\QShr(\bv;\mu)^\top,\\
    \Gophr(\bv;\mu)&\coloneq \QHhr (\bv;\mu)\barGhr(\bv)\QHhr(\bv;\mu)^\top,
\end{aligned}
\end{equation}
where $\QShr, \QHhr$ are as in \eqref{eq:PinQ} and $\barJhr,\barGhr$ are suitable approximations of the operators $\barJ,\barG$ appearing in \eqref{eq:JGop}.
This section is devoted to the derivation of such approximations via MDEIM.
We postpone to \Cref{oss:MDEIMproperties} a discussion on the preservation of the properties of the two operators.

Let $\vect:\R^{N\times N}\to\R^{N^2}$ be the vectorization operator that maps a matrix into the vector obtained by vertically stacking its columns, and let us denote by $\vect^{-1}$ the inverse operation.
The idea of MDEIM for the hyper-reduction of 
$\widehat{Y}$, with $Y\in\{J,G\}$,
is to apply standard DEIM to $\vect(\widehat{Y})$.
More in details, let $\mathbb{P}_Y\coloneq \Psi_Y(P_Y^\top\Psi_Y)^{-1}P_Y^\top\in\mathbb{R}^{N^2\times N^2}$ be the EIM projection, where $\Psi_Y\in\R^{N^2\times \nYmeim}$, $\nYmeim\ll N^2$, is obtained from snapshots of $\vect(\YM(\uf(t;\mu)))$, and $P_Y\in\R^{N^2\times \nYmeim}$ is the EIM interpolation matrix.
Then, $\barYhr$ in \eqref{eq:JGophr} is defined, for $Y\in\{J,G\}$, as
\begin{equation}\label{eq:barYhr}
\barYhr(\bv)=U^\top \vect^{-1}\big(\mathbb{P}_Y \vect(\widehat{Y}(\bv))\big) U,\qquad\forall\,\bv\in\R^N.
\end{equation}
The expression above is computationally as expensive to evaluate as the full order operator, even at a reduced state $U\bw$.
Therefore, its implementation should exploit an affine decomposition, as discussed in the proof of \Cref{lem:Jophr-cost} for $\Jophr$. An analogous result holds for $\Gophr$.


\begin{lemma}\label[lemma]{lem:Jophr-cost}
Assume that each entry of $\JM$ depends on $s_J\leq N$ components of its variable and that each entry of $\bb{s}(\cdot;\mu)$ in \eqref{eq:Nh} depends on $s_{\bb{s}}\leq N$ components of its variable. 
If $\cShr\in\range(U)$, the assembly of the hyper-reduced operator $\Jophr(U\bw;\mu)$ defined in \eqref{eq:JGophr} requires 
$O(n\nJmeim N^2)+O(\nSeim^3)+O(\nJmeim^3)+O(\nSeim d)$ operations independent of $\mu$ and $\bw$,
$O(nN)$ operations per instance of the parameter $\mu$, and
$O(\nJmeim n^2)+O(\nJmeim ns_J)+O(\nJmeim^2)+O(n^2)+O(n\nSeim s_{\bb{s}})$ operations per instance of the vector $\bw\in\R^n$.
\end{lemma}
\begin{proof}
Let $\bw\in\R^n$ and $\mu\in\p$ be fixed.
We distinguish two computational costs for the evaluation of $\Jophr$: 
one associated with the assembly and evaluation of $\PShr(U\bw;\mu)$, and the other with $\barJhr(U\bw)$. The arithmetic complexity of the remaining matrix-matrix multiplication
is at most $O(n^2)$.

Regarding the cost of the evaluation of $\PShr(U\bw;\mu)$ in \eqref{eq:PinQ}, $P=U^\top M\in\R^{n\times N}$ can be computed once with $O(nN^2)$ operations, while the product
$P\cShr\in\R^n$ needs to be computed for every instance of the parameter $\mu$ at a cost of $O(nN)$ operations.
Moreover, as shown in \Cref{oss:GPeim}, the cost of the evaluation of 
$\nabla S_{\hr}(\bw;\mu)=U^\top \nabla \Seim(U\bw;\mu)$ is $O(n\nSeim s_{\bb{s}})$ for each instance of $\bw$ plus $O(\nSeim^3)+O(\nSeim d)$ operations that can be done once.
Moreover, since by hypothesis $\cShr\in\operatorname{Range}(U)$, there exists $\boldsymbol{\xi}(\mu)\in\R^n$ such that $\cShr=U\boldsymbol{\xi}(\mu)$. Then, it holds
$\bShr(U\bw;\mu)
=\boldsymbol{\xi}(\mu)^\top\nabla S_{\hr}(\bw;\mu)$, where
$\boldsymbol{\xi}(\mu)\in\R^n$ can be computed once per each value of the parameter $\mu\in\p$ at an arithmetic complexity of $O(nN)$. 

For the evaluation of $\barJhr$ defined in \eqref{eq:barYhr} we introduce the vector-valued function $\boldsymbol{\eta}:\R^N\to\R^{\nJmeim}$ defined as
$\boldsymbol{\eta}(\bv):=(P_J^\top\Psi_J)^{-1}P_J^\top \vect(\JM(\bv))$, so that
\begin{equation}\label{eq:barJhr-affine}
\barJhr(U\bw)
=\sum_{i=1}^{\nJmeim} \eta_i(U\bw)U^\top \JM_i U,
\end{equation}
where $\JM_i\coloneq\vect^{-1}((\Psi_J)_{:,i})$ denotes the $i$th column of the matrix $\Psi_J$, for any $1\leq i\leq \nJmeim$.
This representation is more amenable to implementation since its evaluation has a computational cost that does not depend on the full dimension $N$. Indeed, for the evaluation of $\boldsymbol{\eta}(U\bw)$, one exploits the sparsity of $\JM$ so that the product $P_J^\top\vect(\JM(U\bw))$ can be evaluated, for each $\bw\in\R^n$, in $O(n\nJmeim s_J)$.
Since the matrix $(P_J^\top\Psi_J)^{-1}\in\R^{\nJmeim\times \nJmeim}$ can be precomputed offline,
the evaluation of $\boldsymbol{\eta}(U\bw)$ has arithmetic complexity
$O(\nJmeim^2)+O(\nJmeim ns_J)$.
Moreover, since the matrices $\JM_i$ are constant, the term $U^\top\JM_i U$ can be computed, for every $1\leq i\leq \nJmeim$, at a cost of $O(n N^2)$ operations. Given $\boldsymbol{\eta}(U\bw)$, the cost to assemble $\barJhr(U\bw)$ for any instance of $\bw$, is then $O(\nJmeim n^2)$.
%
\end{proof}


Note that \Cref{lem:Jophr-cost} tells us that, assuming suitable sparsity of the operators $\JM$, $\GM$, $\bb{s}$, and $\bb{h}$, the hyper-reduced dynamical system \eqref{eq:hrOP-ROM} can be solved online at a cost that is independent of $N$.

\subsection{Properties of the hyper-reduced model}

Owing to the structure of the hyper-reduced operators $\QHhr$ and $\QShr$ and of $\barJhr$, $\barGhr$, the hyper-reduced model \eqref{eq:hrOP-ROM} is (almost-)metriplectic, as shown in the next result, which is the analogous of \Cref{lem:OPmetriplectic}.
\begin{lemma}\label[lemma]{lem:hrproperties}
Let $\mu\in\p$ be fixed and assume that the hyper-reduced model \eqref{eq:hrOP-ROM} is well-posed in $\Dmuhr$ in the sense of \Cref{rmk:OPROMwellposed}.
If $\cX\in\range(U)$, for $X\in\{\Heim,\Seim\}$, and $\barGhr$ defined in \eqref{eq:barYhr} is positive semi-definite then
problem \eqref{eq:hrOP-ROM} is almost-metriplectic with respect to the hyper-reduced Hamiltonian and entropy functions.

Moreover, if $\barJhr\in\R^{n\times n}$ defined in \eqref{eq:barYhr} is constant in $\R^N$, then $\Jophr(\cdot;\mu)$ is a Poisson operator in $\DmuSeim$, and thus \eqref{eq:hrOP-ROM} is metriplectic. 
\end{lemma}
\begin{proof}
By taking the EIM basis $\Psi_J$ via truncated SVD of the snapshot matrix obtained from $\vect(\JM(\uf(t;\mu)))$ under variation of $t$ and $\mu$, its columns are linear combination of the columns of the snapshot matrix. Since each of these columns are vectorization of skew-symmetric matrices, the same holds for the columns of $\Psi_J$, so that $\barJhr$ is skew-symmetric.
A similar reasoning allows to deduce the symmetry of the operator $\barGhr$, and hence of $\Gophr$.



Concerning the degeneracy conditions, one can resort to
\Cref{lem:oblproj}, item (iii), which gives
$$\Jophr(U\bw;\mu)\nabla S_{\hr}(\bw;\mu)=\Gophr(U\bw;\mu)\nabla H_{\hr}(\bw;\mu)=0,\qquad\forall\,\bw\in \Dmuhr.$$
This implies that the hyper-reduced model is almost-metriplectic in $\Dmuhr$.



The fact that problem \eqref{eq:hrOP-ROM} is genuinely metriplectic when $\barJhr$ is constant, follows by the same steps as in the proof of \Cref{lem:OPmetriplectic}.
\end{proof}

The result above establishes that problem \eqref{eq:hrOP-ROM} is almost-metriplectic. As such, the energy is preserved and entropy is produced. However, this refers to the hyper-reduced approximations of the energy and entropy functionals
and it does not ensure that the exact functionals evaluated at the hyper-reduced trajectories will satisfy the conservation/dissipation properties.

\begin{remark}\label[remark]{oss:MDEIMproperties}
In general, there is no guarantee that the operator $\barGhr$ obtained via MDEIM approximation of $\barG$ will preserve positive (semi)-definiteness. In \cite[Section 3.5]{negri2015efficient}, it is suggested that the positivity of the reduced operator can be enforced by solving a convex optimization problem.
However, in addition to the extra computational cost of solving such problem online, the approximation error of the resulting operator $\barGhr$ is not proven to be bounded by the EIM approximation error.
As an alternative to the aforementioned strategy, we propose to project the reduced -- not necessarily positive semi-definite -- matrix $\barGhr$ defined in \eqref{eq:barYhr} onto the positive semi-definite matrices. To this aim, let us consider the spectral decomposition of $\barGhr=W\diag(\lambda_i)W^\top$, then
$$\barGhr^+:=W\diag(\max\{0,\lambda_i\})W^\top,$$
is its nearest positive semi-definite matrix in the Frobenius norm.
The hyper-reduction error is bounded, for all $\bv\in\R^N$, as
\begin{equation*}
\begin{aligned}
\norm{\barGhr^+(\bv)-\barG(\bv)}_2
&\leq \norm{\barGhr^+(\bv)-\barG(\bv)}_F
\leq \norm{\barGhr(\bv)-\barG(\bv)}_F\\
& \leq  \lambda_{\min}^{-1}(M)\norm{(P_G^\top \Psi_G)^{-1}}_2\norm{(I_{N^2}-\Psi_G\Psi_G^\top)\vect(\GM(\bv))}_2.
\end{aligned}
\end{equation*}
and then one can take $\Gophr$ in \eqref{eq:JGophr} as
$$\Gophr(\bv;\mu)=\QHhr (\bv;\mu)\barGhr^+(\bv)\QHhr(\bv;\mu)^\top,\qquad\forall\,\bv\in\R^N.$$

\end{remark}

\subsection{A priori error estimate of the hyper-reduced model}\label{sec:convergence}

Denote by $\norm{\cdot}_M$ the vector/operator norm induced by the symmetric positive-definite matrix $M\in\R^{N\times N}$, 
by $\norm{\cdot}_2$ the Euclidean norm and by $\inp{\cdot,\cdot}_2$ the corresponding inner product.
We establish a priori error estimates in the $L^2$-norm over the space $\mathcal{T}\times\p$, using the $M$-norm in the state variable, and this will be denoted as $\norm{\cdot}_{L^2_M(\mathcal{T}\times\p)}$.
\begin{theorem}\label{thm:convergence-hr}
For any $\mu\in\mathcal{P}$,
let $\big(\uf(\cdot;\mu),\jf(\cdot;\mu),\gf(\cdot;\mu)\big)$ be solution of the full order model \eqref{eq:MF} coupled with
the mutual degeneracy conditions \eqref{eq:degCond}.
Let $\uhr(\cdot;\mu)$ be solution of the hyper-reduced model \eqref{eq:hrOP-ROM} with initial condition $\uhr(0;\mu)=P\uf(0;\mu)$.
Assume that there exist $\cHhr,\cShr\in\range(U)$ such that $\uf(t;\mu)\in \Dmueim$, and that $\uhr(t;\mu)\in \Dmuhr$,
for all $t\in\mathcal{T}$.

Assume that there exist $\aj,\ag,\kH,\kS\geq 0$, independent of $n$, such that, for any $t\in\mathcal T$ and $\mu\in\p$,
\begin{equation}\label{eq:opbound-hr}
\norm{M^{-1}\widehat{J}(\uf(t;\mu))}_{\MtoM}\leq \aj,\qquad
\norm{M^{-1}\widehat{G}(\uf(t;\mu))}_{\MtoM}\leq \ag,
\end{equation}
and
\begin{nalign}\label{eq:fracbound-hr}
\dfrac{\norm{\cHhr}_M \max_{X\in\{H,S\}}\,\norm{M^{-1}\nabla \Xeim(\uf(t;\mu))}_M}{|\bHhr(\uf(t;\mu);\mu)|}&\leq\kH,\\[.5ex]
\dfrac{\norm{\cShr}_M\norm{M^{-1}\nabla \Heim(\uf(t;\mu))}_M}{|\bShr(\uf(t;\mu);\mu)|}&\leq\kS.
\end{nalign}

Let $\bb{e}:=P \uf-\uhr$.
Assume that there exist constants $C_e\in\mathbb R$ and $\Cu\geq0$ such that, for every $\mu\in\p$, and almost every $t\in\mathcal T$, it holds
\begin{equation}\label{eq:boundinp-hr}
\begin{aligned}
   \big\langle \bb{e}(t;\mu),\,
   \Xophr(\uf(t;\mu);\mu)-\Xophr(U\uhr(t;\mu);\mu)\big\rangle_2
   \leq \, C_e\,\norm{\bb{e}(t;\mu)}_2^2&\\
   +\Cu\norm{(I_N-\Pi_U)\uf(t;\mu)}_M\norm{\bb{e}(t;\mu)}_2&.
\end{aligned}
\end{equation}
Finally, assume there exists $C_{\hr}>0$, independent of $n$, such that, for every $\mu\in\p$ and almost every $t\in\mathcal T$, it holds
\begin{equation}\label{eq:boundyz-hr}
    \max_{X\in\{H,S\}}\norm{U^\top \nabla\Xeim(\uf(t,\mu);\mu)}_2 \leq C_{\hr}\,\lambda_{\min}(M),
\end{equation}
where $\lambda_{\min}(M)$ denotes the minimum eigenvalue of the matrix $M$.

Then, the error between the full order and hyper-reduced solutions satisfies
\begin{equation*}
\begin{aligned}
\norm{\bb \uf-U\uhr}&^2_{L_M^2(\mathcal T\times\p)}
    \leq  \big(1+7C(T,C_e)C_u^2\big)\norm{(I_N-\Pi_U)\uf}_{L_M^2(\mathcal{T}\times\p)}^2\\
    & + 7C(T,C_e)\Big(C_H^2\norm{(I_N-\Pi_U)\jf}_{L_M^2(\mathcal{T}\times\p)}^2
    + C_S^2\norm{(I_N-\Pi_U)\gf}_{L_M^2(\mathcal{T}\times\p)}^2\\
    & + \sum_{X\in\{H,S\}} C_{X}^2\norm{M^{-1}\Dx(\uf)^\top (I_d-\mathbb{P}_X^\top )\cNx}_{L_M^2(\mathcal{T}\times\p)}^2 \\
    & + \sum_{Y\in\{J,G\}} C_{Y,\hr}^2
    \norm{(I_{N^2}-\Psi_Y\Psi_Y^\top) \vect(\YM(\uf))}_{L^2(\mathcal{T}\times\p)}^2\Big),
\end{aligned}
\end{equation*}
where the constants $C_H,C_S,C_{J,\hr},C_{G,\hr} \geq 0$ and $C(T,C_e)$ are defined as
\begin{equation*}
\begin{aligned}
& C_H\coloneq\aj + \ag\,\kH(2+\kH),\qquad
C_S\coloneq\ag+2\aj\,\kS,\\[1ex]
& C_{J,\hr}\coloneq C_{\hr}(1+2\kS)\norm{(P_J^\top\Psi_J)^{-1}}_2,\quad
C_{G,\hr}\coloneq C_{\hr}(1+\kH(2+\kH))\norm{(P_G^\top\Psi_G)^{-1}}_2,\\[1ex]
& C(T,C_e)\coloneq
\begin{cases}
\displaystyle\dfrac{e^{2C_e T}-1}{4C_e^2}-\dfrac{T}{2C_e}\,,
&C_e\ne0,\\[2ex]
\displaystyle\frac{T^2}{2},&C_e=0.
\end{cases}
\end{aligned}
\end{equation*}
\end{theorem}
\begin{proof}
Let us decompose the model order reduction error as
\begin{equation*}
\begin{aligned}
    \uf(t;\mu)-U\uhr(t;\mu)&=\left(\uf(t;\mu)-UP\uf(t;\mu)\right)+U\left(P\uf(t;\mu)-\uhr(t;\mu)\right)\\
     &\eqqcolon \Pi_U^{\perp} \uf(t;\mu)+U\bb{e}(t;\mu).
\end{aligned}
\end{equation*}
where $\Pi_U^{\perp}:=I_N-\Pi_U$.
Since the projection error is $M$-orthogonal to $\range(U)$, 
\begin{equation}\label{eq:errMnorm-hr}
    \norm{\uf(t;\mu)-U\uhr(t;\mu)}^2_M = \norm{\Pi_U^{\perp}\uf(t;\mu)}^2_M+\norm{\bb{e}(t;\mu)}^2_2.
\end{equation}
Taking the time derivative of $\bb{e}$ gives
\begin{align*}
    \frac{d}{dt}\bb{e}(t;\mu)&=U^\top\JM(\uf(t;\mu))\jf(t;\mu)+ U^\top\GM(\uf(t;\mu))\gf(t;\mu) - \Xophr(U\uhr(t;\mu);\mu).
\end{align*}
By adding and subtracting the velocity field $\Xophr$ of the hyper-reduced model evaluated at the full-order state $\uf(t;\mu)$, we get
\begin{equation}\label{eq:err12-hr}
\begin{aligned}
    \frac{d}{dt}\bb{e}(t;\mu)= &\, \underbrace{U^\top\JM(\uf(t;\mu))\jf(t;\mu)
    +U^\top\GM(\uf(t;\mu))\gf(t;\mu)- \Xophr(\uf(t;\mu);\mu) }_{=:A}\\
    & +\Xophr(\uf(t;\mu);\mu)-\Xophr(U\uhr(t;\mu);\mu).
\end{aligned}
\end{equation}
Taking the inner product with $\bb{e}$, and using hypothesis \eqref{eq:boundinp-hr} yields
\begin{equation}\label{eq:dte-hr}
\begin{aligned}
    \dfrac12\dfrac{d}{dt}\norm{\bb{e}(t;\mu)}^2_2  &=
    \inp{\bb{e}(t;\mu),A}+\inp{\bb{e}(t;\mu),\Xophr(\uf(t;\mu);\mu)-\Xophr(U\uhr(t;\mu);\mu)}_2\\
    &\leq\norm{\bb{e}(t;\mu)}_2\big(\norm{A}_2+C_e\,\norm{\bb{e}(t;\mu)}_2
   +\Cu\norm{\Pi_U^{\perp}\uf(t;\mu)}_M\big).
\end{aligned}
\end{equation}
Next, we estimate the term $A$ defined in \eqref{eq:err12-hr} distinguishing a part $A_J$ associated with the Hamiltonian dynamics and a part $A_G$ associated with the gradient flow.

Let us introduce the quantities
$$\jeim(t;\mu):=M^{-1}\nabla\Heim(\uf(t,\mu);\mu),\qquad
\geim(t;\mu):=M^{-1}\nabla\Seim(\uf(t,\mu);\mu).$$
The idea is to separate the error in the approximation introduced by the hyper-reduction of the gradient of $H$ and $S$ and the one associated with the hyper-reduction of the operators $J$ and $G$.

In the remaining of the proof, whenever clear from the context, we omit the dependence on $(t;\mu)$.
We start by adding and subtracting in the expression for $A$ the term $\sum_{Y\in\{J,G\}}\big(\barY(\uf)-\barYhr(\uf)\big)P\yeim$,
where $\yeim$ equals either $\jeim$ or $\geim$ according to the value of $Y$.
This results in
\begin{equation*}
\begin{aligned}
A =\sum_{Y\in\{J,G\}}
\Big(U^\top\YM(\uf)\big(\yf-\Pi_U\yeim\big)
+\big(\barY(\uf)-\barYhr(\uf)\big)P\yeim
+ A_Y\Big)
\end{aligned}
\end{equation*}
where $A_Y:=(\barYhr(\uf)-\Yophr(\uf;\mu))P\yeim$.
First, one can observe that $\norm{P\bv}_2\leq\norm{\bv}_M$ for any $\bv\in\R^N$, with the equality whenever $\bv\in\range(U)$. Then, the inequalities in \eqref{eq:opbound-hr} give, for $Y\in\{J,G\}$ and $\bv=\yf-\Pi_U\yeim$,
\begin{equation}\label{eq:PM}
\norm{PM^{-1}\YM(\uf)\bv}_2
    \leq \norm{M^{-1}\YM(\uf)\bv}_M
    \leq \ay \norm{\bv}_M.
\end{equation}


\noindent
\textbf{Bound on $A_{J}$.}
By definition of $\Jophr$, skew-symmetry of $\barJhr$, and \eqref{eq:PinQ}, yields
\begin{equation*}
    \barJhr(\uf)-\Jophr(\uf;\mu)=
    \PShr(\uf;\mu)\barJhr(\uf)-\big(\PShr(\uf;\mu)\barJhr(\uf)\big)^{\top}
\end{equation*}
and
\begin{equation*}
    \PShr(\uf;\mu)\barJhr(\uf)=
    \bShr^{-1}(\uf;\mu)
    U^\top M \cShr\nabla \Seim(\uf;\mu)^\top U\barJhr(\uf).
\end{equation*}
Therefore, the skew-symmetry of $\barJhr$ yields
\begin{align*}
\norm{A_J}_2&=\norm{\big(\barJhr(\uf)-\Jophr(\uf;\mu)\big)P\jeim}_2\\
&\leq
2\dfrac{\norm{P\cShr}_2\norm{P\jeim}_2}{|\bShr(\uf;\mu)|}\norm{\barJhr(\uf)P\geim}_2.
\end{align*}
The degeneracy conditions \eqref{eq:degCond} applied to the solution of the full order model give $\JM(\uf)\gf=0$, and this implies
$\barJ(\uf)P\geim = -U^{\top}\JM(\uf)(\gf-\Pi_U\geim)$.
Hence,
\begin{equation*}
    \norm{\barJhr(\uf)P\geim}_2\leq \norm{U^\top \JM(\uf)(\gf-\Pi_U\geim)}_2+\norm{\big(\barJ(\uf)-\barJhr(\uf;\mu)\big)P\geim}_2.
\end{equation*}
This bound, together with
assumption \eqref{eq:fracbound-hr}, gives
\begin{equation}\label{eq:AJbound}
\begin{aligned}
\norm{A_J}_2
&\leq 2\kS\norm{\barJhr(\uf)P\geim}_2\\
&\leq 2\kS\big(\aj \norm{\gf-\Pi_U\geim}_M+\norm{\big(\barJ(\uf)-\barJhr(\uf;\mu)\big)P\geim}_2\big).
\end{aligned}
\end{equation}

\noindent
\textbf{Bound on $A_G$.}
%
Observe that, by definition of $\Gophr$ and \eqref{eq:PinQ}, it holds
\begin{equation*}
\begin{aligned}
    \barGhr(\uf)-\Gophr(\uf;\mu)=&\,
\PHhr(\uf;\mu)\barGhr(\uf)+\big(\PHhr(\uf;\mu)\barGhr(\uf)\big)^{\top}\\
    & -\PHhr(\uf;\mu)\barGhr(\uf)\PHhr(\uf;\mu)^\top
\end{aligned}
\end{equation*}
where each term can be further expressed using
\begin{equation*}
    \PHhr(\uf;\mu)\barGhr(\uf)
    =\bHhr^{-1}(\uf;\mu)P \cHhr(P\jeim)^\top \barGhr(\uf).
\end{equation*}
Using the two expressions above, one can find the bound
\begin{align*}
\norm{A_G}_2 
&\leq 
\dfrac{\norm{P\cHhr}_2\norm{P\geim}_2}{|\bHhr(\uf;\mu)|}\\&\quad\cdot\left(2+\dfrac{\norm{P\cHhr}_2\norm{P\jeim}_2}{|\bHhr(\uf;\mu)|}\right)
\norm{\barGhr(\uf)P\jeim}_2.
\end{align*}
The fact that $\norm{P\bv}_2\leq\norm{\bv}_M$ for any $\bv\in\R^N$, together with assumption \eqref{eq:fracbound-hr}, gives
\begin{equation*}
\norm{A_G}_2=\norm{\big(\barGhr(\uf)-\Gophr(\uf;\mu)\big)P\geim}_2
\leq \kH(2+\kH)\norm{\barGhr(\uf)P\jeim}_2,
\end{equation*}
where the last term can be bounded in the same way as the corresponding term $\norm{\barJhr(\uf)P\geim}_2$, yielding
\begin{equation}\label{eq:AGbound}
\norm{A_G}_2
\leq \kH(2+\kH)\big(\ag \norm{\jf-\Pi_U\jeim}_M+\norm{\big(\barG(\uf)-\barGhr(\uf;\mu)\big)P\jeim}_2\big).
\end{equation}

Putting together \eqref{eq:PM}, the bound for $A_J$ in \eqref{eq:AJbound}, the bound for $A_G$ in \eqref{eq:AGbound}, and assumption \eqref{eq:boundyz-hr}, one gets
\begin{equation}\label{eq:Abound1}
\begin{aligned}
    \norm{A}_2\leq &\, C_H\norm{\jf-\Pi_U\jeim}_M
    + C_S \norm{\gf-\Pi_U\geim}_M\\
    & + C_{\hr}\lambda_{\min}(M)\Big((1+2\kS)\norm{\barJ(\uf)-\barJhr(\uf;\mu)}_2\\
    &+ (1+\kH(2+\kH))\norm{\barG(\uf)-\barGhr(\uf;\mu)}_2\Big).
\end{aligned}
\end{equation}
The error in the hyper-reduction of the operator $\barY$, for $Y\in\{J,G\}$, can be bounded as
\begin{equation*}
\begin{aligned}
    \norm{\barY(\uf)-\barYhr(\uf;\mu)}_2
    & \leq \norm{U}_2^2\norm{\YM(\uf)-\vect^{-1}(\mathbb{P}_Y\vect(\YM(\uf)))}_2\\
    &\leq \lambda_{\min}^{-1}(M)\norm{\mathbb{P}_Y}_2\norm{(I_{N^2}-\Psi_Y\Psi_Y^\top)\vect(\YM(\uf))}_2,
\end{aligned}
\end{equation*}
while the EIM error on the hyper-reduction of the Hamiltonian and entropy gradients enters the first two terms of \eqref{eq:Abound1} as follows
\begin{equation*}
\begin{aligned}
\norm{\jf-\Pi_U\jeim}_M
&\leq \norm{\Pi_U^{\perp}\jf}_M
  +\norm{U^\top (\nabla H(\uf;\mu)-\nabla \Heim(\uf;\mu))}_2\\
&\leq\norm{\Pi_U^{\perp}\jf}_M+\norm{M^{-1} \Dh(\uf;\mu)^\top(I_d-\mathbb{P}_H^\top)\cNh}_M.
\end{aligned}
\end{equation*}
This implies that the bound for the 2-norm of $A$ in \eqref{eq:Abound1} can be further written as
\begin{equation}\label{eq:Abound2}
\begin{aligned}
    \norm{A}_2\leq &\, C_H\norm{\Pi_U^{\perp}\jf}_M
    + C_S \norm{\Pi_U^{\perp}\gf}_M
    + C_H\, e_{H,\eim}(t;\mu)
    + C_S\, e_{S,\eim}(t;\mu)\\
    & + C_{\hr}(1+2\kS)\norm{(P_J^\top\Psi_J)^{-1}}_2\, e_{J,\hr}(t;\mu)\\
    & + C_{\hr}(1+\kH(2+\kH))\norm{(P_G^\top\Psi_G)^{-1}}_2\, e_{G,\hr}(t;\mu),
\end{aligned}
\end{equation}
where
\begin{align*}
e_{X,\eim}(t;\mu)&\coloneq\lVert M^{-1}\Dx(\uf(t;\mu))^\top (I_d-\mathbb{P}_X^\top )\cNx\rVert_M,& X\in\{H,S\},\\
e_{Y,\hr}(t;\mu)&\coloneq\norm{(I_{N^2}-\Psi_Y\Psi_Y^\top) \vect(\widehat{Y}(\uf(t;\mu)))}_2, & Y\in\{J,G\}.
\end{align*}
Inserting the bound \eqref{eq:Abound2} in the expression for the time derivative of the norm of the error in \eqref{eq:dte-hr} gives
\begin{equation*}
    \dfrac{d}{dt}\norm{\bb{e}(t;\mu)}_2 
    \leq (e_p(t;\mu)+e_{\hr}(t;\mu)) + C_e\,\norm{\bb{e}(t;\mu)}_2,
\end{equation*}
where
\begin{equation}\label{eq:epehr}
\begin{aligned}
    e_p & :=\Cu\norm{\Pi_U^{\perp}\uf(t;\mu)}_M
    +C_H\norm{\Pi_U^{\perp}\jf}_M
    +C_S\norm{\Pi_U^{\perp}\gf}_M,\\
    e_{\hr} & :=C_H\, e_{H,\eim}+C_S\, e_{S,\eim}+C_{J,\hr}\, e_{J,\hr}+ C_{G,\hr}\,e_{G,\hr}.
\end{aligned}
\end{equation}
Since $\bb{e}(0;\mu)=0$, the Gr\"onwall inequality together with the Cauchy-Schwarz inequality yields
\begin{equation}\label{eq:inte}
\begin{aligned}
\int_0^T \norm{\bb{e}(t;\mu)}_2^2\, dt & \leq
\int_0^T \left(\int_0^t e^{C_e(t-s)}(e_p(s;\mu)+e_{\hr}(s;\mu))\,ds\right)^2\, dt\\
&\leq
\int_0^T
\left(\int_0^t e^{2C_e(t-s)}\,ds\right)
\left(\int_0^t (e_p(s;\mu)+e_{\hr}(s;\mu))^2\,ds\right)dt\\
&\leq
C(T,C_e)\int_0^T (e_p(t;\mu)+e_{\hr}(t;\mu))^2\,dt.
\end{aligned}
\end{equation}
The conclusion follows by integrating the $M$-norm of the error in \eqref{eq:errMnorm-hr} over $[0,T]\times\p$, and using \eqref{eq:inte} together with the expression for the error components in \eqref{eq:epehr}.
\end{proof}

\begin{remark}
    In the absence of hyper-reduction, the solution of the reduced model \eqref{eq:OP-ROM} can be shown to converge to the solution of the full order model with a result analogous to \Cref{thm:convergence-hr}. Hypotheses \eqref{eq:fracbound-hr} and \eqref{eq:boundinp-hr} need to be replaced by the ones corresponding to the original Hamiltonian $H$ and entropy $S$, while hypothesis \eqref{eq:boundyz-hr} is not needed. Then the a priori estimate on the error is analogous with no term associated with the hyper-reduction.
\end{remark}

\section{Model order reduction via a block diagonal approach}\label{sec:block}

In this section we consider a specific type of metriplectic systems where some components of the state variables are characterized by a purely Hamiltonian dynamics and the remaining components by a gradient flow dynamics. An example of this situation will be shown in the numerical experiments in \Cref{sec:thermospring}.

To simplify the notation, we consider problem \eqref{eq:MF} with $M$ equal to the identity matrix.
%
%
Let $N_1, N_2\in\mathbb{N}$ be such that $N_1+N_2=N$ and assume that, for any $\bv\in\R^N$, the operators $J$ and $G$ in \eqref{eq:MF} are of the form
\begin{equation}\label{eq:JGstruct}
    J(\bv)=\begin{bmatrix}
        J_1(\bv) & \\
        & 0_{N_2\times N_2}
    \end{bmatrix}\in\R^{N\times N},\qquad
    G(\bv)=\begin{bmatrix}
        G_1 & \\
        & G_2(\bv)
    \end{bmatrix}\in\R^{N\times N},
\end{equation}
where $J_1(\bv)\in\R^{N_1\times N_1}$ is a skew-symmetric and invertible matrix (thus $N_1$ is even), $G_1\in\R^{N_1\times N_1}$ is symmetric positive semi-definite,
and $G_2(\bv)\in\R^{N_2\times N_2}$ is symmetric positive semi-definite.

Based on the partitioning introduced in \eqref{eq:JGstruct}, let us write any generic vector $\bv\in\R^N$ as $\bv=(\bv_1^\top,\bv_2^\top)^\top$ with $\bv_1\in\R^{N_1}$ and $\bv_2\in\R^{N_2}$.
Owing to the invertibility of $J_1(\bv)$, for any $\bv\in\R^N$, the metriplectic degeneracy conditions \eqref{eq:degCond}
imply that
\begin{equation*}\label{eq:G1degCond}
\begin{aligned}
    \nabla_{1}S(\bv;\mu)&=0,&\quad\forall\,\bv\in\R^N,\\
    G_1\nabla_{1}H(\bv;\mu)&=0,\qquad G_2(\bv)\nabla_{2}H(\bv;\mu)=0,& \quad\forall\,\bv\in\R^N,
\end{aligned}
\end{equation*}
where we use the notation $\nabla_{k}$, for $k\in\{1,2\}$, to denote the gradient with respect to the $k$th component of the variable, of length $N_k$.
A metriplectic system with $J$ and $G$ of the form \eqref{eq:JGstruct} reads
\begin{equation*}
\dot{\uf}(t,\mu)=\begin{bmatrix}
    \dot{\uf}_1(t,\mu)\\[.5ex]
    \dot{\uf}_2(t,\mu)
\end{bmatrix}=
\begin{bmatrix}
    J_1(\uf(t,\mu))\nabla_1 H(\uf(t,\mu);\mu)\\[.5ex]
    G_2(\uf(t,\mu))\nabla_2 S(\uf(t,\mu);\mu)
\end{bmatrix}
\end{equation*}
with $S$ that depends only on $\bv_2$,
$\nabla_{1}H(\bv;\mu)\in\ker(G_1)$, and
$\nabla_{2}H(\bv;\mu)\in\ker(G_2(\bv))$,
for any $\bv\in\R^N$. Note that $G_1$ does not play a role in the characterization of the dynamics but only of the degeneracy conditions.

For the reduced basis $U\in\R^{N\times n}$, with $n\leq N$, we consider the following block diagonal form:
for $n\in\mathbb{N}$ fixed, let $n_1, n_2\in\mathbb{N}$ be such that $n_1+n_2=n$ and $n_k\leq N_k$, for $k\in\{1,2\}$, then
\begin{equation}\label{eq:Ublock}
    U=\begin{bmatrix}
        U_1 &\\
         & U_2
    \end{bmatrix}\in\R^{N\times n},\qquad U_1\in\R^{N_1\times n_1},\; U_2\in\R^{N_2\times n_2},
\end{equation}
with $U_1^\top U_1 = I_{n_1}$ and $U_2^\top U_2=I_{n_2}$.
Note that the block structure of $U$ in \eqref{eq:Ublock} is meant to be such that the first block $U_1$ provides the reduction corresponding to the dynamics associated to $\bv_1$, while the second block is associated with the reduction of the dynamics of $\bv_2$.

Concerning the reduced metriplectic conditions \eqref{eq:rdegCond},
using the definition of $\Jr$ and the partitioning in \eqref{eq:JGstruct} gives, for any $\bw\in\R^n$,
\begin{equation*}
    \Jr(\bw)\nabla S_r(\bw;\mu)=\begin{bmatrix}
        U_1^\top J_1(U\bw)U_1 &\\
        & 0_{n_2\times n_2}
    \end{bmatrix}\begin{bmatrix}
        0_{n_1\times 1} \\
        U_2^\top\nabla_{2}S(U\bw;\mu)
    \end{bmatrix}=0.
\end{equation*}
This means that the first of the reduced metriplectic conditions \eqref{eq:rdegCond} is satisfied. Concerning the second condition, one has that, for any $\bw\in\R^n$,
\begin{equation}\label{eq:GrdegCond}
    G_r(\bw)\nabla H_r(\bw;\mu)=
    \begin{bmatrix}
        U_1^\top G_1\Pi_{U_1}\nabla_{1}H(U\bw;\mu)\\
        U_2^\top G_2(U\bw)\Pi_{U_2}\nabla_{2}H(U\bw;\mu)
    \end{bmatrix},
\end{equation}
where $\Pi_{U_k}:=U_kU_k^\top\in\R^{N_k\times N_k}$, for $k\in\{1,2\}$.
Since $G_k$ is positive semi-definite, the above vector vanishes if and only if $\Pi_{U_k}\nabla_{k}H(U\bw;\mu)\in\ker(G_k(U\bw))$ for any $\bw\in\R^n$, $k\in\{1,2\}$. A sufficient condition is given by
\begin{equation}\label{eq:kern-cond}
\Pi_{U_1}\ker(G_1)\subseteq \ker(G_1),\qquad
\Pi_{U_2}\ker(G_2(U\bw))\subseteq \ker(G_2(U\bw)),\qquad\forall\,\bw\in\R^n.
\end{equation}

\paragraph{Construction of the reduced basis}
Assume that $\SM\in\R^{N\times N_S}$ is a snapshot matrix, that is, $\SM$ collects the solution $\uf(t;\mu)$ of the full order model \eqref{eq:MF}, with $M=I_N$, for different values of time $t$ and parameter $\mu$. Based on the partitioning introduced in \eqref{eq:JGstruct}, let us write the snapshot matrix as
    \begin{equation*}
        \SM=\begin{bmatrix}
            \SM_1\\
            \SM_2
        \end{bmatrix}, \qquad \SM_1\in\R^{N_1\times N_S}, \quad \SM_2\in\R^{N_2\times N_S}.
    \end{equation*}
To ensure that the first condition in \eqref{eq:kern-cond} is satisfied, a possible strategy is to choose $U_1$ of the form $U_1=QY$, where $Y\in\R^{N_1\times n_1}$ needs to be suitably constructed and $Q\in\R^{N_1\times N_1}$ is an orthogonal matrix that diagonalizes the symmetric positive semi-definite matrix $G_1$. More in details,
let $r:=\rank(G_1)$ and let $Q=[Q_+\quad Q_0]\in\R^{N_1\times N_1}$, with $Q_+\in\R^{N_1\times r}$ an orthonormal basis of $\range(G_1)$, and $Q_0\in\R^{N_1\times (N_1-r)}$ an orthonormal basis of $\ker(G_1)$, so that
$G_1=Q_+\Lambda Q_+^\top$ with $\Lambda\in\R^{r\times r}$ diagonal positive definite.

For $n_1\in\mathbb{N}$ fixed, let $n_+, n_0\in\mathbb{N}$ be such that $n_++n_0=n_1$, and let $Y\in\R^{N_1\times n_1}$ be of the form
\begin{equation}\label{eq:Ymatrix}
    Y=\begin{bmatrix}
        Y_+ &\\
        & Y_0
    \end{bmatrix},\qquad Y_+\in\R^{r\times n_+},\quad Y_0\in\R^{(N_1-r)\times n_0},
\end{equation}
so that $U_1=[Q_+Y_+\quad Q_0Y_0]$. With this choice of $U_1$, one can easily verify that $\Pi_{U_1}\bb{x}$ belongs to $\ker(G_1)$ for any $\bb{x}\in\ker(G_1)$, and thus the first component $U_1^\top G_1\Pi_{U_1}\allowbreak\nabla_1 H(U\bw;\mu)$ of \eqref{eq:GrdegCond} vanishes for any $\bw\in\R^n$.

In \eqref{eq:Ymatrix}, $Y_+$ is meant to reduce the state components in the range of $G_1$, while $Y_0$ those in its kernel.
Since the space spanned by the columns of $U_1=Q Y$ should provide an approximation space for the first component of the solution vector, a natural idea is to construct $U_1$ from the snapshot matrix $\SM_1$ or, equivalently, construct $Y$ from the matrix $Q^\top \SM_1\in\R^{N_1\times N_S}$. More precisely, $Y$ is taken of the form \eqref{eq:Ymatrix} with $Y_+$ and $Y_0$ obtained via POD or a greedy strategy from $Q_+^\top \SM_1$ and $Q_0^\top \SM_1$, respectively.

Concerning the construction of the matrix $U_2\in\R^{N_2\times n_2}$ in \eqref{eq:Ublock}, several situations can be considered.
If $G_2$ has a kernel that is independent of the state variable, the construction above applies mutatis mutandis to $U_2$.
In particular, if $G_2(\uf)$ is invertible, then $U_2$ can be obtained from a standard POD or greedy strategy applied to the matrix $\SM_2$, with $n_2\leq N_2$ chosen according to the corresponding algorithm's criterion.
%
In the more general case where $G_2(\uf)$ has a state-dependent kernel, there is no general recipe to ensure that the second reduced metriplectic condition is satisfied (equivalently, that \eqref{eq:GrdegCond} vanishes). However, a possible way is to set $U_2=I_{N_2}$ which entails that the second component of the state variable is not reduced.

We point out that the above procedure builds a block diagonal reduced basis matrix which assures almost-metriplecticity of the reduced order model obtained by a standard Galerkin projection of the full order model \eqref{eq:MF} under the assumption that the operators $J$ and $G$ have the form \eqref{eq:JGstruct}. Since the construction of $U$ is entirely carried out in the offline phase, solving the corresponding reduced order model introduces no additional online costs compared with a standard Galerkin reduced order model of the same dimension.

\paragraph{Approximation properties of the reduced basis space}
Consider problem \eqref{eq:MF}-\eqref{eq:degCond} with $M=I_N$, and $J$ and $G$ of the form \eqref{eq:JGstruct} with $G_2$ having a state-independent kernel.
Assume that the reduced basis $U$ is constructed from a POD strategy according to the steps described in the previous paragraph. We use the notation introduced above and we add a subscript $k$ to all quantities associated with the block $U_k$, for $k\in\{1,2\}$.
It is possible to show that the projection error of the snapshot matrix $\SM\in\R^{N\times N_S}$, where $N_S$ is the number of snapshots, into the reduced basis space satisfies
    \begin{equation*}
        \norm{\SM-UU^\top \SM}_F^2\leq
        \sum_{j=m_1+1}^{\rank(\SM_1)}\sigma_j(\SM_1)^2
        +\sum_{j=m_2+1}^{\rank(\SM_2)}\sigma_j(\SM_2)^2
        \leq \sum_{j=\min\{m_1,m_2\}+1}^{\rank(\SM)}\sigma_j(\SM)^2,
    \end{equation*}
where $m_k:=\min\{n_{k,+},n_{k,0}\}$, for $k\in\{1,2\}$, and
$\sigma_j(\SM)$ denotes the $j$th singular value of the matrix $\SM$ assumed to be taken in decreasing order.
Indeed, owing to the definition of $U_1$ and $U_2$, and by the Eckart–Young–Mirsky theorem, it holds, for any $k\in\{1,2\}$,
\begin{equation*}
\begin{aligned}
\norm{\SM_k-U_kU_k^\top \SM_k}_F^2 & =
\sum_{\sigma\in\{0,+\}}
\norm{Q_{k,\sigma}(I-Y_{k,\sigma}Y_{k,\sigma}^\top)Q_{k,\sigma}^\top \SM_k}_F^2\\
& =\sum_{j>n_{k,+}}\sigma_j(Q_{k,+}^\top\SM_k)^2
+\sum_{j>n_{k,0}}\sigma_j(Q_{k,0}^\top\SM_k)^2
\leq \sum_{j>m_k}\sigma_j(\SM_k)^2.
\end{aligned}
\end{equation*}


\section{Numerical experiments}\label{sec:num-exp}
In this section, we present two numerical examples to assess the effectiveness of the proposed reduction strategies and their conservation properties.
A more detailed description of the test cases and of the implementation of the reduction and hyper-reduction methods is provided in the supplementary material.

In the following experiments the approximation of the solution of the generic metriplectic system \eqref{eq:pbm}-\eqref{eq:degCond} is computed from the solution $\uf_{\Lambda}$ of a given (hyper-)reduced model as $\uf_0(\mu)+U\uf_{\Lambda}(t;\mu)$ for all $t\in\mathcal{T}$, where $\uf_0(\mu)$ is the initial condition of \eqref{eq:pbm} corresponding to the parameter $\mu\in\p$. Accordingly, the (hyper-)reduced initial condition is always zero, implying that the reduced Hamiltonian coincides exactly with the full Hamiltonian at the initial time instance, as shown in \Cref{prop:rMFproperties} and subsequent discussion.

We compare the performances of the methods proposed in this work with a projection-based reduced order model (POD-ROM) obtained by Galerkin projection of the full order model \eqref{eq:pbm} onto a global POD basis, and with the energetically consistent reduced order model (EC-ROM) proposed in \cite{gruber2023energetically}.
The reduced order model \eqref{eq:OP-ROM} based on oblique projections will be dubbed OP-ROM, while the one based on the block construction of \Cref{sec:block} will be denoted as BD-ROM.

All implemented methods use the same snapshot matrix to construct the corresponding reduced basis via a POD strategy. In particular, for a fixed $n$, POD-ROM and EC-ROM have the same reduced basis, while the ones for BD-ROM and OP-ROM are constructed to satisfy the corresponding constraints.
%
More in details, the reduced basis of each model is constructed from snapshots of the solution
$(\uf,\jf,\gf)$ of the full order model \eqref{eq:MF}, obtained under variation of time and parameter. In particular,
consider a discrete parameter domain $\p_{\train}\subset\p$ and a set of time instants $\mathcal{T}_{\train}:= \{t_1,\ldots,t_{N_{\train}}\}$, and let $N_S=|\p_{\train}|N_{\train}$.
Once the solutions of the full order model corresponding to the sets above have been computed with a suitable time integration scheme, we collect them in a matrix $\SM\in\R^{N\times 3N_S}$ of the form
\begin{equation}\label{eq:snap-matrix}
    \SM = \begin{bmatrix}
    \ldots\, \uf(\theta_{ij})-\uf_0(\mu_j)\, \ldots & \ldots\, \jf(\theta_{ij})\, \ldots & \ldots\, \gf(\theta_{ij})\,\ldots
\end{bmatrix},
\end{equation}
where $\theta_{ij}=(t_i,\mu_j)\in\mathcal{T}_{\train}\times\p_{\train}$.
Note that, since in both examples considered the entropy is linear, its gradient $\nabla S$ is constant and it is thus included once as the final column of the snapshot matrix, so that $\SM\in\R^{N\times (2N_S+1)}$.

For numerical time integration, linear multistep methods are employed via the LSODA scheme \cite{hindmarsh1983odepack,petzold83} with tolerance of $10^{-14}$, as in \cite{gruber2023energetically}. 
Since this is an adaptive time-stepping, the numerical solutions obtained with different models are evaluated on
a fixed partition of the temporal interval $\mathcal{T}_h:= \{t_0=0,\ldots,t_{N_T}=T\}$,
Via piecewise polynomial interpolation 
of a degree consistent with the integration method.


To assess the convergence of the various methods we consider the following errors.
Given $\mu\in\p$,
let $\uf^{\idx}(\mu)\in\R^N$ be solution of the full order model at time $t_{\idx}$, with $0\leq\tau\leq N_T$, obtained with the aforementioned LSODA solver, and
let $\uf_\Lambda^{\idx}(\mu)\in\R^n$ be solution of a given reduced or hyper-reduced model (indicated by $\Lambda$) at time $t_{\idx}$.
We consider the relative errors
\begin{equation}\label{eq:relerr}
\errmu\coloneq\max_{1\leq \idx\leq N_T}\errmu(t_{\idx}),\quad
\errmu(t_{\idx})\coloneq\frac{\norm{\uf^{\idx}(\mu)-U\uf_\Lambda^{\idx}(\mu)-\uf_0(\mu)}_2}{\norm{\uf^{\idx}(\mu)}_2},\quad t_{\idx}\in\mathcal{T}_h.
\end{equation}
To assess the behavior of the error in the Hamiltonian and entropy functions evaluated at the approximate solution we introduce, for $X\in\{H,S\}$ and $t_{\idx}\in\mathcal{T}_h$, the errors
\begin{equation}\label{eq:relerrF}
\errmuF\coloneq\max_{1\leq \idx \leq N_T}\errmuF (t_{\idx}),\quad
\errmuF(t_{\idx})\coloneq\frac{|X(\uf^{\idx}(\mu)) - X(U\uf_\Lambda^{\idx}(\mu)+\uf_0(\mu))|}{|X(\uf^{\idx}(\mu))|}.
\end{equation}
Finally, to study the conservation of the Hamiltonian along the approximate trajectories we define
\begin{equation}\label{eq:presrelerr}
    \Dmuh(t_{\idx})\coloneq \frac{|H(U\uf_\Lambda^{\idx}(\mu)+\uf_0(\mu)) - H(\uf_0(\mu))|}{|H(\uf_0(\mu))|},\quad t_{\idx}\in\mathcal{T}_h.
\end{equation}

All experiments were presented below were coded in Python 3.12.3 and run on an Acer Swift SFX14-41G with an AMD Ryzen~7 5700U CPU and 8GB of RAM.

\subsection{Two thermo-spring system in entropy variables}\label{sec:thermospring}
As first test, we consider the entropy-based metriplectic formulation of the two thermo-spring system introduced in \cite[Section 4.3.1]{martin2016energy}. The system consists of two point masses connected by two thermo-springs that exchange heat. The first spring connects the first mass to a fixed point, while the second spring connects the two masses. In order to test model order reduction on this problem, we retain its structure but we increase the mechanical dimension $R$. 
The state vector is given by
\begin{equation*}
    \uf=\begin{bmatrix}
        \bb q_1^\top &\bb q_2^\top & \bb p_1^\top &\bb p_2^\top &S_1 &S_2
    \end{bmatrix}^\top
    \in\R^N,
\end{equation*}
where $N\coloneq4R+2$, $\bb q_1,\bb q_2\in\R^R$ are the position vectors, 
$\bb p_1,\bb p_2\in\R^R$ are the corresponding momenta, and the thermodynamic state is represented by the entropies $S_1,S_2$. 
Let $D\subseteq\R^N$ be the set of state vectors in $\R^N$ such that
$\bb q_1\neq0$, and $\bb q_1\neq\bb q_2$, and
let us consider a generic element $\widetilde{\uf}$ of $D$,
whose components will be denoted with the tilde.
The dynamics of this problem can be written in metriplectic form \eqref{eq:pbm}-\eqref{eq:degCond} where
the entropy of the system is linear and given by $S=S_1+S_2$, while the Hamiltonian is the total energy
\begin{equation*}
H(\widetilde{\uf})=\frac12\left(\|\widetilde{\bb{p}}_1\|_2^2
+\|\widetilde{\bb{p}}_2\|_2^2\right)
+e_1(\lambda_1(\widetilde{\uf}),\widetilde{S}_1)
+e_2(\lambda_2(\widetilde{\uf}),\widetilde{S}_2),
\qquad\forall\,\widetilde{\uf}\in D,
\end{equation*}
where 
$\lambda_1(\widetilde{\uf}):=\norm{\widetilde{\bb q}_1}_2$, and $\lambda_2(\widetilde{\uf}):=\norm{\widetilde{\bb q}_2-\widetilde{\bb q}_1}_2$ are the elongations of the springs, and
\begin{equation*}
    e_j(\lambda_j(\widetilde{\uf}),\widetilde{S}_j) :=\frac12\log^2\left({\lambda_j(\widetilde{\uf})}\right)+\log\left({\lambda_j(\widetilde{\uf})}\right)+ \exp\left( {\widetilde{S}_j- \log(\lambda_j(\widetilde{\uf}))}  \right)-1,\quad j\in\{1,2\}.
\end{equation*}
Moreover, the Poisson and semi-metric operators have the following expressions
\begin{equation*}
    J=\begin{bmatrix}
        0_{2R\times 2R} & I_{2R} &0_{2R\times 2}\\
        -I_{2R} & 0_{2R\times 2R} & 0_{2R\times 2} \\
        0_{2\times 2R} & 0_{2\times 2R} & 0_{2\times 2}\\
    \end{bmatrix},\,
    G(\widetilde{\uf})=\begin{bmatrix}
        0_{4R\times4R} & 0_{4R\times1} & 0_{4R\times1}\\
        0_{1\times 4R} & T_2(\widetilde{\uf})/T_1(\widetilde{\uf}) & -1\\
        0_{1\times 4R} & -1 & T_1(\widetilde{\uf})/T_2(\widetilde{\uf})
    \end{bmatrix},
\end{equation*}
where, for any $\widetilde{\uf}\in D$,
\begin{equation*}
    T_j(\widetilde{\uf}):=\frac{\partial e_j}{\partial S_j}(\lambda_j(\widetilde{\uf}), \widetilde{S}_j)= \exp\left({\widetilde{S}_j-\log(\lambda_j(\widetilde{\uf}))}\right),\qquad j\in\{1,2\}.
\end{equation*}

It can be easily verified that
the degeneracy conditions \eqref{eq:degCond} are satisfied,
and therefore the total energy is conserved, and the total entropy is
non-decreasing,
since $T_j(\widetilde{\uf})>0$, $j\in\{1,2\}$, for all $\widetilde{\uf}\in D$.

\paragraph{The block diagonal reduced order model (BD-ROM)}
The two thermo-spring system has the structure discussed in \Cref{sec:block}, with $J$ and $G$ of the form \eqref{eq:JGstruct} with $N_1=4R$ and $N_2=2$, $G_1=0_{4R\times 4R}$, and
\begin{equation*}
    J_1=\begin{bmatrix}
        0_{2R\times 2R} & I_{2R} \\
        -I_{2R} & 0_{2R\times 2R} 
    \end{bmatrix}, \qquad
    G_2(\widetilde{\uf})=\begin{bmatrix}
            T_2(\widetilde{\uf})/T_1(\widetilde{\uf}) & -1\\
            -1 & T_1(\widetilde{\uf})/T_2(\widetilde{\uf})
        \end{bmatrix}.
\end{equation*}
In the following tests, we consider, for this model, a reduced basis of the form \eqref{eq:Ublock} with $U_2=I_2$, since the last two variables correspond to the entropy of the system, while $U_1\in\R^{N_1\times n_1}$ is obtained by applying POD to the first $4R$ components of the shifted snapshots.
To test the convergence properties of the method, in the numerical experiments, we consider fixed values of $n$. For the block-diagonal model this corresponds to fixing $n_1=n-2$, given $n$.

\paragraph{The oblique-projection reduced order model (OP-ROM)}
To derive the OP-ROM in \eqref{eq:OP-ROM}, we set $\ccH=\ccS=2^{-1/2}\nabla S$. Clearly, $\bS=\sqrt{2}$,
and $\bH(\widetilde{\uf})=2^{-1/2}(T_1(\widetilde{\uf})+T_2(\widetilde{\uf}))>0$,
for all $\widetilde{\uf}\in D$.
This implies that $\DmuS=\R^N$ and $\DmuH=D$.
%
For this model, a reduced space of dimension $n$ has been constructed via POD of the snapshot matrix truncated to $n-1$, and the reduced basis has been augmented with the vector $\ccH=\ccS\in\R^N$ and then orthonormalized.

Finally, the energetically consistent reduced order model (EC-ROM) has been implemented with  $k_0=k_1=N$ as suggested by the authors in \cite[Theorem 3.4]{gruber2023energetically}.

\begin{remark}\label{rmk:thermo_springs_hr}
	No hyper-reduction has been introduced for this test case; indeed, despite the nonlinearity of the problem, it is possible to make all the computations needed for the reduced order models independent of the full dimension $N$, after suitable parameter-dependent precomputations.
    The main observation is that the nonlinear terms depend on the full state only through the scalar quantities $\lambda_j$, and the entropy variables $S_j$, $j\in\{1,2\}$. A simple computation, which is carried out in detail in \Cref{sec:thermospring-details}, shows that both $\lambda_1(U\bw+\uf_0(\mu)),\lambda_2(U\bw+\uf_0(\mu))$ can be evaluated in $O(n^2)$ operations per instance of $\bw\in\R^n$, plus $O(nN)$ operations for each test parameter.  
    The previous argument also shows that
    the  EC-ROM of \cite{gruber2023energetically} can be evaluated at a comparable arithmetic complexity.
\end{remark}

\paragraph{Set up of the numerical tests}
For all the numerical experiments we set $R=1000$, so that the dimension of the full order system is
$N=4R+2=4002$.
Let $\bb{x}\in\R^R$ be a vector of $R$ uniformly spaced points in $[0,1]$. The initial conditions are, for $j\in\{1,2\}$,
\begin{align*}
    \bb q_j^0(\mu) &= \exp(\mu_j\,\bb x),
    & \mu_1 &\in[-0.2,5.2],\; \mu_2 \in[-0.1,0.1],\\
    \bb p_j^0(\mu) &= (1+\mu_{2+j}\,\bb x^2)^{-1},
    & \mu_3 &\in[-1,1],\; \mu_4 \in[-3,3],\\
    S_j^0(\mu) &= \mu_{4+j},
    & \mu_5 &\in[0.9,1.1],\, \mu_6 \in[0.1,0.3].
\end{align*}
In order to construct the snapshot matrix, we randomly select $|\p_{\train}|=25$ vectors of parameters according to a uniform distribution in the six-dimensional parameter domain, and such that $\uf_0(\mu_{\train})$ is well-defined and belongs to $D$.
For each training parameter, the full order problem is solved in the temporal interval $\mathcal T=[0,8]$ with the LSODA time integrator and snapshots are collected every $\Delta t=0.01$, so that $N_{\train}=800$.


For post-processing, we use the test parameter $\mu_{\mathrm{test}}=(10.45,\,3.12,\,1.25,\,4.08,\allowbreak\,1.54,\,2.22)$, which lies outside the parameter domain used for training. 
Similar results were obtained with different test parameters.

\paragraph{Numerical results}

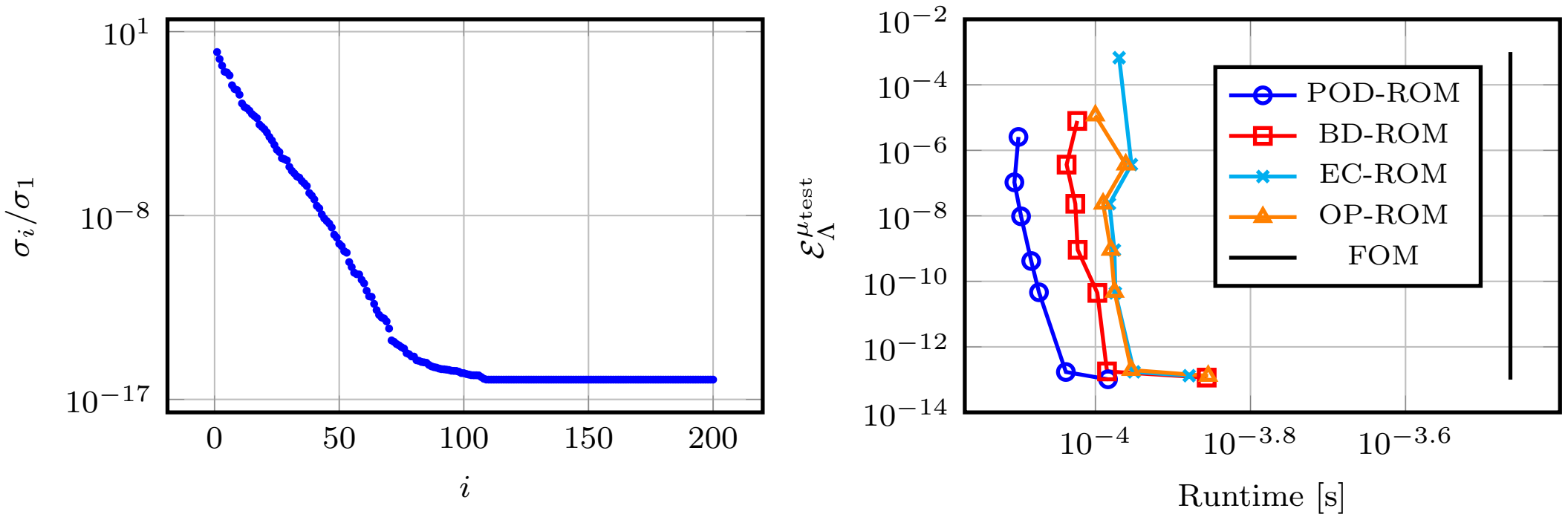
\begin{figure}[htbp]
  \centering
  \begin{tikzpicture}
    \begin{groupplot}[
      group style={group size=2 by 1, horizontal sep=1.8cm},
      height=3.5cm,
      grid=major,
      width=0.34\linewidth, scale only axis
    ]
        \nextgroupplot[
            ymode=log,
            xlabel={$i$},
            ylabel={$\sigma_i/\sigma_1$},
        ]
        \addplot+[
            blue, only marks, mark=*,
            restrict expr to domain={\coordindex}{0:199},
            every mark/.append style={scale=.25},
            unbounded coords=discard,
            filter discard warning=false,
        ]
        table[x=i, y=sigma_normalized]{\sigmatab};

      \nextgroupplot[
        xmode=log,
        ymode=log,
        xlabel={Runtime [s]},
        ylabel={$\errmutest$},
        legend columns=1,
        legend style={
           at={(0.87,0.6)},
           anchor=east,
        },
        ytick distance=100,
      ]

        \addplot+[blue, mark=o]
          table[
            x=POD_rhs_mean, y=POD_state_error,
            row predicate/.code={\ifnum\pgfplotstablerow=0\pgfplotstableuserowfalse\fi}
          ]{\errtab};
        \addlegendentry{POD-ROM}
        \addplot+[red, mark=square]
          table[
            x=BD_rhs_mean, y=BD_state_error,
            row predicate/.code={\ifnum\pgfplotstablerow=0\pgfplotstableuserowfalse\fi}
          ]{\errtab};
        \addlegendentry{BD-ROM}

        \addplot+[cyan, mark=x]
          table[
            x=SP_rhs_mean, y=SP_state_error,
            row predicate/.code={\ifnum\pgfplotstablerow=0\pgfplotstableuserowfalse\fi}
          ]{\errtab};
        \addlegendentry{EC-ROM}

        \addplot+[orange, mark=triangle]
          table[
            x=OP_rhs_mean, y=OP_state_error,
            row predicate/.code={\ifnum\pgfplotstablerow=0\pgfplotstableuserowfalse\fi}
          ]{\errtab};
        \addlegendentry{OP-ROM}
        
        \pgfplotstablegetelem{0}{FOM_rhs_mean}\of{\errtab}
        \edef\fomrhstime{\pgfplotsretval}

        \addplot+[
          black,
          solid,
          mark=none
        ]
        coordinates {
          (\fomrhstime,1e-13)
          (\fomrhstime,1e-3)
        };
        \addlegendentry{FOM}
    \end{groupplot}
  \end{tikzpicture}
  \caption{Two thermo-spring system. Left: Decay of the first 200 normalized singular values of the snapshot matrix $\SM$. Right: Error $\errmutest$ in \eqref{eq:relerr} for different models vs. the mean computational time per right-hand-side evaluation. Each mark corresponds to a different value of the dimension $n$ of the reduced order model, with $n\in\{20, 30, 40, 50, 60, 80, 100\}$.}
  \label{fig:thermospring_singval_error}
\end{figure}

In \Cref{fig:thermospring_singval_error} (right) we compare the performances of various methods in terms of accuracy and computational time. It can be observed that, in the temporal interval $[0,8]$ where snapshots are collected, all methods perform similarly in terms of overall accuracy and consistently with the decay of the singular values of the snapshot matrix \eqref{eq:snap-matrix}, as shown in \Cref{fig:thermospring_singval_error} on the left.

Concerning the computational runtime, 
\Cref{fig:thermospring_singval_error} (right) displays the mean time required for a single evaluation of the right-hand sides of the corresponding systems, while the total computational runtimes 
are reported in \Cref{table:thermospring_times}.
All methods provide a significant speedup compared to evaluating the full order model. As expected, the EC-ROM and OP-ROM are slightly more expensive as the operators involved are approximated via functions that depend on the state.


\pgfplotstableread[col sep=space,header=true]
{Numerical_results/metriplectic_error_vs_n.csv}\errtab

\pgfplotstableread[col sep=space,header=true]
{Numerical_results/metriplectic_time_n20_plot.csv}\timetab

\pgfplotstableread[col sep=space]
{Numerical_results/singular_values_pod.csv}\sigmatab

\pgfplotstableread[col sep=space,header=true]
{Numerical_results/metriplectic_error_vs_n.csv}\thermospringtimes

\newcommand{\ThermoSpringValue}[3][fixed,fixed zerofill,precision=4]{%
  \pgfplotstablegetelem{#2}{#3}\of{\thermospringtimes}%
  \pgfmathprintnumber[#1]{\pgfplotsretval}%
}

\newcommand{\ThermoSpringRow}[1]{%
  \ThermoSpringValue[fixed,precision=0]{#1}{n}
  & \ThermoSpringValue{#1}{POD_runtime}
  & \ThermoSpringValue{#1}{BD_runtime}
  & \ThermoSpringValue{#1}{SP_runtime}
  & \ThermoSpringValue{#1}{OP_runtime}\\%
}


\begin{table}[htbp]
\footnotesize
\caption[Computational times for the FOM and ROMs]{
Two thermo-spring system. Computational times required to solve the reduced order model on one test parameter for various methods and different reduced space dimension $n$. The final time considered is $T=8$. The runtime of the full order model is $\ThermoSpringValue{0}{FOM_runtime}\,\mathrm{s}$.}
\label{table:thermospring_times}
\setlength{\tabcolsep}{7pt}
\renewcommand{\arraystretch}{1.1}

\begin{center}
\begin{tabular}{|c|c|c|c|c|}
\hline
$n$ & POD-ROM & BD-ROM & EC-ROM & OP-ROM\\
\hline
\ThermoSpringRow{0}
\ThermoSpringRow{1}
\ThermoSpringRow{2}
\ThermoSpringRow{3}
\ThermoSpringRow{4}
\ThermoSpringRow{5}
\ThermoSpringRow{6}
\hline
\end{tabular}
\end{center}
\end{table}

\begin{figure}[htbp]
  \centering
  \tikzsetnextfilename{thermospring_hamiltonian_study}
  \begin{tikzpicture}
    \begin{groupplot}[
      group style={group size=2 by 1, horizontal sep=1.8cm},
      height=3.5cm,
      width=0.34\linewidth,
      scale only axis,
      ymode=log,
      grid=major,
    ]

      \nextgroupplot[
        xlabel={$n$},
        ylabel={$\errmutestH$},
        legend columns=1,
        legend pos=north east
      ]
        \addplot+[blue,mark=o] table[x=n,y=POD_energy_error]{\errtab};
        \addlegendentry{POD-ROM}
        \addplot+[red,mark=square] table[x=n,y=BD_energy_error]{\errtab};
        \addlegendentry{BD-ROM}
        \addplot+[cyan,mark=x] table[x=n,y=SP_energy_error]{\errtab};
        \addlegendentry{EC-ROM}
        \addplot+[orange,mark=triangle] table[x=n,y=OP_energy_error]{\errtab};
        \addlegendentry{OP-ROM}

      \nextgroupplot[
        xlabel={Time instant $t_{\idx}$},
        ylabel={$\errmutestH(t_{\idx})$},
        legend columns=1,
        legend style={
           at={(0.97,0.6)},
           anchor=east,
        },
      ]
        \addplot+[blue,dashed,mark=none] table[x=t,y=POD_energy_error]{\timetab};
        \addlegendentry{POD-ROM}
        \addplot+[red,densely dashed,mark=none] table[x=t,y=BD_energy_error]{\timetab};
        \addlegendentry{BD-ROM}
        \addplot+[cyan,dotted,mark=none] table[x=t,y=SP_energy_error]{\timetab};
        \addlegendentry{EC-ROM}
        \addplot+[orange,dashdotted,mark=none] table[x=t,y=OP_energy_error]{\timetab};
        \addlegendentry{OP-ROM}

    \end{groupplot}
  \end{tikzpicture}
  \caption{Two thermo-spring system.
  Left: Error $\errmutestH$ in \eqref{eq:relerrF} between the Hamiltonian
  evaluated at the full order solution and the Hamiltonian evaluated at the
  reconstructed reduced solution vs. the reduced space dimension $n$.
  Right: evolution of the relative error $\errmutestH(t_{\idx})$ in
  \eqref{eq:relerrF}, for $n=20$ fixed.}
  \label{fig:thermospring_hamiltonian_study}
\end{figure}
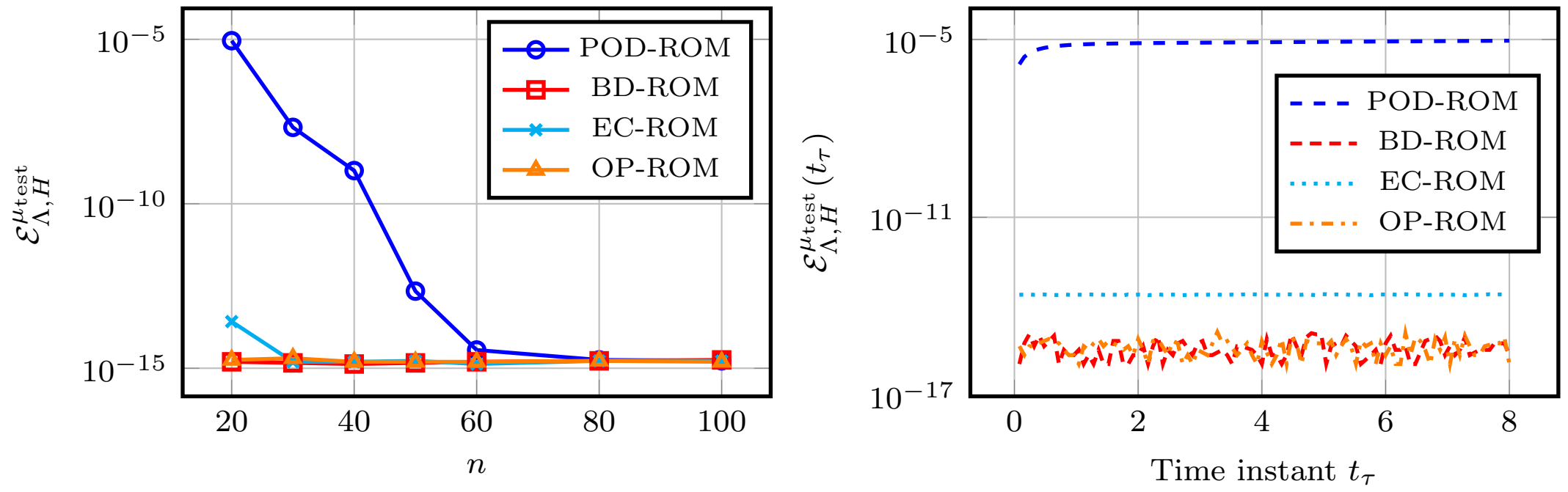

In \Cref{fig:thermospring_hamiltonian_study} we study the error in the Hamiltonian approximation.
First, we observe that the error in the Hamiltonian evaluated at the approximate solution of the structure-preserving models -- BD-ROM, OP-ROM and EC-ROM -- reaches the tolerance of the solver.
The Hamiltonian evaluated at the POD-ROM solution is not preserved, although the error with the Hamiltonian at the full order solution decreases with $n$. This can be associated with the inclusion of the snapshots of the Hamiltonian gradient in the construction of the reduced basis spaces.

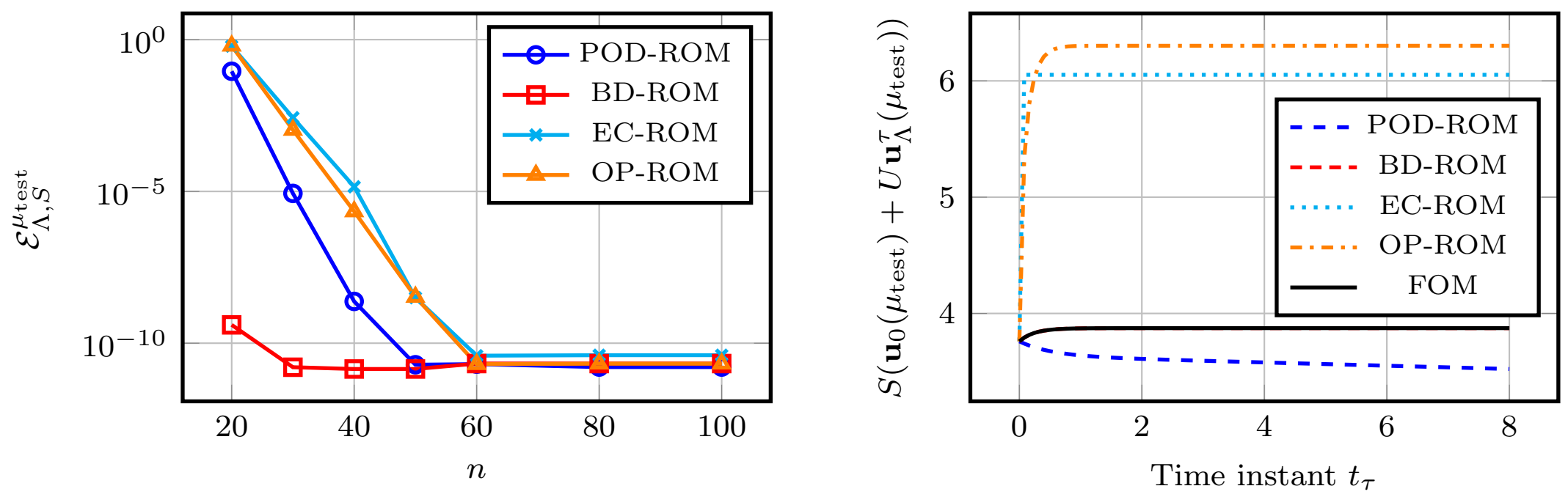
\begin{figure}[htbp]
  \centering
  \tikzsetnextfilename{thermospring_entropy_study}
  \begin{tikzpicture}
    \begin{groupplot}[
      group style={group size=2 by 1, horizontal sep=1.8cm},
      height=3.5cm,
      width=0.34\linewidth,
      scale only axis,
      grid=major,
    ]

      \nextgroupplot[
        ymode=log,
        xlabel={$n$},
        ylabel={$\errmutestS$},
        legend columns=1,
        legend pos=north east,
      ]
        \addplot[blue,mark=o] table[x=n,y=POD_entropy_error]{\errtab};
        \addlegendentry{POD-ROM}
        \addplot[red,mark=square] table[x=n,y=BD_entropy_error]{\errtab};
        \addlegendentry{BD-ROM}
        \addplot[cyan,mark=x] table[x=n,y=SP_entropy_error]{\errtab};
        \addlegendentry{EC-ROM}
        \addplot[orange,mark=triangle] table[x=n,y=OP_entropy_error]{\errtab};
        \addlegendentry{OP-ROM}

      \nextgroupplot[
        ymode=normal,
        xlabel={Time instant $t_{\idx}$},
        ylabel={$S(\uf_0(\mu_\text{test})+U\uf_\Lambda^{\idx}(\mu_\text{test}))$},
        legend columns=1,
        legend style={
            at={(0.97,0.5)},
            anchor=east,
        },
      ]
        \addplot[blue,dashed,mark=none] table[x=t,y=POD_entropy]{\timetab};
        \addlegendentry{POD-ROM}
        \addplot[red,densely dashed,mark=none] table[x=t,y=BD_entropy]{\timetab};
        \addlegendentry{BD-ROM}
        \addplot[cyan,dotted,mark=none] table[x=t,y=SP_entropy]{\timetab};
        \addlegendentry{EC-ROM}
        \addplot[orange,dashdotted,mark=none] table[x=t,y=OP_entropy]{\timetab};
        \addlegendentry{OP-ROM}
        \addplot[black,solid,mark=none] table[x=t,y=FOM_entropy]{\timetab};
        \addlegendentry{FOM}
    \end{groupplot}
  \end{tikzpicture}
  \caption{Two thermo-spring system. Left: error $\errmutestS$ in \eqref{eq:relerrF} between the entropy evaluated at the full order and reconstructed solutions as a function of the reduced dimension $n$. Right: evolution of the entropy $S$ at the reconstructed state for $n=20$. The FOM curve displays $S(\uf^\idx(\mu_\text{test}))$.}
  \label{fig:thermospring_entropy_study}
\end{figure}

\newpage
In \Cref{fig:thermospring_entropy_study} we compare the approximation of the entropy. As expected, the entropy error of the reduced order models decreases as the dimension $n$ of the reduced basis space increases.
BD-ROM achieves an approximation error close to the solver tolerance for all tested values of $n$. Choosing $U_2=I_2$ corresponds to not reducing the entropy variable and guarantees preservation of the entropies as Casimirs of the reversible dynamics. This highlights the benefits of exploiting the specific structure of the Poisson and semi-metric operators when constructing the reduced basis, as done by BD-ROM.

We also observe (\Cref{fig:thermospring_entropy_study} on the right) that the POD-ROM does not preserve the metriplectic structure and this result in an unphysical entropy decrease.

\subsection{Damped thermoelastic rod}\label{sec:rod}
As a second test case we consider the dynamical system stemming from the semi-discretization of a damped thermoelastic rod model, as presented in \cite[Section 5.2]{gruber2023energetically}. The problem describes a one dimensional elastic rod 
which evolves as a damped Hamiltonian system with friction. For $N=2R+1$, the full state vector is given by
$\uf =\begin{bmatrix} \bb q^\top & \bb p^\top & S\end{bmatrix}^\top\in\R^N$, with $\bb q,\bb p\in\R^R$, and $S\in\R$.
Let us consider a generic state vector $\widetilde{\uf}\in\R^N$ whose components will be denoted with a tilde.
The system admits a metriplectic formulation \eqref{eq:pbm}
with
\begin{equation*}
    J=\begin{bmatrix}
        0_{R\times R} & I_R &0_{R\times 1}\\
        -I_R & 0_{R\times R} & 0_{R\times 1}\\
        0_{1\times R} & 0_{1\times R} &0
    \end{bmatrix},\quad G(\widetilde{\uf})=\gamma\begin{bmatrix}
        0_{R\times R} & 0_{R\times R} & 0_{R\times 1}\\
        0_{R\times R} & I_R & -\widetilde{\bb{p}}\\
        0_{1\times R} & -\widetilde{\bb{p}}^\top & \lVert\widetilde{\bb{p}}\rVert^2
    \end{bmatrix},\qquad\forall\,\widetilde{\uf}\in \R^{N},
\end{equation*}
for a given constant $\gamma>0$,
and where the Hamiltonian and entropy are given by
\begin{equation*}
        H(\widetilde{\uf})=\frac{1}{2}\lVert \widetilde{\bb{p}} \rVert_2^2 + \widetilde{S} + V(\widetilde{\uf}),\qquad  S(\widetilde{\uf})=\widetilde{S}, \qquad\forall\,\widetilde{\uf}\in\R^N.
\end{equation*}
In this test we set $V(\widetilde{\uf})\coloneq \sum_{i=1}^R\cos(q_i)$, and consider a parametric initial condition $\uf_0(\mu)$ whose components are denoted with the subscript $0$.
A straightforward computation shows that the degeneracy conditions \eqref{eq:degCond} are satisfied.

Note that $J,G$ do not possess a structure compatible with the block-diagonal reduced basis construction of \Cref{sec:block}, and thus the corresponding reduced order model is not considered here.

\paragraph{The oblique-projection reduced and hyper-reduced  models}
To construct the reduced model \eqref{eq:OP-ROM}, we set $\ccH=\ccS=\bb e_{N}\in\R^N$.
Since the last component of $\nabla X$, for $X\in\{H,S\}$, is equal to one, this choice yields $\bX=1$ in \eqref{eq:beta}.
Similarly to the previous example, the evaluation of the reduced semi-metric operator $\Gr(\bw)$ can be carried out in $O(n^2)$ operations per instance of the reduced state $\bw$ plus an additional $O(nN)$ operations per parameter instance; we refer to the supplementary material, \Cref{sec:rod-details}, for the details. Concerning the 
reduced Hamiltonian gradient,
the gradient-preserving approach of Section \ref{sec:HSeim} can be applied by splitting the Hamiltonian into a quadratic part $\LH$ and the nonlinear part
$\NH({\widetilde{\uf}})=V({\widetilde{\bb{q}}})$.
In this case, it is easy to see that the nonlinear term $\NH$ can be decomposed as
$\NH({\widetilde{\uf}})=\cNh^\top \bb{h}({\widetilde{\bb{q}}}),$
where $\cNh=\bb{1}\in\mathbb{R}^R$ is the vector with all entries equal to $1$ and $h_i({\widetilde{\bb{q}}}):=\cos({q}_i)$, for $1\leq i\leq R$.
Then, the hyper-reduced gradient can be evaluated in $O(n\nHeim +n^2)$ online operations per reduced state instance, where $\nHeim$ is the EIM dimension used for the hyper-reduction of $\nabla H_r,$ plus an additional $O(nN)$ operations per parameter instance.
This is consistent with the findings from \Cref{oss:GPeim}. Indeed, in this case, every entry of $\bb h({\widetilde{\bb{q}}})$ depends on a single entry of ${\widetilde{\bb{q}}}$, and hence $s_{\bb{h}}=1$.
The evaluation of the linear part of the gradient has arithmetic complexity $O(n^2)$.

For the hyper-reduced model \eqref{eq:hrOP-ROM}, the vector $\cHhr$ can be chosen equal to $\bb e_N$ since
$\Heim(\widetilde{\uf})=\LH(\widetilde{\uf})+\bb a_H^\top\pb_H\bb h(\widetilde{\bb{q}})$
and the last component of its gradient is constant.
Thus, $\bHhr(\widetilde{\uf})=\bb e_N^\top \nabla \Heim(\bb v)=1\neq0$ for all $\widetilde{\uf}\in\R^N$.
For this model, a reduced space of dimension $n$ has been constructed via POD of the snapshot matrix truncated to $n-1$, and the reduced basis has been augmented with the vector $\ccHhr=\ccS=\bb e_N\in\R^N$ and then orthonormalized.
For the construction of the approximation space in the
gradient-preserving hyper-reduction of the Hamiltonian gradient, we employ the snapshot matrix
\begin{equation}\label{eq:snap-matrix-Heim}
\SM_{\text{H}} = \begin{bmatrix}
    \ldots\,D_{\bb h}(\bb q_0(\mu_j)+U_qU^\top\left(\uf(\theta_{ij})-\uf_0(\mu_j)\right))U_q\, \ldots 
\end{bmatrix}\in\R^{R\times nN_S},
\end{equation}
where $\theta_{ij}=(t_i,\eta_j)$ are subsamples from $\mathcal{T}_{\train}\times\p_{\train}$, and $U_q\in\R^{R\times n}$ is the matrix containing the first $R$ rows of $U\in\R^{N\times n}$.

\paragraph{Set up of the numerical tests}
In the numerical experiments, we set $R=2000$ so that the full order system has dimension $N=2R+1=4001$. Following the setup of the test in \cite{gruber2023energetically}, the value $\gamma=8$ was used, and the dependence on a parameter was introduced through the initial conditions
\begin{align*}
    \bb q_0(\mu)&=e^{\mu_1\, \bb x},& \mu_1\in[-0.2,5.2],\\
    \bb p_0(\mu)&= (1+\mu_2\,\bb x^2)^{-1},&\mu_2\in[-0.9,1],\\
    S_0(\mu)& =\mu_3,&\mu_3\in[1,3],
\end{align*}
where $\bb x\in\R^R$ is a vector of $R$ uniformly spaced points in $[0,1]$. In order to construct the snapshot matrix, we randomly select $|\p_{\train}|=25$ vectors of parameters according to a uniform distribution in the three-dimensional parameter domain. For each training parameter, the full order problem is solved with the LSODA time integrator on $\mathcal{T}=[0,8]$ and snapshots are collected every $\Delta t=0.02$, so that $N_{\train}=400$. 

For post-processing, we consider the test parameter
$\mu_{\mathrm{test}}=(0.65,-0.1,1.9)$, as in \cite{gruber2023energetically}, reduced dimensions
$n\in\{40,80,120,160\}$,
and EIM dimensions
$\nHeim\in\{40, 80, 150, 200, 400, 600\}$.

\paragraph{Numerical results}

\Cref{table:thermorod_long_term} provides a summary of the results over extended integration intervals. The results show that the approximation quality of the POD-ROM deteriorates outside of the training time window, and this is related with the failure to preserve the geometric structure of the problem.

\begin{figure}[htbp]
\centering
\tikzsetnextfilename{thermorod_hamiltonian_preservation}
\begin{tikzpicture}
\begin{groupplot}[
    group style={group size=2 by 1,horizontal sep=1.8cm},
    height = 3.5cm,
    width=0.34\linewidth,
    scale only axis,
    ymode=log,
    xlabel={Time instant $t_{\idx}$},
    xmin=-2.4,
    xmax=50.4,
    xtick={0,10,20,30,40,50},
    grid=major,
    unbounded coords=discard,
]

\nextgroupplot[
    ylabel={$\DmutestH(t_\tau)$},
    ymin=2e-16,
    ymax=1e-5,
    ytick={1e-15,1e-13,1e-11,1e-9,1e-7,1e-5},
    legend columns=1,
    legend style={
        at={(0.98,0.65)},
        anchor=east,
    },
]

\addplot[blue,dashed,mark=none]
table[x=t,y expr={max(\thisrow{POD-ROM},1e-18)}]{\roddrifttab};
\addlegendentry{POD-ROM}

\addplot[cyan,dotted,mark=none]
table[x=t,y expr={max(\thisrow{SP-ROM},1e-18)}]{\roddrifttab};
\addlegendentry{EC-ROM}

\addplot[orange,dashdotted,mark=none]
table[x=t,y expr={max(\thisrow{OP-ROM},1e-18)}]{\roddrifttab};
\addlegendentry{OP-ROM}

\addplot[black,solid,mark=none]
table[x=t,y expr={max(\thisrow{FOM},1e-18)}]{\roddrifttab};
\addlegendentry{FOM}

\nextgroupplot[
    ylabel={$\mathcal{D}^{\mu_\text{test}}_{\text{OP-HROM},H}(t_{\idx})$},
    ymin=1e-16,
    ymax=1.4,
    ytick={1e-16,1e-13,1e-10,1e-7,1e-4,1e-1},
    legend columns=1,
    legend style={
        at={(0.7,0.35)},
        anchor=east,
    },
]

\addplot[blue,solid,mark=none]
table[x=t,y expr={max(\thisrow{OP-HROM_m40},1e-18)}]{\roddrifttab};
\addlegendentry{\(\nHeim=40\)}

\addplot[red,dashed,mark=none]
table[x=t,y expr={max(\thisrow{OP-HROM_m80},1e-18)}]{\roddrifttab};
\addlegendentry{\(\nHeim=80\)}

\addplot[teal!70!black,densely dashed,mark=none]
table[x=t,y expr={max(\thisrow{OP-HROM_m150},1e-18)}]{\roddrifttab};
\addlegendentry{\(\nHeim=150\)}

\addplot[violet,dotted,mark=none]
table[x=t,y expr={max(\thisrow{OP-HROM_m200},1e-18)}]{\roddrifttab};
\addlegendentry{\(\nHeim=200\)}

\addplot[orange!80!black,dashdotted,mark=none]
table[x=t,y expr={max(\thisrow{OP-HROM_m400},1e-18)}]{\roddrifttab};
\addlegendentry{\(\nHeim=400\)}

\addplot[brown,dashdotdotted,mark=none]
table[x=t,y expr={max(\thisrow{OP-HROM_m600},1e-18)}]{\roddrifttab};
\addlegendentry{\(\nHeim=600\)}
\end{groupplot}
\end{tikzpicture}
\caption{Damped thermoelastic rod. Temporal interval $\mathcal{T}=[0,48]$. Left: evolution of the error $\DmutestH$ in \eqref{eq:presrelerr} between the Hamiltonian evaluated at the solution of the various models and the Hamiltonian at the initial condition of the full order model. The reduced space dimension is set to $n=80$. Right: evolution of $\DmutestH$ for the OP-HROM for different values of 
$\nHeim$.}
\label{fig:thermorod_hamiltonian_preservation}
\end{figure}
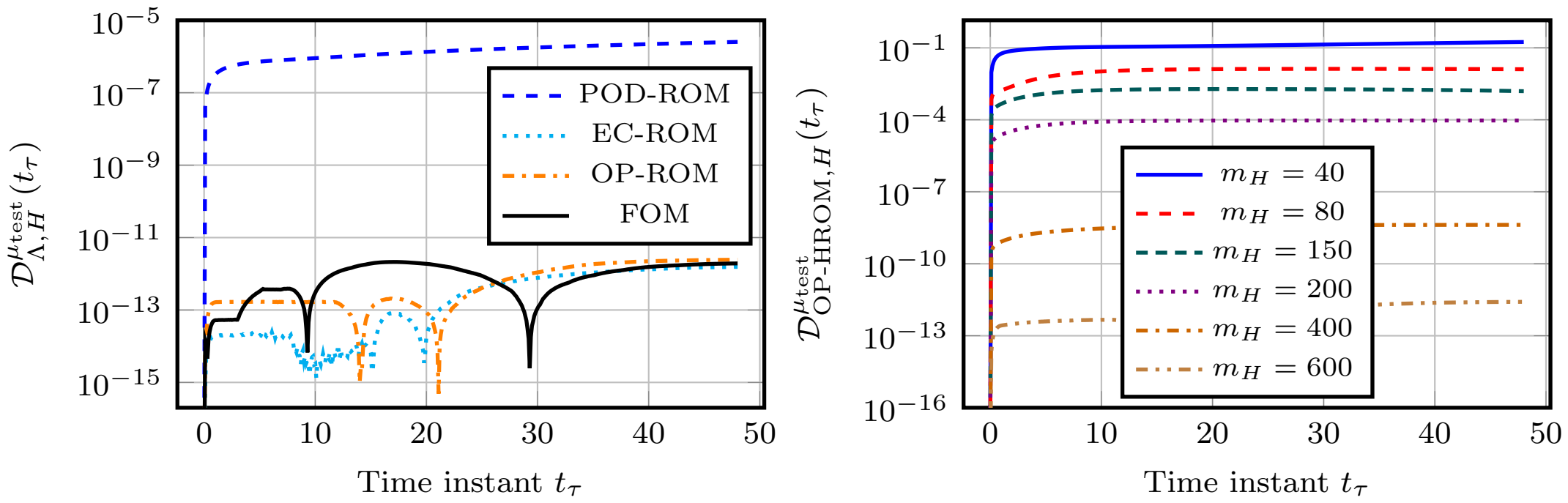

We remark that the hyper-reduced model is metriplectic with Hamiltonian given by an approximation of the original Hamiltonian function, and this can be numerically observed in the conservation error reported in the fourth column of \Cref{table:thermorod_long_term}.
However, as the dimension $\nHeim$ of the hyper-reduced space increases this error decreases, as reported in \Cref{fig:thermorod_hamiltonian_preservation} on the right.

Finally, the computational times reported in the last column of \Cref{table:thermorod_long_term} highlight the effectiveness of reduction and the further runtime saving associated with hyper-reduction. In particular, OP-HROM takes less than a third of the runtime required by the corresponding reduced model without hyper-reduction, while retaining comparable accuracy.

\pgfplotstableread[col sep=comma,header=true]
{Numerical_results/thermoelastic_rod_n40_m150_table.csv}\rodtab

\newcommand{\RodTableCell}[2]{%
  \pgfplotstablegetelem{#1}{#2}\of{\rodtab}%
  \pgfplotsretval%
}
\newcommand{\RodTableValue}[3][]{%
  \pgfplotstablegetelem{#2}{#3}\of{\rodtab}%
  \pgfmathprintnumber[#1]{\pgfplotsretval}%
}

\newcommand{\RodTableResults}[1]{%
  \RodTableCell{#1}{method}
  & \RodTableValue[sci,precision=3]{#1}{state_error}
  & \RodTableValue[sci,precision=3]{#1}{max_relative_hamiltonian_drift}
  & \RodTableValue[fixed,fixed zerofill,precision=4]{#1}{runtime_s}%
}

\newcommand{\RodTableRows}[5]{%
  \multirow{5}{*}{\RodTableCell{#1}{T}}
   & \RodTableCell{#1}{method} & --
   & \RodTableValue[sci,precision=3]{#1}{max_relative_hamiltonian_drift}
   & \RodTableValue[fixed,fixed zerofill,precision=4]{#1}{runtime_s} \\
  \cline{2-5}
   & \RodTableResults{#2} \\
   & \RodTableResults{#3} \\
   & \RodTableResults{#4} \\
   & \RodTableResults{#5} \\
}

\begin{table}[htbp]
\footnotesize
\caption{Damped thermoelastic rod.
Comparison of different methods over longer temporal intervals $[0,T]$ in terms of state error, Hamiltonian preservation error, and total runtime.
The reduced dimension is set to $n=40$ and for the OP-HROM the
hyper-reduced dimension is set to $\nHeim=150$.}
\label{table:thermorod_long_term}
\begin{center}
\begin{tabular}{|c|c|c|c|c|}
\hline
$T$ & Method $\Lambda$ & $\errmutest$ & $\max_\idx \DmutestH(t_\idx)$ & Runtime [s] \\
\hline
\RodTableRows{0}{1}{2}{3}{6}
\hline
\RodTableRows{7}{8}{9}{10}{13}
\hline
\RodTableRows{14}{15}{16}{17}{20}
\hline
\end{tabular}
\end{center}
\end{table}

In \Cref{fig:thermorod_hamiltonian_preservation} we study the Hamiltonian conservation for the different methods.
The relative error in the Hamiltonian conservation (plot on the left) reaches the numerical integration errors only for the structure-preserving models (EC-ROM and OP-ROM).
The right panel shows the Hamiltonian conservation error along the trajectories of the hyper-reduced model, for different values of the EIM dimension $\nHeim$. Although, by \Cref{lem:hrproperties}, OP-HROM preserves the hyper-reduced Hamiltonian rather than the original one, convergence is observed as $\nHeim$ increases.

\pgfplotstableread[col sep=comma,header=true]
{Numerical_results/Numerical_results_thermoelastic_rod_N2000_T48/hamiltonian_drift_n80_plot.csv}\roddrifttab

\Cref{fig:thermorod_entropy} reports the behavior of the entropy for the different models.
In this test case, the OP-ROM provides the most accurate approximation of the entropy (left plot). The curve associated with the OP-HROM shows convergence as $\nHeim$ increases, as expected.
The right panel shows that all the tested ROMs are able to correctly represent the growth of the entropy functional.

\pgfplotstableread[col sep=comma,header=true]
{Numerical_results/Numerical_results_thermoelastic_rod_N2000_T48/entropy_n80_plot.csv}\rodentropytab

\pgfplotstableread[col sep=space,header=true]
{Numerical_results/Numerical_results_thermoelastic_rod_N2000_T48/entropy_hamiltonian_vs_m.csv}\rodentropyerrtab

\pgfplotstablegetelem{18}{E_S_r_f}\of{\rodentropyerrtab}
\edef\RodEntropyPodErr{\pgfplotsretval}
\pgfplotstablegetelem{24}{E_S_r_f}\of{\rodentropyerrtab}
\edef\RodEntropyEcErr{\pgfplotsretval}
\pgfplotstablegetelem{30}{E_S_r_f}\of{\rodentropyerrtab}
\edef\RodEntropyOpErr{\pgfplotsretval}

\begin{figure}[htbp]
  \centering
  \tikzsetnextfilename{thermorod_entropy}
  \begin{tikzpicture}
    \begin{groupplot}[
      group style={
        group name=RodEntropyGroup,
        group size=2 by 1,
        horizontal sep=1.8cm,
      },
      grid=major,
      height = 3.5cm,
      width=0.34\linewidth,
      scale only axis,
      legend columns=3,
      legend style={
        /tikz/every even column/.append style={column sep=0.4em},
      },
    ]

      \nextgroupplot[
        ymode=log,
        xlabel={$\nHeim$},
        ylabel={$\errmutestS$},
        xmin=-15,
        xmax=655,
        xtick={0,200,400,600},
        unbounded coords=discard,
        filter discard warning=false,
        legend pos= north east,
        legend columns=1
      ]
        \addplot[orange,mark=triangle,
          restrict expr to domain={\coordindex}{30:35}]
          table[x=m,y=E_S_hr_f]{\rodentropyerrtab};
        \addlegendentry{OP-HROM}

        \addplot[blue,dashed,mark=none]
          coordinates {(40,\RodEntropyPodErr) (600,\RodEntropyPodErr)};
        \addlegendentry{POD-ROM}

        \addplot[cyan,dotted,mark=none]
          coordinates {(40,\RodEntropyEcErr) (600,\RodEntropyEcErr)};
        \addlegendentry{EC-ROM}

        \addplot[orange,dashdotted,mark=none]
          coordinates {(40,\RodEntropyOpErr) (600,\RodEntropyOpErr)};
        \addlegendentry{OP-ROM}

      \nextgroupplot[
        ymode=normal,
        xlabel={Time instant $t_{\idx}$},
        ylabel={$S(\uf_0(\mu_\text{test})+U\uf_\Lambda^{\idx}(\mu_\text{test}))$},
        xmin=-2.4,
        xmax=50.4,
        xtick={0,10,20,30,40,50},
        legend pos= south east,
        legend columns=1
      ]
        \addplot[black,solid,mark=none] table[x=t,y=FOM]{\rodentropytab};
        \addlegendentry{FOM}
        
        \addplot[blue,dashed,mark=none] table[x=t,y=POD-ROM]{\rodentropytab};
        \addlegendentry{POD-ROM}

        \addplot[cyan,dotted,mark=none] table[x=t,y=SP-ROM]{\rodentropytab};
        \addlegendentry{EC-ROM}

        \addplot[orange,dashdotted,mark=none] table[x=t,y=OP-ROM]{\rodentropytab};
        \addlegendentry{OP-ROM}

    \end{groupplot}
  \end{tikzpicture}

\caption{Damped thermoelastic rod. Temporal interval $\mathcal{T}=[0,48]$.
Left: Error $\errmutestS$ in \eqref{eq:relerrF} between the entropy evaluated at the full order solution and the entropy evaluated at the reconstructed reduced or hyper-reduced solution vs. the hyper-reduced space dimension $\nHeim$.
Right: evolution of the entropy for $n=80$ fixed. The FOM curve displays $S(\uf^\idx(\mu_\text{test}))$.
}
\label{fig:thermorod_entropy}
\end{figure}
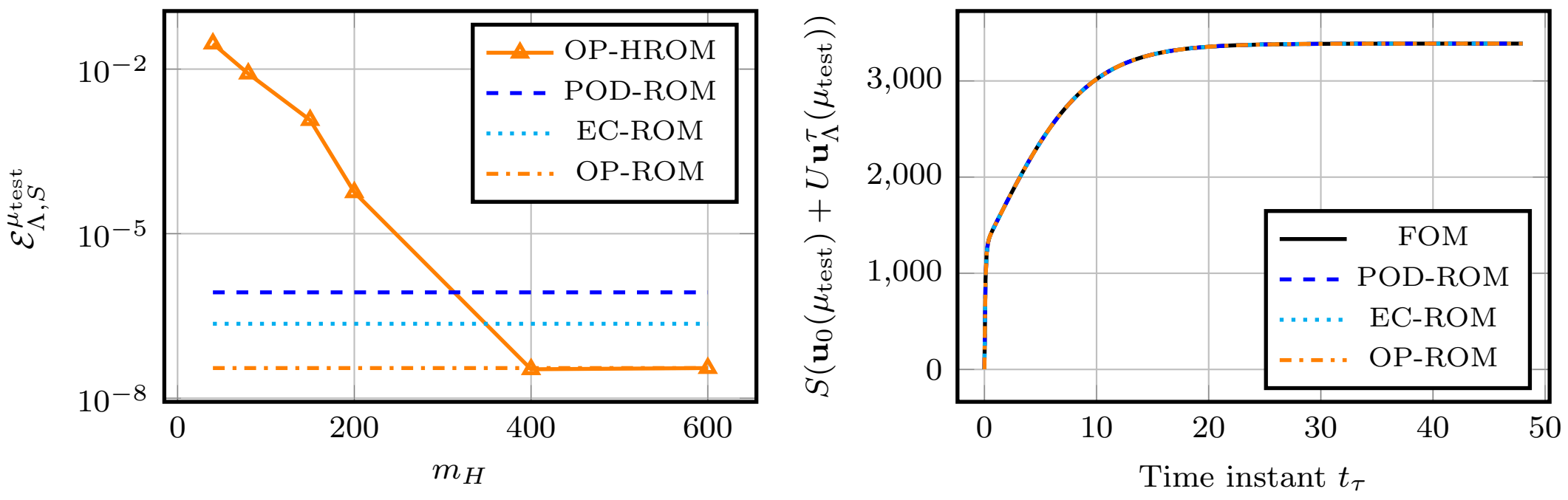

The particular structure of this example allows to efficiently hyper-reduce the EC-ROM of \cite{gruber2023energetically} in a structure-preserving way via the gradient-preserving approach of \cite{pagliantini2023gradient}. The corresponding model is dubbed EC-HROM.
We refer to the supplementary materials for a detailed discussion on the construction.
Since the corresponding reduced bases differ, EC-HROM and OP-HROM are constructed from different snapshot matrices \eqref{eq:snap-matrix-Heim} and thus
use different EIM bases and sets of interpolation indices.
We point out that, for general problems, the EC-ROM requires a factorization of the semi-metric tensor $\GM$, whose complexity scales cubically with $N$.

\Cref{table:thermorod_hrsummary} compares the relative state errors $\errmutest$, defined in \eqref{eq:relerr}, and the computational times. 
For sufficiently large hyper-reduction spaces, both EC-HROM and OP-HROM recover the accuracy of the corresponding ROM while reducing the computational runtime.
Moreover, while their accuracy is comparable over the tested reduced and hyper-reduced spaces, the OP-HROM generally requires less computational time than the EC-HROM.



\pgfplotstableread[col sep=space,header=true]
{Numerical_results/Numerical_results_thermoelastic_rod_N2000_T48/fom_timing_reference.csv}\rodsolvefom

\pgfplotstableread[col sep=space,header=true]
{Numerical_results/Numerical_results_thermoelastic_rod_N2000_T48/non_hyper_reduced_rom_reference.csv}\rodsolverom

\pgfplotstableread[col sep=space,header=true]
{Numerical_results/Numerical_results_thermoelastic_rod_N2000_T48/fixed_n_variable_m.csv}\rodsolvehr

\providecommand{\RodSolveTime}[2]{%
  \pgfplotstablegetelem{#2}{runtime_s}\of{#1}%
  \pgfmathprintnumber[fixed,fixed zerofill,precision=4]{\pgfplotsretval}%
}

\providecommand{\RodSolveDimension}[3]{%
  \pgfplotstablegetelem{#2}{#3}\of{#1}%
  \pgfmathprintnumber[fixed,precision=0]{\pgfplotsretval}%
}

\pgfplotstableread[col sep=space,header=true]
{Numerical_results/Numerical_results_thermoelastic_rod_N2000_T48/fom_timing_reference.csv}\rodsolvefom

\pgfplotstableread[col sep=space,header=true]
{Numerical_results/Numerical_results_thermoelastic_rod_N2000_T48/non_hyper_reduced_rom_reference.csv}\rodsolverom

\pgfplotstableread[col sep=space,header=true]
{Numerical_results/Numerical_results_thermoelastic_rod_N2000_T48/fixed_n_variable_m.csv}\rodsolvehr

\providecommand{\RodSolveTime}[2]{%
  \pgfplotstablegetelem{#2}{runtime_s}\of{#1}%
  \pgfmathprintnumber[fixed,fixed zerofill,precision=4]{\pgfplotsretval}%
}

\providecommand{\RodSolveDimension}[3]{%
  \pgfplotstablegetelem{#2}{#3}\of{#1}%
  \pgfmathprintnumber[fixed,precision=0]{\pgfplotsretval}%
}

\newcommand{\RodSolveError}[3]{%
  \pgfplotstablegetelem{#2}{#3}\of{#1}%
  \pgfmathprintnumber[sci,sci zerofill,precision=2]{\pgfplotsretval}%
}

\newcommand{\RodSolveResults}[4]{%
  \RodSolveError{#1}{#3}{#2} & \RodSolveTime{#1}{#3}
  & \RodSolveError{#1}{#4}{#2} & \RodSolveTime{#1}{#4}%
}

\begin{table}[htbp]
\footnotesize
\caption[]{Damped thermoelastic rod. Temporal interval $\mathcal{T}=[0,48]$. Comparison of the proposed methods and the one of \cite{gruber2023energetically} with and without hyper-reduction, in terms of state error \eqref{eq:relerr}, and total runtime. The reduced dimension is either $n=40$ or $n=80$ and different values of the hyper-reduced space dimension $\nHeim$ are considered. The FOM runtime is $\RodSolveTime{\rodsolvefom}{0}\,\mathrm{s}$.}\label{table:thermorod_hrsummary}
\begin{center}
\begin{tabular}{|c|c|c|c|c|c|c|}
\hline
\multicolumn{3}{|c|}{} & \multicolumn{2}{c|}{EC} & \multicolumn{2}{c|}{OP} \\
\hline
$n$ & Model & $\nHeim$
 & $\errmutest$ & Runtime [s]
 & $\errmutest$ & Runtime [s] \\
\hline
\multirow{7}{*}{\RodSolveDimension{\rodsolverom}{0}{n}} & ROM & --
 & \RodSolveResults{\rodsolverom}{E_r_f}{1}{2} \\
\cline{2-7}
 & \multirow{6}{*}{HROM} & \RodSolveDimension{\rodsolvehr}{0}{m}
 & \RodSolveResults{\rodsolvehr}{E_hr_f}{1}{2} \\
 & & \RodSolveDimension{\rodsolvehr}{3}{m}
 & \RodSolveResults{\rodsolvehr}{E_hr_f}{4}{5} \\
 & & \RodSolveDimension{\rodsolvehr}{6}{m}
 & \RodSolveResults{\rodsolvehr}{E_hr_f}{7}{8} \\
 & & \RodSolveDimension{\rodsolvehr}{9}{m}
 & \RodSolveResults{\rodsolvehr}{E_hr_f}{10}{11} \\
 & & \RodSolveDimension{\rodsolvehr}{12}{m}
 & \RodSolveResults{\rodsolvehr}{E_hr_f}{13}{14} \\
 & & \RodSolveDimension{\rodsolvehr}{15}{m}
 & \RodSolveResults{\rodsolvehr}{E_hr_f}{16}{17} \\
\hline
\multirow{7}{*}{\RodSolveDimension{\rodsolverom}{3}{n}} & ROM & --
 & \RodSolveResults{\rodsolverom}{E_r_f}{4}{5} \\
\cline{2-7}
 & \multirow{6}{*}{HROM} & \RodSolveDimension{\rodsolvehr}{18}{m}
 & \RodSolveResults{\rodsolvehr}{E_hr_f}{19}{20} \\
 & & \RodSolveDimension{\rodsolvehr}{21}{m}
 & \RodSolveResults{\rodsolvehr}{E_hr_f}{22}{23} \\
 & & \RodSolveDimension{\rodsolvehr}{24}{m}
 & \RodSolveResults{\rodsolvehr}{E_hr_f}{25}{26} \\
 & & \RodSolveDimension{\rodsolvehr}{27}{m}
 & \RodSolveResults{\rodsolvehr}{E_hr_f}{28}{29} \\
 & & \RodSolveDimension{\rodsolvehr}{30}{m}
 & \RodSolveResults{\rodsolvehr}{E_hr_f}{31}{32} \\
 & & \RodSolveDimension{\rodsolvehr}{33}{m}
 & \RodSolveResults{\rodsolvehr}{E_hr_f}{34}{35} \\
\hline
\end{tabular}    
\end{center}
\end{table}


\section{Conclusions}\label{sec:conclusions}

We developed two structure-preserving model order reduction strategies for metriplectic systems, using oblique projections compatible with hyper-reduction and block-diagonal bases adapted to the operator structure. We established preservation properties, thermodynamical consistency, and approximation error bounds, and demonstrated accuracy and computational gains on two benchmarks. Further investigation should target
more demanding tests that include nonlinear entropy, state-dependent Poisson operators, and larger dissipative blocks. Another direction is the development of time integrators that preserve the thermodynamic properties at the fully discrete level.

\section*{Acknowledgments}
The authors acknowledge the MIUR Excellence Department Project awarded to the Department of Mathematics,
University of Pisa, CUP I57G22000700001 and the INdAM-GNCS.
The work of A.B. has been supported by the European Union – NextGenerationEU, under the National Recovery and Resilience Plan (NRRP), through the MUR Project CUP I51J23000400007.
A.B. wishes to thank Damiano Lombardi (INRIA Paris), Michiel Hochstenbach (Eindhoven University of Technology), and Mark Peletier (Eindhoven University of Technology) for fruitful discussions on topics related to this work.

\bibliographystyle{plain}
\bibliography{references}

@unpublished{lombardi2024conformal,
      title={{Conformal variational discretisation of infinite dimensional Hamiltonian systems with gradient flow dissipation}}, 
      author={Damiano Lombardi and Cecilia Pagliantini},
      year={2025},
      note = {Preprint, arXiv:2412.06310},
}

@article{barrault2004empirical,
 author = {Barrault, Maxime and Maday, Yvon and Nguyen, Ngoc Cuong and Patera, Anthony T.},
 title = {An `empirical interpolation' method: {Application} to efficient reduced-basis discretization of partial differential equations},
 journal = {C. R. Math. Acad. Sci. Paris},
 issn = {1631-073X},
 volume = {339},
 number = {9},
 pages = {667--672},
 year = {2004}
}

@article{peng2016,
  title={{Symplectic model reduction of Hamiltonian systems}},
  author={Peng, Liqian and Mohseni, Kamran},
  journal={SIAM J. Sci. Comput.},
  volume={38},
  number={1},
  pages={A1--A27},
  year={2016},
  publisher={SIAM}
}

@article{chaturantabut2010nonlinear,
  title={Nonlinear model reduction via discrete empirical interpolation},
  author={Chaturantabut, Saifon and Sorensen, Danny C},
  journal={SIAM J. Sci. Comput.},
  volume={32},
  number={5},
  pages={2737--2764},
  year={2010},
  publisher={SIAM}
}

@article{negri2015efficient,
  title={Efficient model reduction of parametrized systems by matrix discrete empirical interpolation},
  author={Negri, Federico and Manzoni, Andrea and Amsallem, David},
  journal={J. Comput. Phys.},
  volume={303},
  pages={431--454},
  year={2015},
  publisher={Elsevier}
}

@article{pagliantini2023gradient,
  title={Gradient-preserving hyper-reduction of nonlinear dynamical systems via discrete empirical interpolation},
  author={Pagliantini, Cecilia and Vismara, Federico},
  journal={SIAM J. Sci. Comput.},
  volume={45},
  number={5},
  pages={A2725--A2754},
  year={2023},
  publisher={SIAM}
}

@article{gruber2023energetically,
    AUTHOR = {Gruber, Anthony and Gunzburger, Max and Ju, Lili and Wang,
              Zhu},
     TITLE = {Energetically consistent model reduction for metriplectic
              systems},
   JOURNAL = {Comput. Methods Appl. Mech. Engrg.},
    VOLUME = {404},
      YEAR = {2023},
     eid = {115709},
}

@article{afkham2017structure,
  title={{Structure preserving model reduction of parametric Hamiltonian systems}},
  author={Maboudi Afkham, Babak and Hesthaven, Jan S},
  journal={SIAM J. Sci. Comput.},
  volume={39},
  number={6},
  pages={A2616--A2644},
  year={2017},
  publisher={SIAM}
}

@article{hesthaven2021structure,
  title={{Structure-preserving reduced basis methods for Poisson systems}},
  author={Hesthaven, Jan S and Pagliantini, Cecilia},
  journal={Math. Comp.},
  volume={90},
  number={330},
  pages={1701--1740},
  year={2021}
}

@article{Ottinger2015,
  title = {Preservation of thermodynamic structure in model reduction},
  author = {{\"O}ttinger, Hans Christian},
  journal = {Phys. Rev. E},
  volume = {91},
  issue = {3},
  eid = {032147},
  year = {2015},
  month = {Mar},
  publisher = {American Physical Society}
}

@article{Hernandez2021,
title = {Deep learning of thermodynamics-aware reduced-order models from data},
journal = {Comput. Methods Appl. Mech. Engrg.},
volume = {379},
eid = {113763},
year = {2021},
author = {Quercus Hernandez and Alberto Badías and David González and Francisco Chinesta and Elías Cueto},
}

@article{morrison1986paradigm,
  title={{A paradigm for joined Hamiltonian and dissipative systems}},
  author={Morrison, Philip J},
  journal={Phys. D},
  volume={18},
  number={1-3},
  pages={410--419},
  year={1986},
  publisher={Elsevier}
}

@phdthesis{martin2016energy,
  title={Energy-entropy-momentum time integration methods for coupled smooth dissipative problems},
  author={Conde Mart{\'i}n, Sergio},
  address={Madrid, Spain},
  year={2016},
  school={Universidad Polit{\'e}cnica de Madrid}
}

@article{Zhang2022,
    AUTHOR = {Zhang, Zhen and Shin, Yeonjong and Karniadakis, George Em},
     TITLE = {{GFINN}s: {GENERIC} formalism informed neural networks for deterministic and stochastic dynamical systems},
   JOURNAL = {Philos. Trans. Roy. Soc. A},
    VOLUME = {380},
      YEAR = {2022},
    NUMBER = {2229},
     eid = {20210207},
}

@article{gong2017structure,
  title={Structure-preserving {G}alerkin {POD} reduced-order modeling of {Hamiltonian} systems},
  author={Gong, Yuezheng and Wang, Qi and Wang, Zhu},
  journal={Comput. Methods Appl. Mech. Engrg.},
  volume={315},
  pages={780--798},
  year={2017},
  publisher={Elsevier}
}

@article{buchfink2019symplectic,
  title={Symplectic model order reduction with non-orthonormal bases},
  author={Buchfink, Patrick and Bhatt, Ashish and Haasdonk, Bernard},
  journal={Math. Comput. Appl.},
  volume={24},
  number={2},
  eid={43},
  year={2019},
  publisher={MDPI}
}

@article{gugercin12,
title = {Structure-preserving tangential interpolation for model reduction of port-{H}amiltonian systems},
author = {Serkan Gugercin and Rostyslav V. Polyuga and Christopher Beattie and Arjan {van der Schaft}},
journal = {Automatica},
volume = {48},
number = {9},
pages = {1963-1974},
year = {2012},
issn = {0005-1098},
}

@article{polyuga12,
title = {Effort- and flow-constraint reduction methods for structure preserving model reduction of port-{H}amiltonian systems},
author = {Rostyslav V. Polyuga and Arjan J. {van der Schaft}},
journal = {Syst. Control Lett.},
fjournal = {Systems & Control Letters},
volume = {61},
number = {3},
pages = {412-421},
year = {2012},
issn = {0167-6911}
}

@incollection{hindmarsh1983odepack,
    AUTHOR = {Hindmarsh, Alan C.},
     TITLE = {O{DEPACK}, a systematized collection of {ODE} solvers},
 BOOKTITLE = {Scientific computing},
    SERIES = {IMACS Trans. Sci. Comput.},
    VOLUME = {I},
     PAGES = {55--64},
 PUBLISHER = {IMACS, New Brunswick, NJ},
      YEAR = {1983},
}

@article{petzold83,
author = {Petzold, Linda},
title = {Automatic Selection of Methods for Solving Stiff and Nonstiff Systems of Ordinary Differential Equations},
journal = {SIAM J. Sci. Statist. Comput.},
volume = {4},
number = {1},
pages = {136-148},
year = {1983},
}

@article{grmela1997dynamics,
  title = {{Dynamics and thermodynamics of complex fluids.  I. Development of a general formalism}},
  author = {Grmela, Miroslav and {\"O}ttinger, Hans Christian},
  journal = {Phys. Rev. E},
  volume = {56},
  issue = {6},
  pages = {6620--6632},
  year = {1997},
  publisher = {American Physical Society},
}

@inproceedings{gruber2025efficiently,
    author = {Gruber, Anthony and Lee, Kookjin and
              Lim, Haksoo and Park, Noseong and
              Trask, Nathaniel},
    title = {Efficiently Parameterized Neural
             Metriplectic Systems},
    booktitle = {Proceedings of the 13th International Conference on Learning Representations (ICLR)},
    pages = {85828--85852},
    year = {2025}
}

@article{GO97,
    AUTHOR = {\"{O}ttinger, Hans Christian and Grmela, Miroslav},
     TITLE = {Dynamics and thermodynamics of complex fluids. {II}.
              {I}llustrations of a general formalism},
   JOURNAL = {Phys. Rev. E (3)},
  FJOURNAL = {Physical Review E. Statistical, Nonlinear, and Soft Matter
              Physics},
    VOLUME = {56},
      YEAR = {1997},
    NUMBER = {6},
     PAGES = {6633--6655},
      ISSN = {1539-3755}
}

@article{EGB20,
	author = {Eldred, Christopher and Gay-Balmaz, Fran{\c{c}}ois},
	 title = {Single and double generator bracket formulations of multicomponent fluids with irreversible processes},
   journal = {J. Phys. A: Math. Theor.},
  fjournal = {Journal of Physics A: Mathematical and Theoretical},
 publisher = {{IOP} Publishing},
	volume = {53},
      year = {2020},
	number = {39},
 	 pages = {395701}}

@article{suzuki20,
year = {2020},
month = {mar},
publisher = {IOP Publishing},
volume = {52},
number = {1},
pages = {015516},
author = {Suzuki, Yukihito},
title = {A GENERIC formalism for Korteweg-type fluids: I. A comparison with classical theory},
fjournal = {Fluid Dynamics Research},
journal = {Fluid Dyn. Res.}
}

@article {DPZ12,
    AUTHOR = {Duong, Manh Hong and Peletier, Mark A. and Zimmer, Johannes},
     TITLE = {G{ENERIC} formalism of a {V}lasov-{F}okker-{P}lanck equation
              and connection to large-deviation principles},
   JOURNAL = {Nonlinearity},
  FJOURNAL = {Nonlinearity},
    VOLUME = {26},
      YEAR = {2013},
    NUMBER = {11},
     PAGES = {2951--2971},
      ISSN = {0951-7715,1361-6544}
}

\clearpage
\appendix
\section{Example: two thermo-spring system in entropy variables}\label{sec:thermospring-details}

In this part we provide further details about the entropy-based metriplectic formulation of the two thermo-spring system discussed in \Cref{sec:thermospring} and about the implementation of the associated reduced order models.

The system, as described in \cite[Section 4.3.1]{martin2016energy}, consists of two point masses connected by two
thermo-springs which exchange heat. The first spring connects the first
mass to a fixed point, while the second spring connects the two masses. In order to obtain a problem suitable for model order reduction, we retain the structure as described in \cite{martin2016energy} but increase the mechanical dimension.
The state vector is given by
\begin{equation*}
    \uf=\begin{bmatrix}
        \bb q_1^\top &\bb q_2^\top & \bb p_1^\top &\bb p_2^\top &S_1 &S_2
    \end{bmatrix}^\top
    \in\R^N,
\end{equation*}
where $N\coloneq4R+2$, $\bb q_1,\bb q_2\in\R^R$ are the position vectors, 
$\bb p_1,\bb p_2\in\R^R$ are the corresponding momenta, and the thermodynamic state is represented by the entropies $S_1,S_2$. 

Let us consider the following notation. Let $\widetilde{\uf}$ denote a generic state in $\R^N$ and let us use the symbol tilde to refer to its component. Moreover, for a generic element $\bw$ in $\R^n$, we denote with $\wrec(\mu)\coloneq\uf_0(\mu)+U\bw$ its lifting to the full order space shifted by the initial condition $\uf_0(\mu)$ of the full order model \eqref{eq:pbm} with parameter $\mu\in\p$.
%
Let us introduce the sets
\begin{equation*}
\begin{aligned}
    D & \coloneq\{ \widetilde{\uf}\in\R^N\mid\ \widetilde{\bb{q}}_1\neq0,\,\widetilde{\bb{q}}_1\neq\widetilde{\bb{q}}_2\},\\[0.5ex]
    D_r(\mu)&\coloneq\{\bw\in\R^n\mid\ \wrec(\mu)=\uf_0(\mu)+U\bw\in D\}.
\end{aligned}
\end{equation*}
The dynamics of this problem can be written in metriplectic form \eqref{eq:pbm}-\eqref{eq:degCond} as
\begin{equation*}
    \dot{\uf}(t;\mu)=J\nabla H(\uf(t;\mu))+G(\uf(t;\mu))\nabla S(\uf(t;\mu)),\qquad\forall\, t\in\mathcal{T},
\end{equation*}
where
the entropy of the system is linear and given by $S=S_1+S_2$, while the Hamiltonian is the total energy
\begin{equation*}
    H(\widetilde{\uf})=\frac12\left(\|\widetilde{\bb{p}}_1\|_2^2+\|\widetilde{\bb{p}}_2\|_2^2\right)
    +e_1(\lambda_1(\widetilde{\uf}),\widetilde{S}_1)+e_2(\lambda_2(\widetilde{\uf}),\widetilde{S}_2),\qquad\forall\,\widetilde{\uf}\in D,
\end{equation*}
where 
$\lambda_1(\widetilde{\uf}):=\norm{\widetilde{\bb{q}}_1}_2$, and $\lambda_2(\widetilde{\uf}):=\norm{\widetilde{\bb{q}}_2-\widetilde{\bb{q}}_1}_2$ are the elongations of the springs, and, for $j\in\{1,2\}$,
\begin{align*}
    e_j(\lambda_j(\widetilde{\uf}),\widetilde{S}_j) &=\frac12\log^2\left(\frac{\lambda_j(\widetilde{\uf})}{\lambda_j^0}\right)+\alpha_jT_{\rm ref}\log\left(\frac{\lambda_j(\widetilde{\uf})}{\lambda_j^0}\right)\\
    &\quad +c_jT_{\rm ref}\left[ \exp\left( \frac{\widetilde{S}_j-\alpha_j \log(\lambda_j(\widetilde{\uf})/\lambda_j^0)}{c_j}  \right)-1\right].
\end{align*}
Here, for $j\in\{1,2\}$, the constant $\lambda_j^0$ is the natural length of the $j$th spring, $c_j$ is its heat capacity, $\alpha_j$ accounts for the thermo-mechanical coupling,  and $T_\text{ref}$ is the reference temperature at which the unstressed deformation is defined to be zero.

Moreover, the Poisson and semi-metric operators have the following expressions
\begin{align*}
    J&=\begin{bmatrix}
        0_{2R\times 2R} & I_{2R} &0_{2R\times 1} & 0_{2R\times 1}\\
        -I_{2R} & 0_{2R\times 2R} & 0_{2R\times 1} & 0_{2R\times 1}\\
        0_{1\times 2R} & 0_{1\times 2R} & 0 & 0\\
        0_{1\times 2R} & 0_{1\times 2R} & 0 & 0
    \end{bmatrix},\\
    G(\widetilde{\uf})&=k \begin{bmatrix}
        0_{2R\times2R} & 0_{2R\times2R} & 0_{2R\times1} & 0_{2R\times1}\\
        0_{2R\times2R} & 0_{2R\times2R} & 0_{2R\times1} & 0_{2R\times1}\\
        0_{1\times 2R} & 0_{1\times2R} & T_2(\widetilde{\uf})/T_1(\widetilde{\uf}) & -1\\
        0_{1\times 2R} & 0_{1\times2R} & -1 & T_1(\widetilde{\uf})/T_2(\widetilde{\uf})
    \end{bmatrix},
\end{align*}
where $k\geq0$ is the coefficient of thermal conductivity, and, for any $\widetilde{\uf}\in D$,
\begin{equation*}
    T_j(\widetilde{\uf}):=\frac{\partial e_j}{\partial \widetilde{S}_j}(\lambda_j(\widetilde{\uf}), \widetilde{S}_j)= 
    T_{\rm ref}\exp\left(\frac{\widetilde{S}_j-\alpha_j\log(\lambda_j(\widetilde{\uf})/\lambda_j^0)}{c_j}\right),
    \qquad j\in\{1,2\}.
\end{equation*}

The gradient of the Hamiltonian and of the entropy are given, for any $\widetilde{\uf}\in D$, by
    \begin{equation*}
        \nabla H(\widetilde{\uf})=\begin{bmatrix}
            \dfrac{1}{\lVert\widetilde{\bb{q}}_1\rVert_2}\dfrac{\partial e_1}{\partial \lambda_1}(\widetilde{\uf})\widetilde{\bb{q}}_1 + \dfrac{1}{\lVert\widetilde{\bb{q}}_1 - \widetilde{\bb{q}}_2\rVert_2}\dfrac{\partial e_2}{\partial \lambda_2}(\widetilde{\uf})(\widetilde{\bb{q}}_1 -\widetilde{\bb{q}}_2)\\
            \dfrac{1}{\lVert\widetilde{\bb{q}}_1-\widetilde{\bb{q}}_2\rVert_2}\dfrac{\partial e_2}{\partial \lambda_2}(\widetilde{\bb{q}}_2 -\widetilde{\bb{q}}_1)\\
            \widetilde{\bb{p}}_1\\
            \widetilde{\bb{p}}_2\\
            \dfrac{\partial e_1}{\partial \widetilde{S}_1}(\widetilde{\uf})=T_1(\widetilde{\uf})\\
            \dfrac{\partial e_2}{\partial \widetilde{S}_2}(\widetilde{\uf})=T_2(\widetilde{\uf})
        \end{bmatrix},\quad
        \nabla S(\widetilde{\uf})=\begin{bmatrix}
            0_{R\times 1}\\
            0_{R\times 1}\\
            0_{R\times 1}\\
            0_{R\times 1}\\
            1\\
            1
        \end{bmatrix}.
    \end{equation*}
It can be easily verified that
the degeneracy conditions \eqref{eq:degCond} are satisfied,
and therefore the total energy is conserved, and the total entropy is
non-decreasing,
since $T_j(\widetilde{\uf})>0$, $j\in\{1,2\}$, for all $\widetilde{\uf}\in D$.
In particular, if $t\mapsto\uf(t;\mu)$ is solution of the full order model for the parameter $\mu\in\p$, then
\begin{equation*}
    \dot S(\uf(t;\mu)) = k\frac{(T_2(\uf(t;\mu))-T_1(\uf(t;\mu)))^2}{T_1(\uf(t;\mu))T_2(\uf(t;\mu))} \geq 0.
\end{equation*}
In the numerical experiments of \Cref{sec:thermospring}, we set $\lambda_j^0=c_j=\alpha_j=T_{\mathrm{ref}}=k=1$, for $j\in\{1,2\}$, obtaining the simplified expressions reported in the main text.
\paragraph{The block diagonal reduced order model (BD-ROM)}
The two thermo-spring system has the structure discussed in \Cref{sec:block}, that is it can be written as
\begin{equation*}
\dot{\uf}(t,\mu)=\begin{bmatrix}
    \dot{\uf}_1(t,\mu)\\[.5ex]
    \dot{\uf}_2(t,\mu)
\end{bmatrix}=
\begin{bmatrix}
    J_1\nabla_1 H(\uf(t,\mu))\\[.5ex]
    G_2(\uf(t,\mu))\nabla_2 S
\end{bmatrix}
\end{equation*}
with $J$ and $G$ of the form \eqref{eq:JGstruct} with, for any $\widetilde{\uf}\in D$,
\begin{equation}\label{eq:JGspring}
    J_1=\begin{bmatrix}
        0_{2R\times 2R} & I_{2R} \\
        -I_{2R} & 0_{2R\times 2R} 
    \end{bmatrix}, \
    G_2(\widetilde{\uf})=k\begin{bmatrix}
            T_2(\widetilde{\uf})/T_1(\widetilde{\uf}) & -1\\
            -1 & T_1(\widetilde{\uf})/T_2(\widetilde{\uf})
        \end{bmatrix},\ G_1=0_{4R\times 4R},
\end{equation}
so that $N_1=4R$ and $N_2=2$.

For this model, it is meaningful to consider a reduced basis of the form \eqref{eq:Ublock} with $U_2=I_2$, $n_2=N_2=2$, since the last two variables correspond to the entropy of the system. Consequently, $S_r=S$. The block $U_1\in\R^{N_1\times n_1}$ can be obtained by applying POD to the first $4R$ components of the snapshot matrix.

Concerning the reduced metriplectic conditions \eqref{eq:rdegCond},
the first one is trivially satisfied, owing to the block structure of $U$.
For the second condition \eqref{eq:GrdegCond}, let $\wrec=\uf_0(\mu)+U\bw$, with $\bw\in\R^n$. Using the fact that $G_1=0_{4R\times 4R}$ and $\Pi_{U_2}=I_2$, results in
\begin{align*}
    G_r(\bw)\nabla H_r(\bw;\mu)&=
    \begin{bmatrix}
        0_{n_1\times 1}\\
        G_2(\wrec)\nabla_{2}H(\wrec)
    \end{bmatrix}\\
    &=
    k\begin{bmatrix}
    0_{n_1\times 1}\\
        \frac{T_2(\wrec)}{T_1(\wrec)}T_1(\wrec)-T_2(\wrec)\\
        -T_1(\wrec)+\frac{T_1(\wrec)}{T_2(\wrec)}T_2(\wrec)
    \end{bmatrix}=0_{n\times1}.
\end{align*}

\paragraph{The oblique-projection reduced order model (OP-ROM)}
To derive the OP-ROM in \eqref{eq:OP-ROM}, we set
$\ccH=\ccS=2^{-1/2}\nabla S$. Clearly,
$$\bS=\sqrt{2},\qquad
\bH(\widetilde{\uf})=\dfrac{1}{\sqrt{2}}(T_1(\widetilde{\uf})+T_2(\widetilde{\uf}))>0,\qquad\forall\,\widetilde{\uf}\in D.$$
This implies that $\DmuS=\R^N$ and $\DmuH=D$, and, for any $\wrec=\uf_0(\mu)+U\bw$,
$$\QS=I_n-\dfrac12\nabla S_r\nabla S_r^\top,\qquad
\QH(\wrec)=I_n-\dfrac{1}{T_1(\wrec)+T_2(\wrec)}\nabla S_r\nabla H_r(\bw)^\top.$$

The metriplectic reduced order model reads: given $\mu\in\p$ and $\ur(0;\mu)=0$, find $\ur(\cdot;\mu)\in C^1(\mathcal{T};\R^n)$ such that
\begin{equation}\label{eq:OP-ROMstring}
    \dur(t;\mu)=\Jop\nabla H_r(\ur(t;\mu);\mu)+\Gop(U\ur(t;\mu)+\uf_0(\mu))\nabla S_r,\qquad \forall\,t\in\mathcal{T},
\end{equation}
with
\begin{equation*}
    \Jop=\QS U^\top J U\QS^\top,\quad \Gop(\wrec)= \QH(\wrec)U^\top G(\wrec) U\QH(\wrec)^\top,
\end{equation*}
for all $\wrec=\uf_0(\mu)+U\bw\in D$.

For this model, a reduced space of dimension $n$ satisfying the assumptions of \Cref{lem:OPmetriplectic} can be constructed via POD of the snapshot matrix truncated to $n-1$, then augmented with the vector $\ccH=\ccS\in\R^N$ and finally orthonormalized.
The resulting reduced order model is thus metriplectic.

\paragraph{Efficient evaluation of the reduced order model}
We show that, for this test case, the reduced order model OP-ROM in \eqref{eq:OP-ROMstring} can be solved at a cost independent of the dimension of the full order model, despite the nonlinearity in $\nabla H_r$ and in $\Gop$.

Let us introduce the following notation for the reduced basis $U\in\R^{N\times n}$,
$$U=\begin{bmatrix}
        U_{q_1}^\top & U_{q_2}^\top & U_{p_1}^\top & U_{p_2}^\top & U_{S_1}^\top & U_{S_2}^\top
    \end{bmatrix}^\top$$
where $U_{q_1},U_{q_2},U_{p_1},U_{p_2}\in\R^{R\times n}$, and $U_{S_1},U_{S_2}\in\R^{1\times n}$.

By definition of $H$, the gradient of the reduced Hamiltonian reads
\begin{align*}
\nabla H_r(\bw) &= U^\top \nabla H(\wrec(\mu))
= U_{q_1}^\top \bb{q}_1^{(r)}(\mu)\sum_{j=1,2}\lambda_j^{-1}(\wrec(\mu))\,\frac{\partial e_j}{\partial \lambda_j}(\wrec(\mu))\\
&\quad - \lambda_2^{-1}(\wrec(\mu))\,\frac{\partial e_2}{\partial \lambda_2}(\wrec(\mu))\,\left(U_{q_1}^\top\bb{q}_2^{(r)}(\mu) + U_{q_2}^\top \big( \bb{q}_1^{(r)}(\mu)-\bb{q}_2^{(r)}(\mu)\big)\right)\\
&\quad + U_{p_1}^\top \bb{p}_1^{(r)}(\mu)+ U_{p_2}^\top\bb{p}_2^{(r)}(\mu)
+ T_1(\wrec(\mu))U_{S_1}^\top
+ T_2(\wrec(\mu))U_{S_2}^\top,
\end{align*}
where the superscript $r$ refers to the components of the vector $\wrec(\mu)=\uf_0(\mu)+U\bw$.
For $i,j\in\{1,2\}$, it trivially holds
$$U_{q_i}^\top \bb{q}_j^{(r)}(\mu)=U_{q_i}^\top \big(\bb{q}_{0,j}(\mu)+U_{q_j}\bb{w}\big),$$
where the matrix product $U_{q_i}^\top U_{q_j}\in\R^{n\times n}$ can be computed offline, and the operations involving the initial condition can be performed once per instance of the parameter, so that a total of 
$O(n^2)$ online operations are required per instance of the reduced state $\bw\in\R^n$.
Similar considerations can be made for the remaining terms, since
\begin{equation*}
    \lambda_1(\wrec(\mu))^2=\lVert \bb{q}_{0,1}(\mu) + U_{q_1}\bw\rVert_2^2=\norm{\bb{q}_{0,1}(\mu)}_2^2 + 2\bb{q}_{0,1}(\mu)^\top U_{q_1}\bw  +  \norm{U_{q_1}\bw}_2^2,
\end{equation*}
and similarly for $\lambda_2(\wrec(\mu))^2$, and the fact that the term $T_j$, $j\in\{1,2\}$, is a scalar function of $\lambda_j$ and $S_j$.

In order to cheaply evaluate the semi-metric operator 
$\Gop$, one can exploit the block structure \eqref{eq:JGstruct} with blocks given in \eqref{eq:JGspring}.
Denoting by $U_S\in\R^{2\times n}$ the part of the reduced basis $U_S\coloneq\begin{bmatrix}
U_{S_1}^\top & U_{S_2}^\top
\end{bmatrix}^\top\in\R^{2\times n}$, it holds
    \begin{equation*}
        U^\top G(\uf_0(\mu)+U\bw)U=U_S^\top G_2(\uf_0(\mu)+U\bw)U_S,\qquad\forall\,\bw\in D_r(\mu).
    \end{equation*}
The evaluation at the reconstructed state entails
the computation of $T_1$ and $T_2$ at the lifted state variable, which can be done in $O(n^2)$ operations per state instance; 
the matrix products 
can be performed in $O(n^2)$ operations per state instance.

    Note that, since in this case the initial condition is parameter-dependent, an additional $O(nN)$ operations are needed once per parameter instance in order to precompute the quantities that depend on the initial condition. Hence evaluating the reduced order model requires $O(n^2)$ online operations per instance of the reduced state and $O(nN)$ operations per instance of the parameter $\mu$.

\medskip
The previous arguments also show that the EC-ROM of \cite{gruber2023energetically} can be solved with arithmetic complexity independent of the full dimension $N$. Indeed, from Theorem 3.4 of \cite{gruber2023energetically}, the reduced conservative part requires only the evaluation of the reduced Hamiltonian gradient, which can be computed in $O(n^2)$ operations per instance of the reduced state.
Moreover,  since $G(\widetilde{\uf})$ has rank one for any $\widetilde{\uf}\in D$, the entropy term of the EC-ROM can be written as
    \begin{equation*}
        G_{\text{EC}}(\bw;\mu)\nabla S_r:=k\left(\nabla S_r^\top \bb m(\bb{w})\right)\bb m(\bb{w}),\qquad\forall\,\bb{w}\in D_r(\mu),
    \end{equation*}
    where
    \begin{equation*}
        \bb m(\bb{w})\coloneq \frac{\left( U_{S_2}\nabla H_r(\bb{w})  \right)U_{S_1}^\top - \left( U_{S_1}\nabla H_r(\bb{w}) \right)U_{S_2}^\top}{\sqrt{T_1(\uf_0(\mu)+U\bb{w})T_2(\uf_0(\mu)+U\bb{w})}},\qquad\forall\,\bb{w}\in D_r(\mu).
    \end{equation*}
    Since computing $\bb m(\bb{w})\in\R^n$ only requires the temperatures $T_1(\uf_0(\mu)+U\bb{w}),T_2(\uf_0(\mu)+U\bb{w})$ and the Hamiltonian gradient $\nabla H_r(\bb{w}),$ the entropy term can be assembled in $O(n^2)$ online operations. As before, an additional $O(nN)$ cost is incurred once per parameter instance.

\section{Example: the damped thermoelastic rod}
\label{sec:rod-details}

This supplementary material contains further details on the damped thermoelastic rod system discussed in \Cref{sec:rod} and the implementation of the associated (hyper-)reduced order models.

The problem models a 1D elastic rod with space coordinate $x\in[0,L]$ which evolves as a damped Hamiltonian system with friction. In this case, the full state vector is given by
$\uf = \begin{bmatrix}
    \bb q^\top & \bb p^\top & S
\end{bmatrix}^\top\in\R^N$ with $\bb{q},\bb{p}\in\R^R$, $S\in\R$, and $N=2R+1$.
The system can be expressed in metriplectic formulation \eqref{eq:pbm}-\eqref{eq:degCond} as
\begin{equation*}
    \dot{\uf}(t;\mu)= J\nabla H(\uf(t;\mu))+G(\uf(t;\mu))\nabla S,\quad \forall\,t\in\mathcal{T},
\end{equation*}
where, for any $\widetilde{\uf} = \begin{bmatrix}
    \widetilde{\bb q}^\top & \widetilde{\bb p}^\top & \widetilde{S}
\end{bmatrix} ^\top\in\R^{N}$,
\begin{equation}\label{eq:JG}
    J=\begin{bmatrix}
        0_{R\times R} & I_R &0_{R\times 1}\\
        -I_R & 0_{R\times R} & 0_{R\times 1}\\
        0_{1\times R} & 0_{1\times R} &0
    \end{bmatrix},\qquad G(\widetilde{\uf})=\gamma\begin{bmatrix}
        0_{R\times R} & 0_{R\times R} & 0_{R\times 1}\\
        0_{R\times R} & I_R & -\widetilde{\bb{p}}\\
        0_{1\times R} & -\widetilde{\bb{p}}^\top & \norm{\widetilde{\bb{p}}}_2^2
    \end{bmatrix},
\end{equation}
with $\gamma\geq 0$, and 
\begin{equation*}
    H(\widetilde{\uf})=H(\widetilde{\bb{q}}, \widetilde{\bb{p}}, \widetilde{S})=\frac{1}{2}\lVert \widetilde{\bb{p}} \rVert_2^2 + \widetilde{S} + V(\widetilde{\bb{q}}),\qquad V(\widetilde{\bb{q}})=\sum_{i=1}^R\cos(\widetilde{q}_i).
\end{equation*}
In this example, the mass matrix is the identity $M=I_N$. A straightforward computation shows that $J\nabla S=0$ and $G(\widetilde{\uf})\nabla H(\widetilde{\uf})=0$ for all $\widetilde{\uf}\in\R^{N}$, and thus the system is indeed metriplectic.

\paragraph{The oblique-projection reduced order model}
To construct the reduced model \eqref{eq:OP-ROM}, we set $\ccH=\ccS=\bb e_{N}\in\R^N$.
Since the last component of $\nabla X$, for $X\in\{H,S\}$, is equal to one, this choice yields $\bX=1$ in \eqref{eq:beta}, so that
$$\QS=I_n-U^\top\bb{e}_N\nabla S^\top U,\qquad
\QH(U\bw+\uf_0(\mu))=I_n-U^\top\bb{e}_N\nabla H_r(\bw)^\top ,\quad\forall\,\bw\in\R^n.
$$
Similarly to the previous example, the evaluation of the reduced semi-metric operator $\Gr(\bw)$ can be carried out in $O(n^2)$ operations per instance of the reduced state $\bw$, plus an additional $O(nN)$ operations per parameter instance.

\paragraph{Hyper-reduction of the OP-ROM}

We split the Hamiltonian into a quadratic part $\LH$ and a nonlinear part $\NH$ defined as
$$\LH(\widetilde{\uf}):=\dfrac12\norm{\widetilde{\bb{p}}}_2^2+\widetilde{S},\qquad \NH(\widetilde{\uf})=V(\widetilde{\bb{q}}),\qquad \forall\,\widetilde{\uf}\in\R^N.$$
In this case, the nonlinear term $\NH$ admits the decomposition
$$\NH(\widetilde{\uf})=\cNh^\top \bb{h}(\widetilde{\bb q}),\qquad \forall\,\widetilde{\uf}\in\R^N,$$
where $\cNh=\bb{1}\in\mathbb{R}^R$ is the vector with all entries equal to $1$ and $h_i(\widetilde{\bb q}):=\cos(\widetilde{q}_i)$, for $1\leq i\leq R$.
Note that here $d=R$.

Let $U\in\R^{N\times n}$ be such that $U^\top U=I_n$, and let $$\uf_0(\mu) = \begin{bmatrix}
    \bb q_0(\mu)^\top & \bb p_0(\mu)^\top & S_0(\mu)
\end{bmatrix}^\top\in\R^N$$ be the initial condition of the full order model.  We consider the following representation
\begin{equation}\label{eq:U}
    U=\begin{bmatrix}
        U_q\\
        U_p\\
        u_S^\top
    \end{bmatrix},\qquad \mbox{with }\quad
    U_q\in\R^{R\times n},\,
    U_p\in\R^{R\times n},\,
    u_S\in\R^{n}.
\end{equation}
If the reduced basis $U\in\R^{N\times n}$ is split as in \eqref{eq:U}, the hyper-reduced Hamiltonian is given by
\begin{nalign}\label{eq:Hhr-rod}
    H_{\hr}(\bb{w};\mu)&= \LH(\uf_0(\mu)+U\bb{w})+ \NHhr(\bb{w};\mu)\\
    &= \dfrac12\norm{\bb p_0(\mu)+U_p\bb{w}}_2^2+ S_0(\mu)+u_S^\top \bb{w} 
    +\bb{1}^\top\pb_H \bb{h}(\bb q_0(\mu)+U_q\bb{w}),\quad\forall\bb{w}\in\R^n,
\end{nalign}
where $\pb_H\in\R^{R\times R}$ is a suitable EIM projection built as discussed in \Cref{sec:HSeim} and $\Heim(\widetilde{\uf})=\LH(\widetilde{\uf})+\bb{1}^\top\pb_H \bb{h}(\widetilde{\bb{q}})$, for any $\widetilde{\uf}\in\R^N$.

The hyper-reduced nonlinear gradient term approximation is thus given by
\begin{align*}
    \nabla\NHhr(\bb w;\mu)&=U_q^\top \Dh(\bb q_0 (\mu)+U_q\bb w)^\top \mathbb{P}_H^\top \bb 1\\
    &=-\underbrace{U_q^\top \operatorname{diag}(\sin(\bb q_0 (\mu)+U_q\bb w))P_H}_{I_1}\underbrace{(P_H^\top \Xi_H)^{-\top}\Xi_H^\top\bb 1}_{I_2},\qquad\forall\,\bb w\in\R^n,
\end{align*}
where $\Dh$ is the Jacobian matrix of $\bb h$. The term $I_2$ does not depend on the state, and can thus be computed offline; on the other hand, the term $I_1$ can be efficiently evaluated online, since
\begin{equation*}
    I_1 = \begin{bmatrix}
        \sin(\bb{q}_{0,i_1}(\mu)+U_{q_{i_1}}\bb w)U_{q_{i_1}}^\top & \cdots & \sin(\bb{q}_{0,i_{\nHeim}}(\mu)+U_{q_{i_{\nHeim}}}\bb w) U_{q_{i_{\nHeim}}}^\top
    \end{bmatrix}
\end{equation*}
where $i_1,\ldots, i_{\nHeim}$ are the EIM interpolation indices, $\bb{q}_{0,j}$ is the $j$th component of $\bb q_0 $ and $U_{q_j}\in\R^{1\times n}$ denotes the $j$th row of $U_q$. Evaluating the $\nHeim$ terms $U_{q_{i_1}}\bb w,\ldots,U_{q_{i_{\nHeim}}}\bb w$ requires $\nHeim$ scalar products; thus the cost of evaluating the $I_1$ term is $O(n\nHeim)$, independent of $N$. The hyper-reduced Hamiltonian gradient can thus be expressed, for all $\bw\in\R^n$, as
\begin{align*}
    \nabla  H_{\hr}(\bb w;\mu)&= -U_q^\top \operatorname{diag}(\sin(\bb q_0(\mu)+U_q\bb w))P_H(P_H^\top \Xi_H)^{-\top}\Xi_H^\top\bb 1\\
    &\quad + U_p^\top (\bb p_0(\mu)+U_p\bb w)+u_S,
\end{align*}
which can be evaluated with an online cost of $O(n\nHeim + n^2)$ algebraic operations per instance of the reduced state. Again, since the initial condition is parameter dependent, the precomputable terms involving the initial condition need to be recomputed for each parameter instance, but only once per time simulation. This amounts to an extra $O(nN)$ cost once per parameter instance.
Note that this complexity analysis is consistent with the findings from \Cref{oss:GPeim}. Indeed, in this case, every entry of $\bb h(\widetilde{\bb q})$ depends on a single entry of $\widetilde{\bb q}$, and hence $s_{\bb{h}}=1$, giving an online cost for the evaluation of the hyper-reduced gradient of the nonlinear term of $O(n\nHeim s_{\bb{h}})=O(n\nHeim)$, per state instance.

Regarding the semi-metric operator, one could derive $\barGhr$ via EIM approximation of $\vect(G)$ where $G$ is defined in \eqref{eq:JG}.
However, it is possible to assemble and evaluate the reduced semi-metric operator $G_r$ with arithmetic complexity independent of $N$, as follows.
From the definition of $G$ one can easily verify that
\begin{equation}\label{eq:Gbar}
    \barG(\widetilde{\uf})=U^\top G(\widetilde{\uf})U=\gamma\left(U_p^\top U_p - U_p^\top \widetilde{\bb p} u_S^\top - u_S \widetilde{\bb p}^\top U_p + \norm{\widetilde{\bb p}}_2^2 u_S u_S^\top\right),\quad\forall\,\widetilde{\uf}\in\R^{N},
\end{equation}
Then, for any $\bb{w}\in\R^n$,
\begin{align*}
    \barG(\uf_0(\mu)+U\bb{w})&=U^\top G(\uf_0(\mu)+U\bb{w})U\\
    &=\gamma \Big(U_p^\top U_p - U_p^\top \bb{p}_0(\mu) u_S^\top
    -U_p^\top U_p\bb{w} u_S^\top- u_S \bb{p}_{0}(\mu)^\top U_p\\
    &\quad -u_S \bb{w}^\top U_p^\top U_p+ \left( \norm{\bb{p}_0(\mu)}_2^2 + 2\bb{p}_0(\mu)^\top U_p\bb{w} + \bb{w}^\top U_p^\top U_p \bb{w} \right)u_S u_S^\top\Big),
\end{align*}
which can be evaluated with $O(nN)+O(n^2)$ operations per instance of the parameter $\mu$, and $O(n^2)$ operations per instance of the state $\bb{w}\in\R^n$.

The hyper-reduced model \eqref{eq:hrOP-ROM} reads:
given $\mu\in\p$ and $\uhr(0;\mu)=0\in\R^n$, find $\uhr(\cdot;\mu)\in C^1(\mathcal{T};\R^n)$ such that
\begin{equation*}
    \duhr(t;\mu)=\Jop\nabla H_{\hr}(\uhr(t;\mu);\mu)+\Gop(U\uhr(t;\mu)+\uf_0(\mu);\mu)U^\top \nabla S,\qquad\forall\,t\in\mathcal{T}.
\end{equation*}
The function $H_{\hr}$ is given in \eqref{eq:Hhr-rod}, while the hyper-reduced Poisson and semi-metric tensors are defined, for any $\widetilde{\uf}\in \R^N$, as
\begin{equation*}
\begin{aligned}
    \Jop(\widetilde{\uf};\mu)&\coloneq \QS(\widetilde{\uf};\mu)\barJ(\widetilde{\uf})\QS(\widetilde{\uf};\mu)^\top,\\
    \Gop(\widetilde{\uf};\mu)&\coloneq \QHhr (\widetilde{\uf};\mu)\barG(\widetilde{\uf})\QHhr(\widetilde{\uf};\mu)^\top,
\end{aligned}
\end{equation*}
where the operators $\QS$ and $\QHhr$ are defined as in \eqref{eq:PinQ}. Since the last component of both $\nabla \Heim$ and $\nabla S$ is equal to one at all times, setting $\cHhr=\cS$ equal to the $N$th vector $\bb{e}_{N}$ of the canonical basis in $\R^{N}$, gives 
$$\bHhr(\widetilde{\uf};\mu)=\bS(\widetilde{\uf};\mu)=1\qquad\forall\,\widetilde{\uf}\in\R^N,\,\mu\in\p.$$
In the following we assume that the reduced basis matrix $U\in\R^{N\times n}$ employed in the OP-ROM and OP-HROM is such that $\bb e_N\in\operatorname{Range}(U)$. Observe that
$\barJ=U^\top J U=U_q^\top U_p-U_p^\top U_q\in\R^{n\times n}$ and, since $\bS(\cdot;\mu)\equiv 1$, then
$$\PS(\widetilde{\uf};\mu) := U^{\top} M\cS \nabla S(\widetilde{\uf};\mu)^\top U= U^{\top} \bb{e}_{N} \bb{e}_{N}^\top U = u_S u_S^\top,\qquad\forall\,\widetilde{\uf}\in\R^N.$$
Moreover, since $U U^\top \bb{e}_{N}=Uu_S=\bb{e}_{N}$, it holds $u_S^\top u_S = 1$ and $U_qu_S=U_p u_S=0$.
This implies that $\PS(\cdot;\mu)\barJ=0$, and, hence,
\begin{equation*}
    \Jop(\widetilde{\uf};\mu)\coloneq(I_n-\PS(\widetilde{\uf};\mu))\barJ(\widetilde{\uf})(I_n-\PS(\widetilde{\uf};\mu))^\top=
    \barJ,\qquad\forall\,\widetilde{\uf}\in\R^N.
\end{equation*}
Moreover, since $\bHhr\equiv 1$, we have
$$\PHhr(\widetilde{\uf};\mu) := U^{\top} M\cHhr \nabla \Heim(\widetilde{\uf})^\top U
= u_S \nabla \Heim(\widetilde{\uf})^\top U,\quad\forall\,\widetilde{\uf}\in\R^N$$
so that $\Gop$ reads, for any $\widetilde{\uf}\in\R^N$,
\begin{equation}\label{eq:Gop-rod}
    \Gop(\widetilde{\uf};\mu)=\Gop(\widetilde{\uf})=\big(I_n- u_S \nabla \Heim(\widetilde{\uf})^\top U\big)U^\top G(\widetilde{\uf})U\big(I_n-U^\top\nabla \Heim(\widetilde{\uf})u_S^\top\big).
\end{equation}

\paragraph{Hyper-reduction of the EC-ROM of \cite{gruber2023energetically}}
For the numerical experiments of \Cref{table:thermorod_hrsummary},
we consider, for the sake of comparison, the EC-ROM of \cite{gruber2023energetically} with a hyper-reduction of the nonlinear terms analogous to the one discussed in the previous paragraph.
Indeed, from \cite[Theorem 3.4]{gruber2023energetically} it can be easily verified that the conservative part of the dynamics can be written in terms of $J_r$, $\nabla H_r$, and $\nabla S_r$, which, in this test case, can be either precomputed or hyper-reduced through the gradient-preserving EIM technique recalled in \Cref{sec:HSeim}. 
Next we show that the entropy term can be written in terms of $\barG$, $\nabla H_r$, and $\nabla S_r$ and, thus, the discussion from the previous paragraph holds in this case as well.

Using the notation from \cite{gruber2023energetically}, let us write the semi-metric operator $G$ as
    \begin{equation}
        G(\widetilde{\uf})=\sum_{\alpha=1}^{R} \malpha(\widetilde{\uf})\malpha(\widetilde{\uf})^\top,\qquad \malpha(\widetilde{\uf})\coloneq\sqrt\gamma\begin{bmatrix}
            0_{R\times 1}\\
            \bb e_\alpha\\
            -\widetilde{p}_\alpha
        \end{bmatrix},\quad\forall\,\widetilde{\uf}\in\R^{N}.
    \end{equation}
If we denote with $G_{\text{EC}}$ the semi-metric operator in the EC-ROM, by \cite[Theorem 3.4]{gruber2023energetically}, one has, for all $\bw\in\R^n$,
\begin{equation*}
    G_{\text{EC}}(\bw;\mu)\nabla S_r =\sum_{\alpha=1}^R\left(\nabla S_r^\top\Ahat(\uf_0(\mu)+U\bw)\nabla H_r(\bw)\right) \Ahat(\uf_0(\mu)+U\bw)\nabla H_r(\bw),
\end{equation*}
where
\begin{equation*}
    \Ahat(\widetilde{\uf})\coloneq\ahat(\widetilde{\uf}) u_S^\top- u_S\ahat(\widetilde{\uf})^\top,\qquad
    \ahat (\widetilde{\uf})\coloneq U^\top \malpha(\widetilde{\uf}), \qquad\forall\,\widetilde{\uf}\in\R^{N};
\end{equation*}
Let $B(\widetilde{\uf})\coloneq \sum_{\alpha=1}^R\ahat(\widetilde{\uf})\ahat(\widetilde{\uf})^\top$, for any $\widetilde{\uf}\in\R^{N}$, and
$c(\bw;\mu)\coloneq u_S^\top\nabla  H_r(\bw;\mu)$, for any $\bw\in\R^n$.
Recalling that $\nabla S_r=u_S$, a straightforward computation yields
\begin{align*}
    G_{\text{EC}}(\bw;\mu)\nabla S_r&= c(\bw;\mu)^2 B(\uf_0(\mu)+U\bw)u_S \\
    &\quad-c(\bw;\mu) \lVert u_S\rVert_2^2 B(\uf_0(\mu)+U\bw)\nabla H_r(\bw;\mu)\\
    &\quad - \left(c(\bw;\mu) u_S^\top
    -\lVert u_S\rVert_2^2  \nabla H_r(\bw;\mu)^\top\right) B(\uf_0(\mu)+U\bw)\nabla  H_r(\bw;\mu) u_S.
    \end{align*}
Moreover, the identity
    \begin{align*}
B(\widetilde{\uf})&=\sum_{\alpha=1}^R\ahat(\widetilde{\uf})\ahat(\widetilde{\uf})^\top
=\sum_{\alpha=1}^RU^\top \malpha(\widetilde{\uf})\malpha(\widetilde{\uf})^\top U
= U^\top G(\widetilde{\uf}) U,\qquad\forall\,\widetilde{\uf}\in\R^{N},
\end{align*}
implies  that $B(\uf_0(\mu)+U\bw)=G_r(\bw;\mu)$, for any $\bw\in\R^n$.
This means that the term 
$G_{\text{EC}}(\bw;\mu)\nabla S_r$
can be efficiently computed via hyper-reduction of $\nabla H_r$, as discussed in the previous paragraph.

\paragraph{Verification of the hypotheses of Theorem \ref{thm:convergence-hr}}
We verify that this example satisfies the hypotheses of the convergence \Cref{thm:convergence-hr} under the hypothesis that, for any parameter $\mu\in\p$ the initial condition of the full order model satisfies the uniform bound
$$\norm{\bb{q}_0(\mu)}_2^2+\norm{\bb{p}_0(\mu)}_2^2\leq P_0^2,$$
for some $P_0\in\R$.
In order to simplify the notation, here we consider the unshifted reconstruction of the full order state, i.e., $\uf(t;\mu)$ is approximated by $U\uhr(t;\mu)$, so that the initial condition of \eqref{eq:hrOP-ROM} is $\uhr(0;\mu)=U^\top\uf_0(\mu)$.

First, we observe that the operator $\Xophr$ is defined in this example as
\begin{equation}\label{eq:Xophr-rod}
    \Xophr(\widetilde{\uf})=\barJ U^\top \nabla \Heim(\widetilde{\uf})+\Gop(\widetilde{\uf})U^\top\nabla S,\qquad \forall\,\widetilde{\uf}\in \R^N.
\end{equation}
Since the operator $\Xophr$ is of the form above, we need to verify only hypotheses \eqref{eq:opbound-hr}, \eqref{eq:boundinp-hr} and the first equation in \eqref{eq:fracbound-hr} in the form
\begin{equation*}
\norm{\bb{e}_N}_M \max\{\norm{M^{-1}\nabla \Heim(\uf(t;\mu))\|_M,\|M^{-1}\nabla S}_M\}\leq\kH,
\end{equation*}
with $\uf$ solution of the full order model.

To verify hypothesis \eqref{eq:opbound-hr}, observe that, for any $\widetilde{\uf} = (\widetilde{\bb q}^\top, \widetilde{\bb p}^\top, \widetilde{S})^\top\in\R^{N}$, it holds
$$\norm{M^{-1}\JM(\widetilde{\uf} )}_{M\to M}=\norm{\J}_2=1,
\qquad
\norm{M^{-1}\GM(\widetilde{\uf} )}_{M\to M}=\norm{\G(\widetilde{\uf} )}_2=\gamma(1+\norm{\widetilde{\bb p}}_2^2).
$$
To bound the $2$-norm of $G$, we take the mechanical energy
$$E(\widetilde{\uf}):=\dfrac12\norm{\widetilde{\bb p}}_2^2+V(\widetilde{\bb q})$$
and observe that it satisfies
$\dot E(\uf(t;\mu))=-\gamma\norm{\bb p(t;\mu)}_2^2\leq 0$ for $\uf$ solution of the full order model.
If the initial condition $\uf_0(\mu) = (\bb q_0^\top(\mu), \bb p_0^\top(\mu), S_0(\mu))\in\R^{N}$ satisfies
$\norm{\bb p_0(\mu)}_2\leq P_0$
uniformly in $\mu$ for some $P_0\in\R$, then
$$\dfrac12\norm{\bb p(t;\mu)}_2^2
\leq E(\bb q_0(\mu),\bb p_0(\mu))-V(\bb q(t;\mu))
\leq \dfrac12 P_0^2 + 2R,$$
since $-R\leq V(\bb z)\leq R$ for any $\bb z\in\R^R$.
Thus, $\aj=1$ and $\ag=\gamma(1+P_0^2+4R)$ in \eqref{eq:opbound-hr}.

Concerning the bound in \eqref{eq:boundinp-hr}, observe that
the definition of the hyper-reduced Hamiltonian, with $\Dh(U_q\bb{w})=-\operatorname{diag}(\sin(U_q\bb{w}))$, for any $\bb{w}\in\R^n$, gives
$$U^\top \nabla \Heim(\widetilde{\uf})=U_q^\top \Dh(\widetilde{\bb{q}})^\top\mathbb{P}_H ^\top\bb{1} +U_p^\top \widetilde{\bb{p}}+u_S,\qquad\forall\,\widetilde{\uf}\in\R^N,$$
meaning that the first term of $\Xophr$ in \eqref{eq:Xophr-rod} is equal to $\barJ f_{\hr}(\uf)$ with $f_{\hr}(\widetilde{\uf}):=U_q^\top \Dh(\widetilde{\bb q})^\top\mathbb{P}_H^\top\bb{1}+U_p^\top \widetilde{\bb{p}}$.
If we evaluate $f_{\hr}$ at the solution $\uf$ of the full order model and use the bound for the norm of $\bb{p}$, we get
\begin{equation}\label{eq:normfhr}
\norm{f_{\hr}(\uf)}^2_2\leq\norm{\operatorname{diag}(\sin(\bb{q}))}_2^2\norm{\mathbb{P}_H ^\top\bb{1}}^2_2+\norm{\bb p}^2_2\leq R\norm{\mathbb{P}_H}_2^2+\norm{\bb p}_2^2\leq R\norm{\mathbb{P}_H}_2^2+ P_0^2+4R.
\end{equation}

For the second term of $\Xophr(\cdot;\mu)$, we use the expression of $\Gop$ derived in \eqref{eq:Gop-rod}.
Since $U^\top \nabla S(\uf;\mu)=u_S$,
standard calculations show that the second term of $\Xophr$ reads
$$\Gop(\widetilde{\uf})u_S=\gamma\big(-U_p^\top U_p f_{\hr}(\widetilde{\uf})+u_S f_{\hr}(\widetilde{\uf})^\top U_p^\top U_p f_{\hr}(\widetilde{\uf})\big),\qquad\forall\widetilde{\uf}\in\R^N.$$
Putting all terms together yields, for all $\widetilde{\uf}\in\R^N$,
\begin{equation*}
    \begin{aligned}
    \Xophr(\widetilde{\uf};\mu)=\Xophr(\widetilde{\uf})&=\barJ U^\top\nabla \Heim(\widetilde{\uf})+\Gop(\widetilde{\uf})u_S\\
    & =(\barJ-\gamma U_p^\top U_p) f_{\hr}(\widetilde{\uf})+
    \gamma u_S f_{\hr}(\widetilde{\uf})^\top U_p^\top U_p f_{\hr}(\widetilde{\uf}).
    \end{aligned}
\end{equation*}
To find the bound in \eqref{eq:boundinp-hr}, we observe that, for any $\widehat{\uf},\widetilde{\uf}\in\R^{N}$,
\begin{equation*}
    \begin{aligned}
    \Xophr(\widehat{\uf};\mu)-\Xophr(\widetilde{\uf};\mu)=
    &\,(\barJ-\gamma U_p^\top U_p) (f_{\hr}(\widehat{\uf})-f_{\hr}(\widetilde{\uf}))\\
    &+ \gamma u_S f_{\hr}(\widehat{\uf})^\top U_p^\top U_p f_{\hr}(\widehat{\uf})
    -\gamma u_S f_{\hr}(\widetilde{\uf})^\top U_p^\top U_p f_{\hr}(\widetilde{\uf}).
    \end{aligned}
\end{equation*}
Since $U_p^\top U_p$ is symmetric we have
$$f_{\hr}(\widehat{\uf})^\top U_p^\top U_p f_{\hr}(\widehat{\uf})
- f_{\hr}(\widetilde{\uf})^\top U_p^\top U_p f_{\hr}(\widetilde{\uf})
= (f_{\hr}(\widehat{\uf})-f_{\hr}(\widetilde{\uf}))^\top U_p^\top U_p (f_{\hr}(\widehat{\uf})+f_{\hr}(\widetilde{\uf})).$$
Then, using the fact that
$$\norm{\barJ}_2=\norm{U_q^\top U_p-U_p^\top U_q}_2\leq 1,\qquad\norm{U_p^\top U_p}_2\leq 1,$$
gives, for any $\widehat{\uf},\widetilde{\uf}\in\R^{N}$,
\begin{equation*}
\begin{aligned}
\norm{\Xophr(\widehat{\uf};\mu)-\Xophr(\widetilde{\uf};\mu)}_2^2 & 
\leq \norm{f_{\hr}(\widehat{\uf})-f_{\hr}(\widetilde{\uf})}_2^2\big((1+\gamma)^2 \\
&\quad +2\gamma^2 \norm{u_S}_2^2 \norm{U_p^\top U_p}^2_2\big(\norm{f_{\hr}(\widehat{\uf})}_2^2+\norm{f_{\hr}(\widetilde{\uf})}_2^2\big)\big).
\end{aligned}
\end{equation*}
To bound the term $\norm{f_{\hr}(\widehat{\uf})}_2$ one can use  \eqref{eq:normfhr}, if $\widehat{\uf}=\uf$ is the solution of the full order model. When $\widehat{\uf}=U\uhr$, with $\uhr$ solution of the hyper-reduced model, one can only conclude that
$$\norm{f_{\hr}(U\uhr)}_2^2\leq R\norm{\mathbb{P}_H}_2^2 + \norm{U_p\uhr}_2^2.$$

Let $E_{\hr}$ denote the hyper-reduction of the mechanical energy $E$, that is
$$E_{\hr}(\bb{w}):=\dfrac12\norm{U_p\bb w}_2^2+ \NHhr(\bb{w})
    = \dfrac12\norm{U_p\bb{w}}_2^2
    +\bb{1}^\top\pb_H \bb{h}(U_q\bb{w}),\qquad\forall\,\bb{w}\in\R^n.$$
Then, since $u_S^\top f_{\hr}=0$ and $\barJ$ is skew-symmetric, it holds
\begin{equation*}
    \begin{aligned}
        \dfrac{d}{dt}E_{\hr}(\uhr) &=
        \nabla E_{\hr}(\uhr)^\top\duhr=
        f_{\hr}(U\uhr)^\top\Xophr(U\uhr)\\
        & =f_{\hr}(U\uhr)^\top\big(\barJ-\gamma U_p^\top U_p\big)f_{\hr}(U\uhr)\\
        &\quad+\gamma f_{\hr}(U\uhr)^\top u_S f_{\hr}(U\uhr)^\top U_p^\top U_p f_{\hr}(U\uhr)\\
        & =-\gamma \norm{U_p f_{\hr}(U\uhr)}_2^2\leq 0.
    \end{aligned}
\end{equation*}
This implies that $E_{\hr}(\uhr(t;\mu))\leq E_{\hr}(\uhr(0;\mu))$.
Moreover, since $-\norm{\mathbb{P}_H^\top \bb{1}}_1\leq \bb{1}^\top\pb_H \bb{h}(\bb{z})\leq \norm{\mathbb{P}_H^\top \bb{1}}_1$, for any $\bb{z}\in\R^R$, then
\begin{align*}
    \dfrac12 \norm{U_p\uhr(t;\mu)}_2^2&\leq 
\dfrac12 \norm{U_p\uhr(0;\mu)}_2^2+
\bb{1}^\top\pb_H \big(\bb{h}(U_q\uhr(0;\mu))
- \bb{h}(U_q\uhr(t;\mu))\big)\\
&\leq \dfrac12 P_0^2 +2\norm{\mathbb{P}_H^\top \bb{1}}_1,
\end{align*}
if we assume that, for any parameter $\mu\in\p$ for which the hyper-reduced model is tested,
$$\norm{U_p\uhr(0;\mu)}_2^2=\norm{U_p(U_q^\top \bb{q}_0(\mu)+U_p^\top\bb{p}_0(\mu))}_2^2\leq
\norm{\bb{q}_0(\mu)}_2^2+\norm{\bb{p}_0(\mu)}_2^2\leq P_0^2.$$
This implies that
$$\norm{f_{\hr}(U\uhr)}_2^2\leq R\norm{\mathbb{P}_H}_2^2 + P_0^2 + 4\norm{\mathbb{P}_H^\top \bb{1}}_1.$$
Finally, we observe that,
for $\widehat{\uf} = (\widehat{\bb q}^\top, \widehat{\bb p}^\top, \widehat{S})^\top,
\widetilde{\uf} = (\widetilde{\bb q}^\top, \widetilde{\bb p}^\top, \widetilde{S})^\top\in\R^{N}$, it holds
\begin{equation*}
\begin{aligned}
\norm{f_{\hr}(\widehat{\uf})-f_{\hr}(\widetilde{\uf})}_2^2 &=
\left\lVert U^\top
\begin{pmatrix}
    \big(\Dh(\widehat{\bb{q}})-\Dh(\widetilde{\bb{q}})\big)\mathbb{P}_H^\top\bb{1}\\
    \widehat{\bb{p}}-\widetilde{\bb{p}}\\
    0
\end{pmatrix}\right\rVert_2^2\\
&\leq \norm{\Dh(\widehat{\bb{q}})-\Dh(\widetilde{\bb{q}})}_F^2\norm{\mathbb{P}_H^\top\bb{1}}_2^2+\norm{\widehat{\uf}-\widetilde{\uf}}_2^2\\
&\leq (1+R\norm{\mathbb{P}_H}_2^2) \norm{\widehat{\uf}-\widetilde{\uf}}_2^2.
\end{aligned}
\end{equation*}
This implies that, for
$$C:=\sqrt{(1+\gamma)^2+2\gamma^2(2P_0^2+4R+2R\norm{\mathbb{P}_H}_2^2+4\norm{\mathbb{P}_H^\top\bb{1}}_1)}\sqrt{1+R\norm{\mathbb{P}_H}_2^2},$$
it holds
\begin{equation*}
\begin{aligned}
   \big\langle \bb{e}(t;\mu),\,
   \Xophr(\uf(t;\mu);\mu)&-\Xophr(U\uhr(t;\mu);\mu)\big\rangle_2
    \leq C \norm{\bb{e}(t;\mu)}_2 \norm{\uf-U\uhr}_2\\
   & \leq \norm{\bb{e}(t;\mu)}_2\big(C_e\,\norm{\bb{e}(t;\mu)}_2
   +C_u\norm{(I_N-\Pi_U)\uf(t;\mu)}_M\big),
\end{aligned}
\end{equation*}
which corresponds to the bound \eqref{eq:boundinp-hr} 
with $C_u=C_e=C$.

Condition \eqref{eq:fracbound-hr} requires us to find $\kH>0$ such that
\begin{equation*}
\norm{\bb{e}_N}_M \max\{\norm{M^{-1}\nabla \Heim(\uf(t;\mu))}_M,\norm{M^{-1}\nabla S(\uf(t;\mu))}_M\}\leq\kH.
\end{equation*}
Using the fact that $M$ is the identity matrix, and the definition of $S$, we get that $\norm{M^{-1}\nabla S(\uf(t;\mu))}_M=\norm{\nabla S}_2=1$, while
$$\norm{\nabla \Heim(\uf(t;\mu))}_2^2=\norm{\Dh(\bb{q}(t;\mu))^\top\mathbb{P}_H^\top\bb{1}}_2^2+\norm{\bb{p}(t;\mu)}_2^2+1
\leq R\norm{\mathbb{P}_H}_2^2 +P_0^2 + 4R +1,$$
with a reasoning similar to \eqref{eq:normfhr}.
Then, condition \eqref{eq:fracbound-hr} holds with $$\kH=\sqrt{R\norm{\mathbb{P}_H}_2^2 +P_0^2 + 4R +1}.$$

\end{document}